\documentclass[11pt]{article}
\usepackage{url,fullpage,amsmath,graphicx,amssymb,algorithmic,fancyvrb,xspace}
\usepackage[pdftex]{color}
\usepackage{pdfcolmk}  % <-- fixes bug of colors not spanning pagebreaks

\usepackage[pdftex]{hyperref}
\usepackage{algorithm}
\DefineVerbatimEnvironment{code}{Verbatim}{tabsize=0,xleftmargin=0.25in}

\usepackage[small,bf]{caption}
\usepackage{subfig}

\usepackage{cite}  % puts citation numbers in order, changes 14,15,16 to 14-16.

\newtheorem{theorem}{Theorem}[section]

\newcommand\cA{\ensuremath{\mathcal{A}}}
\newcommand\bb{\ensuremath{y}}   % specific b, e.g. data.  Using ''y'' for now.
\newcommand{\order}{\mathcal{O}}              % big O notation
\newcommand{\T}{*}                            % for the adjoint/transpose

\newcommand\cAA{A}
\newcommand\cWW{W}

\newcommand\dual{\lambda}

\newcommand\eg{e.g.\xspace}
\newcommand\ie{i.e.\xspace}

\newcommand\thalf{{\textstyle\frac{1}{2}}}

\newcommand\dom{\operatorname{\textrm{dom}}}

\DeclareMathOperator*{\argmax}{arg\,max}        % puts subscripts in the right place
\DeclareMathOperator*{\argmin}{arg\,min}
\newcommand{\bmat}[1]{\begin{bmatrix} #1 \end{bmatrix}}
\floatname{algorithm}{Listing}
\newcommand\cK{\mathcal{K}}

\newcommand\cP{\mathcal{P}}

\newcommand{\lag}{\ensuremath{\mathcal{L}}}   % Lagrangian

\newcommand{\gs}{g_\text{sm}}             % smooth part of dual objective

\newcommand{\R}{\mathbb{R}}

\newcommand{\<}{\langle}
\renewcommand{\>}{\rangle}

\newcommand{\goto}{\rightarrow}
\newcommand{\sgn}{\textrm{sgn}}

\newcommand{\SRB}[1]{\textcolor{blue}{SRB: #1}}
\newcommand{\MCG}[1]{\textcolor{magenta}{MCG: #1}}
\newcommand{\CHANGE}[1]{\textcolor{red}{CHANGE: #1}}
\renewcommand{\vec}[1]{{\boldsymbol{#1}}}

\newcommand{\rank}{\operatorname{rank}}

\newcommand{\supp}[1]{\operatorname{supp}(#1)}

\newcommand{\Id}{\text{\em I}}

\numberwithin{equation}{section}
\def \endprf{\hfill {\vrule height6pt width6pt depth0pt}\medskip}
\newenvironment{proof}{\noindent {\bf Proof} }{\endprf\par}

\newcommand{\iprod}[2]{\left\langle #1 , #2 \right\rangle}

\newcommand{\ST}{\operatorname{SoftThreshold}}
\newcommand{\SVT}{\operatorname{SoftThresholdSingVal}}
\newcommand{\Shrink}{\operatorname{Shrink}}
\newcommand{\Trunc}{\operatorname{Trunc}}
\newcommand{\CTrunc}{\operatorname{CTrunc}}

\title{Templates for Convex Cone Problems\\with Applications to Sparse
  Signal Recovery}

\author{Stephen Becker$^{1}$, Emmanuel J. Cand\`es$^{2}$ and Michael Grant$^{1}$\\
  \vspace{-.1cm}\\
  $^1$Applied and Computational Mathematics, Caltech, Pasadena, CA 91125\\
  \vspace{-.3cm}\\
  $^2$Departments of Mathematics and of Statistics, Stanford
  University, Stanford,
  CA 94305}

\date{September 2010}

\begin{document}
\maketitle

\begin{abstract}
  This paper develops a general framework for solving a variety of
  convex cone problems that frequently arise in signal processing,
  machine learning, statistics, and other fields. The approach works as
  follows: first, determine a conic formulation of the problem; 
  second, determine its dual; third, apply smoothing; and
  fourth, solve using an optimal first-order method. A merit of this
  approach is its flexibility: for example, all compressed
  sensing problems can be solved via this approach. These include
  models with objective functionals such as the total-variation norm, $\|Wx\|_1$
  where $W$ is arbitrary, or a combination thereof. In addition, the
  paper also introduces a number of technical contributions such as a
  novel continuation scheme, a novel approach for controlling the step
  size, and some new results showing that the smooth and
  unsmoothed problems are sometimes formally equivalent. Combined with
  our framework, these lead to novel, stable and computationally
  efficient algorithms. For instance, our general implementation is
  competitive with state-of-the-art methods for solving intensively
  studied problems such as the LASSO. Further, numerical experiments
  show that one can solve the Dantzig selector problem, for which no
  efficient large-scale solvers exist, in a few hundred iterations.
  Finally, the paper is
  accompanied with a software release.
  This software is not a single,
  monolithic solver; rather, it is a suite of programs and routines
  designed to serve as building blocks for constructing complete
  algorithms.
\iffalse
  \SRB{[Proposing to use this sentence instead:]
  The software is a suite of routines designed to serve as building
  blocks for constructing customized solvers; additionally,
  pre-constructed solvers to the LASSO and other common problems
  are provided.}
  \MCG{[I frankly see no problem with duplicating language between
  the abstract and main text. Furthermore my preference is not to mention 
  LASSO more than once lest we draw undue attention to it. Therefore
  I prefer the sentence as it is written; or perhaps 
  something like this:] This software is not a single,
  monolithic solver; rather, it is a suite of programs and routines
  designed to serve as building blocks for constructing complete
  algorithms, and includ}
\fi  
\end{abstract}

{\bf Keywords.} Optimal first-order methods, Nesterov's accelerated
descent algorithms, proximal algorithms, conic duality, smoothing by
conjugation, the Dantzig selector, the LASSO, nuclear-norm
minimization.

\section{Introduction}
\label{sec:intro}

\subsection{Motivation}
\label{sec:motivation}

This paper establishes a general framework for constructing optimal
first-order methods for solving certain types of convex optimization
programs that frequently arise in signal and image processing, statistics,
computer vision, and a variety of other fields.\footnote{The meaning of
  the word `optimal' shall be made precise later.} 
%To be as concrete
%as possible, suppose 
In particular,  % taking Ewout's suggesting, since we do use ``suppose'' a lot
we wish to recover an unknown vector $x_0 \in
\R^n$ from the data $y \in \R^m$ and the model
\begin{equation}
  \label{eq:model}
  y = A x _0 + z;
\end{equation}
here, $A$ is a known $m \times n$ design matrix and $z$ is a noise
term. To fix ideas, suppose we find ourselves in the increasingly
common situation where there are fewer observations/measurements than
unknowns, \ie, $m < n$. While this may seem {\em a priori} hopeless,
an impressive body of recent works has shown that accurate estimation
is often possible under reasonable sparsity constraints on $x_0$. % taking Ewout's suggestion
One
practically and theoretically effective estimator is the {\em Dantzig
  selector} introduced in \cite{DS}. The idea of this procedure is
rather simple: find the estimate which is consistent with the observed
data and has minimum $\ell_1$ norm (thus promoting
sparsity). Formally, assuming that the columns of $A$ are
normalized,\footnote{There is a slight modification when the columns
  do not have the same norm, namely, $\|D^{-1}A^*(y-Ax)\|_\infty \le
  \delta$, where $D$ is diagonal and whose diagonal entries are the
  $\ell_2$ norms of the columns of $A$.} the Dantzig selector is the
solution to the convex program
\begin{equation}
\label{eq:ds}
  \begin{array}{ll}
    \text{minimize}   & \quad \|x\|_1 \\
    \text{subject to} & \quad \|A^\T(y-Ax)\|_\infty \le \delta, 
  \end{array}
\end{equation}
where $\delta$ is a scalar. Clearly, the constraint is a data fitting
term since it asks that the correlation between the residual vector $r
= y - Ax$ and the columns of $A$ is small. Typically, the scalar
$\delta$ is adjusted so that the true $x_0$ is feasible, at least with
high probability, when the noise term $z$ is stochastic; that is,
$\delta$ obeys $\|A^\T z\|_\infty \le \delta$ (with high
probability). Another effective method, which we refer to as the
LASSO~\cite{lasso:tib}  % adding in reference as Ewout suggested (note: this is a new reference,
%   but I think we should have already had it, so this is good).
(also known as basis pursuit denoising, or BPDN),
assumes a different fidelity term and is the solution to
\begin{equation}
  \label{eq:lasso}
  \begin{array}{ll}
    \text{minimize}   & \quad \|x\|_1 \\
    \text{subject to} & \quad \|y-Ax\|_2 \le \epsilon, 
  \end{array}
\end{equation}
where $\epsilon$ is a scalar, which again may be selected so that the
true vector is feasible.  Both of these estimators are generally able to
accurately estimate nearly sparse vectors and it is, therefore, of
interest to develop effective algorithms for each that can deal with
problems involving thousands or even millions of variables and
observations. 

There are of course many techniques, which are perhaps more complicated than
\eqref{eq:ds} and \eqref{eq:lasso}, for recovering signals or images
from possibly undersampled noisy data. Suppose for instance that we
have noisy data $y$ \eqref{eq:model} about an $n\times n$ image $x_0$;
that is, $[x_0]_{ij}$ is an $n^2$-array of real numbers. Then to
recover the image, one might want to solve a problem of
this kind:
\begin{equation}
  \label{eq:W+TV}
  \begin{array}{ll}
    \text{minimize}   & \quad \|W x \|_1 + \lambda \|x\|_{\text{TV}}\\
    \text{subject to} & \quad \|y-Ax\|_2 \le \epsilon, 
  \end{array}
\end{equation}
% taking Ewout's suggestions
where $W$ is some (possibly nonorthogonal) transform such as an undecimated
wavelet transform enforcing sparsity of the image in this domain, and
$\|\cdot\|_{\text{TV}}$ is the total-variation norm defined as 
\[
\|x\|_{TV} := \sum_{i,j} \sqrt{|x[i+1,j] - x[i,j]|^2 + |x[i,j+1] - x[i,j]|^2}.
\]
\iffalse
\[
\|x\|_{TV} := \sum_{i,j} \sqrt{|(D_1 x) [i,j]|^2 + |(D_2 x)[i,j]|^2}, 
\]
where $D_1$ and $D_2$ are the horizontal and vertical differences
$(D_1 x)[i,j] = x[i+1,j] - x[i,j]$, $(D_2 x)[i,j] = x[i,j+1] -
x[i,j]$. 
\fi
The motivation for \eqref{eq:W+TV} is to look for a sparse
object in a transformed domain while reducing artifacts due to
sparsity constraints alone, such as Gibbs oscillations, by means of
the total-variation norm \cite{ROF,EdgePreserving,PracticalCS}.  The
proposal \eqref{eq:W+TV} appears computationally more involved than
both \eqref{eq:ds} and \eqref{eq:lasso}, and our goal is to develop
effective algorithms for problems of this kind as well.

To continue our tour, another problem that has recently attracted a
lot attention concerns the recovery of a low-rank matrix $X_0$ from
undersampled data
\begin{equation}
\label{eq:model2}
y = \cA(X_0) + z,
\end{equation}
where $\cA: \R^{n_1 \times n_2} \goto \R^{m}$ is a linear
operator supplying information about $X_0$. An important example
concerns the situation where only some of the entries of $X_0$ are
revealed, $\cA(X_0) = [X_0]_{ij} : (i,j) \in E \subset [n_1]
\times [n_2]$, and the goal is to predict the values of all the
missing entries. It has been shown \cite{CR08,CT09,Gross} that an
effective way of recovering the missing information from $y$ and the
model \eqref{eq:model2} is via the convex program
\begin{equation}
\label{eq:sdp}
  \begin{array}{ll}
    \text{minimize}   & \quad \|X\|_*\\
    \text{subject to} & \quad X \in {\cal C}.
 \end{array}
\end{equation}
Here, $\|X\|_*$ is the sum of the singular values of the matrix $X$, a
quantity known as the {\em nuclear norm} of $X$. ($\|X\|_*$ is also
the dual of the standard operator norm $\|X\|$, given by the largest
singular value of $X$). Above, ${\cal C}$ is a data
fitting set, and might be $\{X : \cA(X) = y\}$ in the
noiseless case, or $\{X : \|\cA^\T(y - \cA(X))\|_\infty \le \delta\}$
(Dantzig selector-type constraint) or $\{X : \|y - \cA(X)\|_2 \le
\epsilon\}$ (LASSO-type constraint) in the noisy setup.  We are again
interested in computational solutions to problems of this type.

\subsection{The literature}

There is of course an immense literature for solving problems of the
types described above. Consider the LASSO, for example.  Most of the
works
\cite{FPC2,FPCAS,BregmanOsher,Bregman,GPSR,sparsa,FISTA,CoordinateDescent}
are concerned with the unconstrained problem
\begin{equation}
  \label{eq:lasso_unc}
  \begin{array}{ll}
    \text{minimize}   & \quad \thalf \|Ax - b\|_2^2 + \lambda \|x\|_1,
  \end{array}
\end{equation}
which differs from \eqref{eq:lasso} in that the hard constraint 
$\|Ax - b\|_2 \le \epsilon$ is replaced with a quadratic
penalty $\thalf\lambda^{-1}\|Ax - b\|_2^2$. There are far fewer methods specially adapted to
\eqref{eq:lasso}; let us briefly discuss some of them. SPGL1
\cite{vandenberg_friedlander_2009} is an excellent solver specifically
designed for \eqref{eq:lasso}. The issue is that at the moment, it
cannot handle important variations such as
\begin{equation*}
  \begin{array}{ll}
    \text{minimize}   & \quad \|Wx\|_1 \\
    \text{subject to} & \quad \|y-Ax\|_2 \le \epsilon, 
  \end{array}
\end{equation*}
where $W$ is a nonorthogonal transform as in \eqref{eq:W+TV}. The main
reason is that SPGL1---as with almost all first-order methods for
that matter---relies on the fact that the proximity operator associated
with the $\ell_1$ norm,
\begin{equation}
  \label{eq:prox-l1}
   x(z;t) \triangleq \argmin_{x \in \R^n} \quad \thalf t^{-1} \|x - z\|_2^2 + \|x\|_1,
\end{equation}
is efficiently computable via soft-thresholding.
This is not the case, however, when $\|x\|_1$ is
replaced by a general term of the form $\|Wx\|_1$. 
NESTA \cite{NESTA}
can efficiently deal with an objective functional of the form
$\|Wx\|_1$, but it requires repeated projections onto the feasible
set; see also \cite{CSALSA} for a related approach.
Hence, NESTA is efficient when $AA^\T$ is a projector or, more
generally, when the eigenvalues of $AA^\T$ are well clustered. Other
types of algorithms such as LARS \cite{LARS} are based on homotopy
methods, and compute the whole solution path; \ie, they find the solution
to \eqref{eq:lasso_unc} for all values of the regularization parameter
$\lambda$ and, in doing so, find the solution to the constrained problem
\eqref{eq:lasso}. These methods do not scale well with problem size, however,
especially when the solution is not that sparse.

Turning to the Dantzig selector, solution algorithms are scarce.  The
standard way of solving \eqref{eq:ds} is via linear programming
techniques \cite{l1magic} since it is well known that it can be recast
as a linear program \cite{DS}. Typical modern solvers rely on interior-point methods which are somewhat problematic for large scale problems,
since they do not scale well with size. Another way of solving
\eqref{eq:ds} is via the new works \cite{DASSO,JustinDS}, which use
homotopy methods inspired by LARS to compute the whole solution path
of the Dantzig selector. These methods, however, are also unable to cope
with large problems. Another alternative is adapting SPGL1 to this
setting, but this comes with the caveat that it does not handle slight
variations as discussed above.

Finally, as far as the mixed norm problem \eqref{eq:W+TV} is
concerned, we are not aware of efficient solution algorithms. One can
always recast this problem as a second-order cone program (SOCP) which one 
could then solve via an interior-point method; but again, this is problematic for
large-scale problems. 

\subsection{Our approach}

In this paper, we develop a template for solving a variety of problems
such as those encountered thus far. The template proceeds as follows:
first, determine an equivalent \emph{conic formulation};
second, determine its \emph{dual}; third, apply \emph{smoothing}; and
fourth, solve using an \emph{optimal first-order method}.

\subsubsection{Conic formulation}
\label{subsec:conic}

In reality, our approach can be applied to general models expressed in the
following canonical form:
\begin{equation}
\label{eq:stdform}
  \begin{array}{ll}
    \text{minimize}   & \quad f(x) \\
    \text{subject to} & \quad \cA(x) + b \in \cK. 
  \end{array}
\end{equation}
The optimization variable is a vector $x\in\R^n$, and the objective
function $f$ is convex, possibly
extended-valued, and not necessarily smooth. The constraint is
expressed in terms of a linear operator $\cA:\R^n \rightarrow \R^m$, a
vector $b\in\R^m$, and a closed, convex cone $\cK\subseteq\R^m$. 
We shall
call a model of the form (\ref{eq:stdform}) that is equivalent to a
given convex optimization model $\mathcal{P}$ a \emph{conic form} for
$\mathcal{P}$. 

The conic constraint $\cA(x)+b\in\cK$ may seem specialized, but in
fact any convex subset of $\R^n$ may be represented in this fashion;
and models involving complex variables, matrices, or other vector
spaces can be handled by defining appropriate isomorphisms.  Of
course, some constraints are more readily transformed into conic form
than others; included in this former group are linear equations,
linear inequalities, and convex inequalities involving norms of affine
forms. Thus virtually every convex compressed sensing model may be
readily converted.  Almost all models admit multiple conic forms, and each
results in a different final algorithm.

For example, the Dantzig selector (\ref{eq:ds})
can be mapped to conic form as follows:
\begin{equation}
\label{eq:conic-ds}
	f(x) \rightarrow \|x\|_1, \qquad
	\cA(x) \rightarrow (A^\T Ax,0), \qquad
	b \rightarrow (-A^\T y,\delta), \qquad
	\cK \rightarrow \mathcal{L}_\infty^n, 
\end{equation}
where $\mathcal{L}_\infty^n$ is the epigraph of
the $\ell_\infty$ norm: $\mathcal{L}_\infty^n = \{(y,t) \in \R^{n+1}: \|y\|_\infty \le t\}$.

\subsubsection{Dualization}

The conic form (\ref{eq:stdform}) does not immediately lend itself to
efficient solution using first-order methods for two reasons: first,
because $f$ may not be smooth; and second, because projection onto the
set $\{x\,|\,\cA(x)+b\in\cK\}$, or even the determination of a single
feasible point, can be expensive.  We propose to resolve these issues
by solving either the dual problem, or a carefully chosen approximation
of it.  Recall that the dual of our canonical form \eqref{eq:stdform}
is given by
\begin{equation}
  \label{eq:conic-dual}
  \begin{array}{ll}
    \text{maximize}   & \quad g(\lambda) \\
    \text{subject to} & \quad \lambda\in\cK^*, 
  \end{array}
\end{equation}
where $g(\lambda)$ is the Lagrange dual function
\[
	g(\lambda) \triangleq \inf_x \mathcal{L}(x,\lambda) = \inf_x f(x) - \< \lambda, \cA(x) + b \>,
\]
and $\cK^*$ is the dual cone defined via
\[
\cK^* = \{\lambda \in \R^m : \<\lambda,x \> \ge 0 \text{ for all }
x \in \cK\}. 
\] 

The dual form has an immediate benefit that for the problems of interest,
projections onto the dual cone are usually tractable and
computationally very efficient. For example,
consider the projection
of a point onto the feasible set $\{x : \|Ax - y\|_2 \le \epsilon\}$
of the LASSO, an operation which may be expensive. However, one can
recast the constraint as $\cA(x)+b\in\cK$ with
\begin{equation}
\label{eq:conic-lasso}
\cA(x) \rightarrow (Ax,0), 
\qquad b \rightarrow (-y,\epsilon) \qquad
	\cK \rightarrow \mathcal{L}_2^m,
\end{equation}
where $\mathcal{L}_2^m$ is the second order cone ${\cal L}_2^m =
\{(y,t) \in \R^{m+1} : \|y\|_2 \le t\}$. This cone is self dual,
\ie, $({\cal L}_2^m)^* = {\cal L}_2^m$, and projection onto ${\cal
  L}_2^m$ is trivial: indeed, it is given by 
\begin{equation}
  \label{eq:projK2}
  (y,t) \mapsto \begin{cases} (y,t), & \|y\|_2 \le t,\\
    c (y,\|y\|_2), &  - \|y\|_2 \le  t \le \|y\|_2,\\
    (0,0), & t \le -\|y\|_2, 
\end{cases} \quad c =   \frac{\|y\|_2 + t}{2 \|y\|_2}. 
\end{equation}
And so we see that by eliminating the affine mapping, the
projection computation has been greatly simplified. Of course,
not every cone projection admits as simple a solution as
(\ref{eq:projK2}); but as we will show, all of the cones of
interest to us do indeed.

\subsubsection{Smoothing}

Unfortunately, because of the nature of the problems under study, the
dual function is usually not differentiable either, and direct
solution via subgradient methods would converge too slowly. Our
solution is inspired by the smoothing technique due to Nesterov
\cite{Nesterov05}.  We shall see that if one modifies the primal
objective $f(x)$ and instead solves
%and minimizes instead  % Ewout remarked about this.  We should at least not use ``minimize'' since the equation also has ``minimize''
\begin{equation}
\label{eq:stdform_mu}
  \begin{array}{ll}
    \text{minimize}   & \quad f_\mu(x) \triangleq f(x) + \mu d(x)\\
    \text{subject to} & \quad \cA(x) + b \in \cK, 
  \end{array}
\end{equation}
where $d(x)$ is a strongly convex function to be defined later and
$\mu$ a positive scalar, then the dual problem takes the form
\begin{equation}
  \label{eq:conic-dual_mu}
  \begin{array}{ll}
    \text{maximize}   & \quad g_\mu(\lambda) \\
    \text{subject to} & \quad \lambda\in\cK^*, 
  \end{array}
\end{equation}
where $g_\mu$ is a smooth approximation of $g$. This approximate
model can now be solved using first-order methods. As a general rule, higher
values of $\mu$ improve the performance of the underlying solver,
but at the expense of accuracy. Techniques such as continuation
can be used to 
% Ewout commented about this, but I think we are OK
recover the accuracy lost, 
however, so the precise
trade-off is not so simple.

In many cases, the smoothed dual can be reduced
to an unconstrained problem of the form
\begin{equation}
 \label{eq:smoothed-dual-comp}
 \begin{array}{ll}
  \text{maximize}   & \quad - \gs(z) - h(z),
\end{array}
\end{equation}
with optimization variable $z \in \R^m$, where $\gs$ is
convex and smooth and $h$ convex, nonsmooth, and possibly extended-valued.
For instance, for the Dantzig selector (\ref{eq:ds}), $h(z) = \delta \|z\|_1$.
As we shall see, this so-called \emph{composite form}
can also be solved efficiently using optimal first-order methods.
In fact, the reduction to composite form often simplifies some
of the central computations in the algorithms.

\subsubsection{First-order methods}
\label{sec:firstorder}

Optimal first-order methods are proper descendants of the classic
projected gradient algorithm. For the smoothed dual problem
(\ref{eq:conic-dual_mu}), a prototypical projected gradient algorithm 
begins with a point $\lambda_0\in\cK^*$, and generates 
updates for $k=0,1,2,...$ as follows:
\begin{equation}
	\label{eq:pg1}
 	\lambda_{k+1} \leftarrow \argmin_{\lambda\in\cK^*} \|\lambda_k + t_k \nabla g_\mu(\lambda_k) - \lambda \|_2,
\end{equation}
given step sizes $\{t_k\}$. The method has also been extended
to composite problems like (\ref{eq:smoothed-dual-comp}) \cite{Wright_talk,Nesterov07,tseng_2008};
the corresponding iteration is
\begin{equation}
	\label{eq:pgc}
 	z_{k+1} \leftarrow \argmin_y \gs(z_k) + \langle \nabla \gs(z_k), z - z_k \rangle + \tfrac{1}{2t_k} \| z - z_k \|^2 + h(z). 
\end{equation}
Note the use of a general norm $\|\cdot\|$ and the inclusion of
the nonsmooth term $h$. We call the minimization
in (\ref{eq:pgc}) a \emph{generalized projection}, because it
reduces to a standard projection (\ref{eq:pg1}) if the norm is
Euclidean and $h$ is an indicator function.
%But the generalized  % changing per Ewout's suggestion -- starting a sentecence with ``but'' isn't great
This generalized
form allows us to construct efficient algorithms for a wider
variety of models.

For the problems under study, the step sizes
$\{t_k\}$ above can be chosen so that $\epsilon$-optimality (that is,
$\sup_{\lambda\in\cK^*} g_\mu(\lambda)-g_\mu(\lambda_k)\leq\epsilon$)
can be achieved in
$\order(1/\epsilon)$ iterations \cite{NesterovBook}.
In 1983, Nesterov reduced this cost to $\order(1/\sqrt{\epsilon})$
using a slightly more complex iteration
\begin{equation}
	\label{eq:pga}
 	\lambda_{k+1}\leftarrow \argmin_{\lambda\in\cK^*} \|\nu_k + t_k \nabla g_\mu(\nu_k) - \lambda \|_2, \quad \nu_{k+1}\leftarrow \lambda_{k+1} + \alpha_k ( \lambda_{k+1} - \lambda_k ), 
\end{equation}
where $\nu_0=\lambda_0$ and the sequence $\{\alpha_k\}$ is constructed
according to a particular recurrence relation.
Previous work
by Nemirosvski and Yudin had established $\order(1/\sqrt{\epsilon})$ complexity as
the best that can be achieved for this class of problems \cite{NY83}, 
so Nesterov's modification is indeed optimal. 
Many alternative first-order
methods have since been developed
\cite{Nes88,Nesterov05,Nesterov07,tseng_2008,auslender_teboulle_2006,LLM09},
including methods that support generalized projections.
We examine these methods in more detail in \S
\ref{sec:implement}. 

We have not yet spoken about the complexity of computing
$g_\mu$ or $\gs$ and their gradients.
For now, let us highlight the fact that
%$\nabla g_\mu(\lambda) = b - \cA(x(\lambda))$, where
% Accepting Ewout's change:
$\nabla g_\mu(\lambda) = - \cA(x(\lambda)) -b$, where
\begin{equation}
	\label{eq:recoverx}
	x(\lambda) \triangleq \argmin_x \mathcal{L}_\mu(x,\lambda) = \argmin_x f(x) + \mu d(x) - \langle \cA(x)+b,\lambda \rangle,
\end{equation}
and $d(x)$ is a selected proximity function. In the common case that
$d(x)=\thalf\|x-x_0\|^2$,
the structure of (\ref{eq:recoverx}) is identical to that of a
generalized projection. Thus we see that the ability
to efficiently minimize the sum of a linear term, a proximity function,
and a nonsmooth function of interest is the fundamental computational
primitive involved in our method. Equation (\ref{eq:recoverx}) also
reveals how to recover an approximate primal solution as
$\lambda$ approaches its optimal value.

\subsection{Contributions}

The formulation of compressed sensing models in conic form is not
widely known.   % changing ``But'' to ``Yet'' per Ewout's suggestion
Yet the convex optimization modeling framework
\texttt{CVX} \cite{cvx} converts \emph{all} models into conic form;
and the compressed sensing package $\ell_1$-\textsc{Magic}
\cite{l1magic} converts problems into second-order cone programs
(SOCPs).  Both systems utilize interior-point methods instead of
first-order methods, however. As mentioned above, the smoothing step
is inspired by \cite{Nesterov05}, and is similar in structure to
traditional Lagrangian augmentation. As we also noted, first-order methods
have been a subject of considerable research.

Taken separately, then, none of the components in this approach is % Ewout's change
new. % Ewout's change: But to However
However their combination and application to solve compressed
sensing problems leads to effective algorithms that have not
previously been considered. For instance, applying our methodology to
the Dantzig selector gives a novel and efficient algorithm
(in fact, it gives several novel algorithms, depending on
which conic form is used).  Numerical experiments presented% removed ``being'', as Ewout suggests.  I agree.  Also change Introduced to Presented
later in the paper show that one can solve the Dantzig selector
problem with a reasonable number of applications of $A$ and its
adjoint; the exact number depends upon the desired level of accuracy. In
the case of the LASSO, our approach leads to novel algorithms
which are competitive with state-of-the-art methods such as SPGL1.

Aside from good empirical performance, we believe that the primary
merit of our framework lies in its flexibility.  To be sure, all the
compressed sensing problems listed at the beginning of this paper, and of
course many others, can be solved via this approach. These include
total-variation norm problems, $\ell_1$-analysis problems involving
objectives of the form $\|Wx\|_1$ where $W$ is neither orthogonal
nor diagonal, and so on. In each case, our framework allows us to 
construct an effective algorithm, thus providing a computational solution
to almost every problem arising in sparse signal or low-rank matrix
recovery applications.

Furthermore, in the course of our investigation, we have developed a
number of additional technical contributions. For example, we will
show that certain models, including the Dantzig selector, exhibit 
an \emph{exact penalty} property: the exact solution to the original
problem is recovered even when some smoothing is applied. We have also
developed a novel continuation scheme that allows us to employ more
aggressive smoothing to improve solver performance while still
recovering the exact solution to the unsmoothed problem.  The
flexibility of our template also provides opportunities to employ
novel approaches for controlling the step size.

\subsection{Software}

This paper is accompanied with a software release \cite{Templates},
including a detailed user guide which gives many additional
implementation details not discussed in this paper. 
Since most compressed sensing
problems can be easily cast into standard conic form, our software
provides a powerful and flexible computational tool for solving a
large range of problems researchers might be interested in
experimenting with.

The software is not a single, monolithic solver; rather, it is a
suite of programs and routines designed to
serve as building blocks for constructing complete algorithms.
Roughly speaking, we can divide the routines into three levels.
On the first level is a suite of routines that implement a
variety of known first-order solvers, including the standard
projected gradient algorithm and known optimal variants by
Nesterov and others. On the second level are wrappers designed
to accept problems in conic standard form (\ref{eq:stdform})
and apply the first-order solvers to the smoothed dual problem.
Finally, the package includes a variety of routines to
directly solve the specific models described in this paper
and to reproduce our experiments.

We have worked to ensure that each of the solvers is as easy to use as possible, by providing sensible defaults for line search, continuation,
and other factors. At the same time, we have sought to give the user
flexibility to insert their own choices for these components. We also
want to provide the user with the opportunity to compare the
performance of various first-order variants on their particular
application. We do have some general views about which algorithms
perform best for compressed sensing applications, however, and will
share some of them in \S\ref{sec:numerics}.

\subsection{Organization of the paper}

In \S\ref{sec:conic}, we continue the discussion of conic
formulations, including a derivation of the dual conic formulation and
details about the smoothing operation.  Section \ref{sec:dantzig}
instantiates this general framework to derive a new algorithm for the
Dantzig selector problem. In \S\ref{sec:examples}, we provide further selected
instantiations of our framework including the LASSO, total-variation
problems, $\ell_1$-analysis problems, and common 
%minimum nuclear-norm  % Accepting Ewout's suggestion
nuclear-norm minimization
problems.  In \S\ref{sec:implement}, we review a variety of
first-order methods and suggest improvements.  Section
\ref{sec:numerics} presents numerical results illustrating both the
empirical effectiveness and the flexibility of our approach. Section
\ref{sec:software} provides a short introduction to the software
release accompanying this paper. Finally, the appendix proves the
exact penalty property for general linear programs, and thus for
Dantzig selector and basis pursuit models, and describes a unique
approach we use to generate test models so that their
exact solution is known in advance.

\section{Conic formulations}
\label{sec:conic}

\subsection{Alternate forms}
\label{subsec:alternate}

In the introduction, we presented our standard conic form
(\ref{eq:stdform}) and a specific instance for Dantzig selector
in (\ref{eq:conic-ds}).  As we said then, conic forms are rarely
unique; this is true even if one disregards simple scalings of the cone
constraint. For instance, we may express the Dantzig selector
constraint as an intersection of linear inequalities, $-\delta \vec{1}
\preceq A^\T(y-Ax) \preceq \delta \vec{1}$, suggesting the following
alternative:
\begin{equation}
	\label{eq:ds-lp}
	f(x) \rightarrow \|x\|_1, \qquad
	\cA(x) \rightarrow \bmat{ -A^\T A \\ A^\T A } x, \qquad
	b \rightarrow \bmat{ \delta \vec{1} + A^\T y \\ \delta \vec{1} - A^\T y }, \qquad
	\cK \rightarrow \R^{2n}_+. 
\end{equation}
We will return to this alternative later in \S\ref{sec:alternate-ds}.
In many instances, a conic form may involve the manipulation of the objective
function as well. For instance, if we first transform (\ref{eq:ds}) to
\begin{equation*}
	\begin{array}{ll}
          \text{minimize}   & \quad t\\
          \text{subject to} & \quad \|x\|_1 \leq t \\ 
          & \quad \|A^\T(y-Ax)\|_\infty \leq \delta, 
	\end{array}
\end{equation*}
then yet another conic form results:
\begin{equation}
	f(x,t)  \rightarrow t, \qquad
	\cA(x,t) \rightarrow (x,t,A^\T Ax,0), \qquad
	b \rightarrow (0,0,-A^\T y,\delta), \qquad
	\cK \rightarrow \mathcal{L}_1^n\times\mathcal{L}_\infty^n,
\end{equation}
where $\mathcal{L}_1^n$ is the epigraph of the $\ell_1$ norm,
$\mathcal{L}_1^n = \{(y,t) \in \R^{n+1}: \|y\|_1 \le t\}$.

Our experiments show that different conic formulations yield different
levels of performance using the same numerical algorithms. Some are
simpler to implement than others as well. Therefore, it is worthwhile
to at least explore these alternatives to find the best choice for a
given application.
 
\subsection{The dual}

To begin with, the conic Lagrangian associated with \eqref{eq:stdform}
is given by
\begin{equation}
	\label{eq:Lagrangian}
	\lag(x,\lambda) = f(x) - \langle \lambda, \cA(x)+b \rangle, 
\end{equation}
where $\lambda\in\R^m$ is the Lagrange multiplier, constrained to lie
in the dual cone $\cK^\T$. 
The dual function $g:\R^m\rightarrow(\R\cup-\infty)$ is, therefore,
\begin{equation}
  g(\lambda) = \inf_x {\cal L}(x,\lambda) = - f^\T(\cA^\T(\lambda)) - \<b,\lambda\>.
\end{equation}
Here, $\cA^\T:\R^m\rightarrow\R^n$ is the adjoint of the linear
operator $\cA$ and $f^\T:\R^n\rightarrow(\R\cup+\infty)$ is the convex
conjugate of $f$ defined by  
\[
f^\T(z)=\sup_x \langle z,x\rangle-f(x).
\] 
Thus the dual problem is given by 
\begin{equation}
  \label{eq:conic-dual_withAffine}
  \begin{array}{ll}
    \text{maximize}   & \quad -f^\T(\cA^\T(\lambda)) - \<b,\lambda\>\\
    \text{subject to} & \quad \lambda\in \cK^\T.
  \end{array}
\end{equation}
% \EJC{Why is the affine case interesting and deserves to be
%   explicated?} \MCG{I included it because it occurs whenever we
%   have chosen to move the nonlinear objective into a constraint---for
%   instance, with $\ell_1$-analysis. It is not interesting right here per se,
%   but I do think that the simplification it provides in the next section
%   (see your next comment) is interesting. I am happy to delete these 
%   however if you think it does not add enough to the discussion.}
%   In the case where our objective is an affine function
% $f(x)= \<c_0,x\> + d_0$, the dual becomes
% \begin{equation}
%   \label{eq:conic-dual_affineObj}
%   \begin{array}{ll}
%     \text{maximize}   & \quad d_0 -\<b,\lambda\>\\
%     \text{subject to} & \quad \cA^*(\lambda) = c_0 \\
%                       & \quad \lambda\in \cK^\T.
%   \end{array}
% \end{equation}

Given a feasible primal/dual pair $(x,\lambda)$,
the \emph{duality gap} is the difference between their respective objective
values. The non-negativity of the duality gap is easily verified:
\begin{equation}
	f(x)-g(\lambda) = f(x)+f^\T(\cA^\T(\lambda))+\<b,\lambda\> \geq \langle x, \cA^\T(\lambda) \rangle + \<b,\lambda\>
	= \langle \cA(x) + b, \lambda \rangle \geq 0.
\end{equation}
The first inequality follows from the definition of conjugate
functions, while the second follows from the definition of the dual
cone. If both the primal and dual are strictly feasible---as is the
case for all problems we are interested in here---then the minimum
duality gap is exactly zero, and there exists an optimal pair
$(x^\star,\lambda^\star)$ that achieves $f(x^\star) = g(\lambda^\star)
= {\cal L}(x^\star,\lambda^\star)$. It is important to note that
the optimal points are not necessarily unique; more about this in
\S\ref{subsec:smoothing}. But any optimal primal/dual pair
will satisfy optimality conditions
\begin{equation}
	\label{eq:optcond}
	\cA(x^\star)+b\in\cK, \quad \lambda^\star\in\cK^\T, \quad \langle \cA(x^\star)+b,\lambda^\star \rangle = 0, \quad
	\cA^\T(\lambda^\star) \in \partial f(x^\star),
\end{equation}
where $\partial f$ refers to the subgradient of $f$.

\subsection{The differentiable case}
\label{sec:smoothcase}

The dual function is of course concave; and its derivative (when it exists) is given by
\begin{equation}
	\label{gderivative}
	\nabla g(\lambda) = - \cA(x(\lambda)) - b, \quad x(\lambda) \in \argmin_x {\cal L}(x,\lambda).
\end{equation}
It is possible that the minimizer $x(\lambda)$ is not unique, so in  % accepting Ewout's changes
order to be differentiable, all such minimizers must yield
the same value of $-\cA(x(\lambda))-b$.

If $g$ is finite and differentiable on the entirety of $\cK^\T$,
then it becomes trivial to locate an initial dual point (\eg, $\lambda=0$);
and for many genuinely useful cones $\cK^\T$, it becomes trivial to
project an arbitrary point $\lambda\in\R^m$ onto this feasible set. If
the $\argmin$ calculation in (\ref{gderivative}) is computationally
practical, we may entertain the construction of a projected gradient
method for solving the dual problem \eqref{eq:conic-dual_withAffine} directly; 
\ie, without our proposed smoothing step. Once an optimal
dual point $\lambda^\star$ is recovered, an optimal solution to the original
problem (\ref{eq:stdform}) is recovered by solving $x^\star\in\argmin_x
{\cal L}(x,\lambda^\star)$.

Further suppose that $f$ is strongly convex; that is, it satisfies for
some constant $m_f > 0$,
\begin{equation}
  f((1-\alpha) x + \alpha x') \le  (1-\alpha) f(x) + \alpha f(x')  -  m_f\alpha(1-\alpha) \|x - x'\|_2^2/2
\end{equation}
for all $x, x' \in \dom(f)$ and $0 \le \alpha \le 1$. Then assuming
the problem is feasible, it admits a unique optimal solution. The
Lagrangian minimizers $x(\lambda)$ are unique for all
$\lambda\in\R^n$; so $g$ is differentiable everywhere.  Furthermore,
\cite{Nesterov05} proves that the gradient of $g$ is Lipschitz
continuous, obeying
\begin{equation}
	\label{eq:glipschitz}
        \| \nabla g(\lambda') - \nabla g(\lambda) \|_2 \leq m_f^{-1}\|\cA\|^2\|\lambda'-\lambda\|_2, 
\end{equation}
where $\|\cA\|=\sup_{\|x\|_2 = 1}\|\cA(x)\|_2$ is the induced operator
norm of $\cA$.  So when $f$ is strongly convex, then provably convergent,
accelerated gradient methods in the Nesterov style are possible.

\subsection{Smoothing}
\label{subsec:smoothing}

Unfortunately, it is more likely that $g$ is not differentiable
(or even finite) on all of $\cK^\T$. So we consider a smoothing
approach similar to that proposed in \cite{Nesterov05} to solve an
approximation of our problem. Consider the following perturbation of
(\ref{eq:stdform}):
\begin{equation}
\label{eq:dsmu2}
  \begin{array}{ll}
    \text{minimize}   & \quad f_\mu(x) \triangleq f(x) + \mu d(x)\\
    \text{subject to} & \quad \cA(x) + b \in \cK
  \end{array}
\end{equation}
for some fixed smoothing parameter $\mu>0$ and a strongly convex function $d$ obeying
\begin{equation}
	d(x)\geq d(x_0)+\thalf\|x-x_0\|^2
\end{equation}
for some fixed point $x_0\in\R^n$.
Such a function is usually called a proximity function.

The new objective $f_\mu$ is strongly convex with $m_f=\mu$, so the full benefits
described in \S\ref{sec:smoothcase} now apply.  The Lagrangian and
dual functions become\footnote{One can also observe the identity
  $g_\mu(\lambda) = \sup_z - f^\T(\cA^\T(\lambda)-z) - \mu
  d^\T(\mu^{-1} z) - \<b,\lambda\>.$}
\begin{gather}
  \label{eq:dsmu2-lagr}
  {\cal L}_\mu(x,\lambda) = f(x) + \mu d(x) - \langle \lambda, \cA(x) + b \rangle \\
\label{eq:dsmu2-dfunc}
  g_\mu(\lambda)   \triangleq \inf_x {\cal L}_\mu(x,\lambda)   = - (f + \mu
  d)^\T(\cA^\T(\lambda)) - \<b, \lambda\>.
\end{gather}
% \EJC{Why do we need the second equation in (2.15)? Also, I have the
%   same comment as before concerning the affine objective below.} \MCG{I
%   admit it is hard to see in the generic derivation, but ultimately if
%   we wanted to write out the smoothed dual function for a specific model,
%   such as the Dantzig selector, this second equation would be the most useful
%   to that effort. For instance, when $f=\|x\|_1$, $f^*$ is the indicator
%   function of the $\ell_\infty$ ball. Therefore, $f^*(\cA^*(\lambda)-z)$ translates
%   to an domain constraint $\|\cA^*(\lambda)-z\|_\infty\leq 1$.
%   Furthermore, writing this second equation makes it more
%   obvious how I get the smoothed dual formula for the affine case. Again, I'm
%   happy to delete these if you feel it best. Or perhaps I can clarify earlier
%   that some useful models do indeed have an affine objective.}
One can verify that for the affine objective case $f(x)\triangleq
\<c_0, x\> + d_0$, the dual and smoothed dual function take the form
\begin{align*}
%g(\lambda) & = d_0 - I_{\ell_\infty}(\cA^*(\lambda)-c_0) - \<b,\lambda\>,\\
g(\lambda) & = d_0 - I_{\{0\}}(\cA^*(\lambda)-c_0) - \<b,\lambda\>,\\
g_\mu(\lambda) & = d_0 - \mu d^*(\mu^{-1}(\cA^*(\lambda)-c_0)) -
\<b,\lambda\>,
\end{align*}
where $I_{\{0\}}$ is the indicator function of the set $\{0\}$; that is, 
\[
I_{\{0\}}(y) \triangleq \begin{cases} 0, & y=0, \\ +\infty, & y \neq 0. \end{cases}
\]
%where $I_{\ell_\infty}$ is the indicator function of the $\ell_\infty$
%norm ball; that is,
%\[
	%I_{\ell_\infty}(y) \triangleq \begin{cases} 0, & \|y\|_\infty \leq 1, \\ +\infty, & \|y\|_\infty > 1. \end{cases}
%\]
The new optimality conditions are
\begin{equation}
	\label{eq:optcondsm1}
	\cA(x_\mu)+b\in\cK, \quad \lambda_\mu\in\cK^\T, \quad \langle \cA(x_\mu)+b,\lambda_\mu \rangle = 0, \quad
	\cA^\T(\lambda_\mu) - \mu \nabla d(x_\mu) \in \partial f(x_\mu).
\end{equation}
Because the Lipschitz bound (\ref{eq:glipschitz}) holds, first-order
methods may be employed to solve (\ref{eq:dsmu2}) with provable
performance.  The iteration counts for these methods are proportional
to the square root of the Lipschitz constant, and therefore
proportional to $\mu^{-1/2}$.  There is a trade-off between the
accuracy of the approximation and the performance of the algorithms
that must be explored.\footnote{In fact, even when the original
  objective is strongly convex, further adding a strongly convex term
  may be worthwhile to improve performance.}

For each $\mu>0$, the smoothed model obtains a single minimizer
$x_\mu$; and the trajectory traced by $x_\mu$ as $\mu$ varies
converges to an optimal solution $x^\star\triangleq \lim_{\mu\rightarrow
  0_+} x_\mu$. Henceforth, when speaking about the (possibly
non-unique) optimal solution $x^\star$ to the original model, we will be
referring to this uniquely determined value. Later we will show that
for some models, including the Dantzig selector, the approximate model
is \emph{exact}: that is, $x_\mu=x^\star$ for sufficiently small but
nonzero $\mu$.

\begin{figure}
\centering
\includegraphics[width=0.6\textwidth]{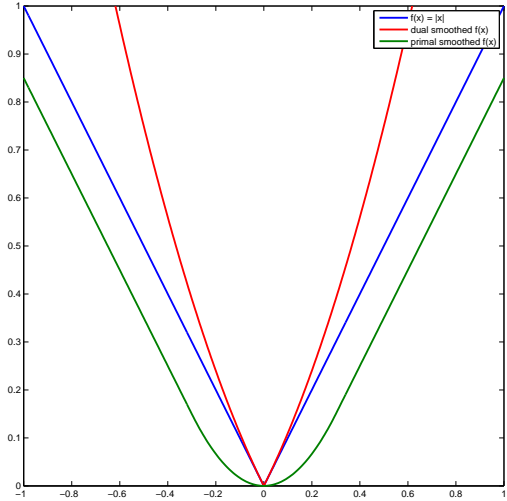}\\[-20pt]
\caption{The original objective $f(x)=|x|$ (blue), our modification (red), and Nesterov's smoothing (green).}
\label{fig:smoothing}
\end{figure}

Roughly speaking, the smooth dual function $g_\mu$ is what we
would obtain if the Nesterov smoothing method described in
\cite{Nesterov05} were applied to the dual function $g$.
It is worthwhile to explore how things would differ if the Nesterov
approach were applied directly to the primal
objective $f(x)$. Suppose that $f(x)=\|x\|_1$ and $d(x)=\thalf\|x\|_2^2$.
The Nesterov approach yields a smooth approximation $f^{\text{N}}_\mu$
whose elements can be described by the formula
\begin{equation}
	\left[ f^{\text{N}}_\mu(x) \right]_i = \sup_{|z| \leq 1} z x_i - \thalf \mu z^2 =
	\begin{cases} \thalf \mu^{-1} x_i^2, & |x_i| \leq \mu, \\ |x_i|-\thalf \mu, & |x|\geq \mu, \end{cases} \quad i= 1, 2, \dots, n. 
\end{equation}
Readers may recognize this as the Huber penalty function with half-width $\mu$;
a graphical comparison with $f_\mu$ is provided in Figure~\ref{fig:smoothing}.
Its smoothness
may seem to be an advantage over our choice $f_\mu(x)=\|x\|_1+\thalf\|x\|_2^2$, but the
difficulty of projecting onto the set $\{x\,|\,\cA(x)+b\in\cK\}$ remains;
so we still prefer
to solve the dual problem. Furthermore, the quadratic behavior of
$f^{\text{N}}_\mu$ around $x_i=0$ eliminates the tendency
towards solutions with many zero values. 
In contrast, $f_\mu(x)$ maintains the
sharp vertices from the $\ell_1$ norm that are known to
encourage sparse solutions. 

\subsection{Composite forms}
\label{sec:composite}

In most of the cases under study, the dual variable can be 
partitioned as $\lambda\triangleq(z,\tau)\in\R^{m-\bar{m}}\times\R^{\bar{m}}$ such that
the smoothed dual $g_\mu(z,\tau)$ is linear in $\tau$; that is,
\begin{equation}
	g_\mu(\lambda) = - \gs(z) - \< v_\mu, \tau \>
\end{equation}
for some smooth convex function $\gs$ and a constant vector
$v_\mu\in\R^{\bar{m}}$. An examination of the Lagrangian
$\mathcal{L}_\mu$ (\ref{eq:dsmu2-lagr}) reveals a precise condition under which
this occurs: when the linear operator $\cA$ is of the form
$\cA(x)\rightarrow(\bar\cA(x),\vec{0}_{\bar{m}\times 1})$, 
as seen in the conic constraints for the Dantzig selector (\ref{eq:conic-ds})
and LASSO (\ref{eq:conic-lasso}).
If we partition $b=(\bar{b},b_\tau)$ accordingly,
then evidently $v_\mu=b_\tau$. 

Under such conditions, it is more natural to work with 
$\bar{\cA}$, $\bar{b}$, and $b_\tau$ directly, and exploit some
useful simplifications. Specifically, let us define the function
\begin{equation}
	\label{eq:hdefine}
	h:\R^{m - \bar{m}}\rightarrow(\R\cup+\infty), \quad 
	h(z) \triangleq \inf \{ \<b_\tau,\tau\> \,|\, (z,\tau)\in\cK^\T \}. 
\end{equation}
Then the dual problem reduces to a nonsmooth unconstrained maximization
\[
	\begin{array}{ll}
		\text{maximize} & \bar{g}_\mu(z) \triangleq - \gs(z) - h(z).
	\end{array}		
\]
The gradient of $\gs$ is $\nabla\gs(z)=\bar\cA(x(z))+\bar{b}$, where
$x(z)$ is the minimizer of a reduced Lagrangian
\begin{gather}
	\bar{\mathcal{L}}_\mu(x,z) = f(x) + \mu d(x) - \< z, \bar\cA(x) + \bar b\>.
\end{gather}

\subsection{Projections}
\label{sec:projections}

A standard gradient projection step for the smoothed dual problem is
% SRB, Ewout: we need to change this!  Add in a lambda, and be careful with minus signs, and do we really mean t_k^{-1} ?
\begin{equation}
	\label{eq:stdproj}
	\lambda_{k+1} \leftarrow \argmin_{\lambda\in\cK^*} \| \lambda - \lambda_k - t_k \nabla g_\mu(\lambda_k) \|_2.
\end{equation}
For the composite version of the same problem, % accepting Ewout's change
the corresponding generalized projection is
\begin{equation}
	\label{eq:pgcr}
 	z_{k+1} \leftarrow \argmin_z \gs(z_k) + \langle \nabla \gs(z_k), z - z_k \rangle + \tfrac{1}{2t_k} \| z - z_k \|^2 + h(z).
\end{equation}
Integrating the definition of $\nabla\gs(z)$ into  (\ref{eq:pgcr}) and simplifying yields a two-sequence recursion:
\begin{equation}
	\label{eq:gptwo}
	\begin{aligned}
	x_k &\leftarrow \argmin_x f(x) + \mu d(x) - \< \bar\cA^*(z_k), x \> \\
 	z_{k+1} &\leftarrow \argmin_z h(z) +  \tfrac{1}{2t_k} \| z - z_k \|^2 + \langle \bar\cA(x_k)+\bar b, z \rangle.  \\
 	\end{aligned} 	
\end{equation}
Note the similarity in computational structure of the two formulae. This similarity
is even more evident in the common scenario where
$d(x)=\thalf\|x-x_0\|^2$ for some fixed $x_0\in\R^n$.

Let $\Sigma$ be the matrix composed 
of the first $m-\bar{m}$ rows of the $m\times  m$ identity matrix, so that
$\Sigma \lambda\equiv z$ for all $\lambda=(z,\tau)$.
Then (\ref{eq:pgcr}) can also be written in terms of $\cK^\T$ and $g_\mu$:
\begin{equation}
	\label{eq:pgc2}
 	z_{k+1} \leftarrow \Sigma \cdot \argmax_{\lambda\in\cK^\T} g_\mu(\lambda_k) + 
 	\langle \nabla g_\mu(\lambda_k), \lambda - \lambda_k \rangle 
 	- \tfrac{1}{2t_k} \| \Sigma (\lambda - \lambda_k) \|^2,
\end{equation}
where $\lambda_k\triangleq(z_k,h(z_k))$. If $m=0$ and the
norm is Euclidean, then $\Sigma=I$ and the standard
projection (\ref{eq:stdproj}) is recovered. So (\ref{eq:pgcr})
is indeed a true generalization, as claimed in \S\ref{sec:firstorder}.

The key feature of the composite approach, then,
is the removal of the linear variables $\tau$ from the proximity term.
Given that they are linearly involved to begin with, this
yields a more accurate approximation of the dual function,
so we might expect a composite approach to yield improved performance.
In fact, the theoretical predictions of the number of iterations required
to achieve a certain level of accuracy are identical; and in
our experience, any differences in practice seem minimal
at best. The true advantage to the composite
approach is that generalized projections more readily admit analytic
solutions and are less expensive to compute.

\section{A Novel Algorithm for the Dantzig selector}
\label{sec:dantzig}

We now weave together the ideas of the last two sections to develop a
novel algorithm for the Dantzig selector problem \eqref{eq:ds}. 

\subsection{The conic form}

We use the standard conic formulation \eqref{eq:stdform} with the mapping
\eqref{eq:conic-ds} as discussed in the introduction, which results in the model
\begin{equation} 
  \label{eq:ds-conic}
  \begin{array}{ll}
    \text{minimize}   & \quad \|x\|_1 \\
    \text{subject to} & \quad (A^\T (y-Ax),\delta)\in\mathcal{L}_\infty^n, 
  \end{array}
\end{equation}
where $\mathcal{L}_\infty^m$ is the epigraph of the
$\ell_\infty$ norm. The dual variable $\lambda$, therefore, will
lie in the dual cone $(\mathcal{L}_\infty^n)^*=\mathcal{L}_1^n$, 
the epigraph of the $\ell_1$ norm. Defining $\lambda=(z,\tau)$, the
conic dual (\ref{eq:conic-dual_withAffine}) is
\begin{equation}
	\label{eq:ds-dual}
  \begin{array}{ll}
    \text{maximize}   & \quad - I_{\ell_\infty}(-A^\T Az) - \<A^\T y,z\> - \delta \tau \\
    \text{subject to} & \quad (z,\tau)\in\mathcal{L}_1^n,
  \end{array}
\end{equation}
where $f^*=I_{\ell_\infty}$ is the indicator function of the
$\ell_\infty$ norm ball as before.  Clearly the optimal value of
$\tau$ must satisfy $\|z\|_1=\tau$,\footnote{We assume $\delta > 0$ here; if $\delta = 0$,
the form is slightly different.} so eliminating it yields
\[
  \begin{array}{ll}
    \text{maximize}   & \quad - I_{\ell_\infty}(-A^\T Az) - \<A^\T y,z\> - \delta \|z\|_1.\\
  \end{array}
\]
In both cases, the dual objectives are not smooth, so the smoothing approach
discussed in \S\ref{subsec:smoothing} will indeed be necessary.

\subsection{Smooth approximation}\label{sec:smoothApprox}

We augment the objective with a strongly convex term
\begin{equation} 
  \label{eq:dsmu}
  \begin{array}{ll}
    \text{minimize}   & \quad \|x\|_1 + \mu d(x) \\
    \text{subject to} & \quad (A^\T (y-Ax),\delta)\in\cK\triangleq\mathcal{L}_\infty^n. 
  \end{array}
\end{equation}
The Lagrangian of this new model is
\[
\lag_\mu(x;z,\tau) = \|x\|_1 + \mu d(x) - \langle z,A^\T
(y-Ax)\rangle-\delta\tau.
\]
Letting $x(z)$ be the unique minimizer of $\lag_\mu(x;z,\tau)$, 
the dual function becomes
\[
g_\mu(z,\tau) = \|x(z)\|_1 + \mu d(x(z)) - \langle
z,A^\T(y-Ax(z))\>  - \tau \delta. 
\]
Eliminating $\tau$ per
\S\ref{sec:composite} yields a composite form 
$\bar{g}_\mu(z)=-\gs(z)-h(z)$ with
\[
  \gs(z) = - \|x(z)\|_1 - \mu d(x(z)) + \langle z,A^\T(y-Ax(z))\>, \quad
  h(z) = \delta \|z\|_1.
\] 
The gradient of $\gs$ is $\nabla \gs(z) = A^\T(y-Ax(z))$.

The precise form of $x(z)$ and $\nabla\gs$  depend of course on
our choice of proximity function $d(x)$.
For our problem, the simple convex quadratic
\[
d(x) = \thalf \|x - x_0\|_2^2,  
\]
for a fixed center point $x_0 \in \R^n$, works well, and guarantees
that the gradient is Lipschitz continuous with a constant at
most $\mu^{-1} \|A^\T A\|^2$. With this choice,
$x(z)$ can be expressed in terms of the \emph{soft-thresholding} operation
which is a common fixture in algorithms for sparse recovery.  
For scalars $x$ and $\tau \ge 0$, define
\[
  \ST(x,\tau) =  \sgn(x) \cdot \max\{|x|-\tau,0\} =
\begin{cases} 
  x + \tau, & x \leq -\tau, \\ 0, & |x| \leq \tau, \\ x - \tau, & x
  \geq \tau.
\end{cases}
\]
When the first input $x$ is a vector, the soft-thresholding operation is to be
applied componentwise. Armed with this definition, the formula for $x(z)$ becomes
\[
x(z) = \ST(x_0 - \mu^{-1} A^\T Az,\mu^{-1}). 
\]
If we substitute $x(z)$ into the formula for $\gs(z)$ and simplify carefully, we find that
\[
\gs(z) = - \thalf \mu^{-1} \| \ST( \mu x_0 - A^\T Az,1) \|_2^2 + \<
A^\T y, z \> + c,
\]
where $c$ is a term that depends only on constants $\mu$ and $x_0$. In other words, to within
an additive constant, $\gs(z)$ is a smooth approximation of the nonsmooth
term $I_{\ell_\infty}(-A^\T Az)+\<A^\T y,z\>$ from (\ref{eq:ds-dual}), and indeed it converges to that function
as $\mu\rightarrow 0$.

For the dual update, the generalized projection is
\begin{equation}
	\label{eq:ds-dual_update}
	z_{k+1} \leftarrow \argmin_z \gs(z_k) + \<\nabla\gs(z_k),z-z_k\> + \tfrac{1}{2t_k} \|z-z_k\|_2^2 + \delta \|z\|_1. 
\end{equation}
A solution to this minimization can also be expressed in terms of the soft thresholding operation:
\[
  z_{k+1} \leftarrow \ST( z_{k} - t_k A^\T(y-Ax(z_{k})), t_k \delta).
%z_{k+1} \leftarrow \ST( z_{k-1} - t_k A^\T(y-Ax(z_{k-1})), t_k \delta). % Ewout noticed this
\]

\subsection{Implementation}
To solve the model presented in \S\ref{sec:smoothApprox},
we considered first-order projected gradient solvers.
% Adding the above sentences, taking Ewout's suggestions
After some experimentation,
we concluded that the Auslender and Teboulle first-order variant
\cite{auslender_teboulle_2006,tseng_2008} is a good choice for this model.
We discuss this and other variants in more detail in \S\ref{sec:implement}, 
so for now we will simply
present the basic algorithm in Listing~\ref{alg:ds-alg}
below. Note that the dual update used 
differs slightly from (\ref{eq:ds-dual_update}) above:
the gradient $\nabla\gs$ is evaluated at $y_k$,
not $z_k$, and the step size in the generalized
projection is multiplied by $\theta_k^{-1}$.
Each iteration requires two applications of both $A$ and $A^\T$,  % Ewout's change
and is computationally inexpensive when a fast matrix-vector
multiply is available. 
\begin{algorithm}[H]
\caption{Algorithm for the smoothed Dantzig selector}
\label{alg:ds-alg}
\begin{algorithmic}[1]
\REQUIRE $z_0,x_0\in\R^n$, $\mu>0$, step sizes $\{t_k\}$
  \STATE $\theta_0\leftarrow 1$, $v_0\leftarrow z_0$
  \FOR{$k=0,1,2,\dots$}
  \STATE $y_k \leftarrow (1-\theta_k) v_k + \theta_k z_k$
  \STATE $x_k \leftarrow \ST(x_0 - \mu^{-1} A^\T A y_k, \mu^{-1})$.
  \STATE $z_{k+1} \leftarrow \ST(z_k - \theta_k^{-1} t_k A^\T(y-Ax_k), \theta_k^{-1} t_k \delta)$
  \STATE $v_{k+1} \leftarrow (1-\theta_k) v_k + \theta_k z_{k+1}$
  \STATE $\theta_{k+1} \leftarrow 2/(1+(1+4/\theta_k^2)^{1/2})$
  \ENDFOR
\end{algorithmic}  
\end{algorithm}

It is known that for a fixed step size $t_k\triangleq t \leq
\mu/\|A^*A\|_2^2$, the above algorithm converges in the sense that
$\bar{g}_\mu(z^*) - \bar{g}_\mu(z_k) = \order(k^{-2})$. Performance
can be improved through the use of a backtracking line search for
$z_k$, as discussed in \S\ref{sec:stepsize} (see also
\cite{FISTA}). Further, fairly standard arguments in convex analysis
show that the sequence $\{x_k\}$ converges to the unique solution to
\eqref{eq:dsmu}.

\subsection{Exact penalty} \label{sec:exactRecovery}

Theorem 3.1 in \cite{SVT} can be adapted to show that as $\mu \goto
0$, the solution to \eqref{eq:dsmu} converges to a solution to
\eqref{eq:ds}. But an interesting feature of the Dantzig selector
model in particular is that if $\mu < \mu_0$ for $\mu_0$ sufficiently
small, the solutions to the Dantzig selector and to its perturbed
variation \eqref{eq:dsmu} coincide; that is, $x^\star=x^\star_\mu$. In
fact, this phenomenon holds for any linear program (LP).

%\begin{theorem}\label{thm:exactRecovery}
    %There exists $\mu_0 > 0$ such that for all $\mu < \mu_0$,
    %the solution to \eqref{eq:dsmu} is also a solution to \eqref{eq:ds}.
%\end{theorem}
\begin{theorem}[Exact penalty]\label{thm:exactRecovery}
  Consider an arbitrary LP with objective $\<c,x\>$ and having an
  optimal solution (\ie, the optimal value is not $-\infty$) and let
  $Q$ be a positive semidefinite matrix.  Then there is a $\mu_0 > 0$
  such that if $0 < \mu \le \mu_0$, any solution to the perturbed
  problem with objective $\<c,x\> + \thalf \mu \<x-x_0, Q(x-x_0)\>$ is
  a solution to LP.  Among all the solutions to LP, the solutions to
  the perturbed problem are those minimizing the quadratic penalty. In
  particular, in the (usual) case where the LP solution is unique, the
  solution to the perturbed problem is unique and they coincide.
\end{theorem}

The theorem is proved in Appendix~\ref{app:exact}. There are also two related results in recent literature: a proof of the special case of noiseless basis pursuit \cite{YinPenalty}, and a more general proof \cite{FrieTsen:2007} that allows for a range of penalty functions. Our theorem is a special case of \cite{FrieTsen:2007}, but we present the proof since it uses a different technique than \cite{FrieTsen:2007}, is directly applicable to our method in this form, and provides useful intuition.
The result, when combined with our continuation techniques in \S\ref{sec:continuation}, is also a new look at
known results about the finite termination of the proximal point algorithm when applied to linear programs~\cite{Poljak74,Bertsekas75}.

As a consequence of the theorem, the Dantzig selector and noiseless basis pursuit, which are both linear programs, have the exact penalty property. To see why it holds for the Dantzig selector, recast it as the LP 
%To see why the exact penalty property holds for the Dantzig selector,
%recast it as the LP
\begin{equation*} 
  \label{eq:ds-lp2}
  \begin{array}{ll}
    \text{minimize}   & \quad \<\vec{1}, u\> \\
    \text{subject to} & \quad -u \le x \le u\\
    & \quad -\delta \vec{1} \le A^\T (y-Ax) \le \delta \vec{1}, 
  \end{array}
\end{equation*}
with optimization variables $(x,u) \in \R^{2n}$. Then the perturbation
$\thalf \mu \|x-x_0\|^2$ corresponds to a quadratic perturbation with a
diagonal $Q \succeq 0$ obeying $Q_{ii} = 1$ for $1 \le i \le n$ and
$Q_{ii} = 0$ for $n+1 \le i \le 2n$.

\begin{figure} 
 \centering
%\hspace{-7mm}
% Put a border around the image to make it easy to see how much to clip
%\setlength\fboxsep{0pt}
%\fbox{
 \includegraphics[trim = 10mm 2mm 20mm 8mm, clip, width=3.5in]{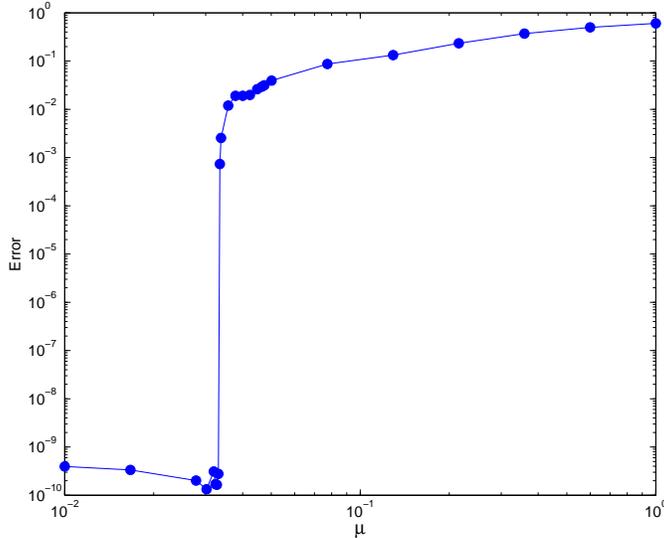}
    %}
%[-16pt]
\\[-12pt]
\caption{\small Demonstration of the exact penalty property for the
  Dantzig selector. For $\mu\le \mu_0\approx 0.025$, the solution to
  the smoothed model coincides with the original model to solver
  precision.}
\label{fig:exactRecovery}
\end{figure}

An illustration of this phenomenon is provided in Figure
\ref{fig:exactRecovery}. For this example, a Dantzig selector
model was constructed for a data set built from a DCT measurement
matrix of size $512\times 4096$. The exact solution $x^\star$ was
constructed
to have 60 dB of dynamic range and 129 nonzeros, using the 
techniques of Appendix~\ref{sec:testProblem}. The figure
plots the smoothing error $\|x^\star-x_\mu\|$ as a function of $\mu$;
below approximately $\mu\approx 0.025$, the error drops rapidly
to solver precision.

% This is new (October 2010):
Unfortunately, compressed sensing models that are not equivalent to linear programs do not in general have the exact penalty property. For instance, it is possible to construct a counter-example for LASSO where the constraint is of the form $\|y-Ax\|\leq\epsilon$. However, the result still suggests that for small $\epsilon$ and $\mu$, the LASSO and its smoothed version are similar.

%It turns out that this \emph{exact penalty} property
%does not hold for all compressed sensing models, but
%it does hold for both the Dantzig selector and noiseless basis
%pursuit. The result was already proved for the special case of noiseless basis pursuit \cite{YinPenalty},
%but it is in fact true for a larger class of problems.

\subsection{Alternative models} % changing from Alternate, as Ewout suggestesd
\label{sec:alternate-ds}

As previously mentioned, different conic formulations result in
different algorithms. To illustrate, consider the first alternative 
(\ref{eq:ds-lp}) proposed in \S\ref{subsec:alternate}, which
represents the Dantzig selector constraint via linear inequalities.
The conic form is
\begin{equation} 
  \label{eq:ds-lp-conic}
  \begin{array}{ll}
    \text{minimize}   & \quad \|x\|_1 \\
    \text{subject to} & \quad \begin{bmatrix} \delta\vec{1}+A^\T (y-Ax) \\ \delta\vec{1}-A^\T (y-Ax) \end{bmatrix} \in\R_+^{2n}. 
  \end{array}
\end{equation}
The dual variable $\lambda\triangleq(\lambda_1,\lambda_2)$ must also
lie in $\R_+^{2n}$. The Lagrangian of the smoothed model is
\[
	\mathcal{L}_\mu(x;\lambda) = \|x\|_1 + \thalf\mu\|x-x_0\|_2^2 - \<\lambda_1,\delta\vec{1}+A^\T(y-Ax)\> - \<\lambda_2,\delta\vec{1}-A^\T(y-Ax)\> \\
\]
and its unique minimizer is given by the soft-thresholding operation
\[
	x(\lambda) = \ST(x_0-\mu^{-1}A^\T A(\lambda_1-\lambda_2),\mu^{-1}).
\]
We cannot eliminate any variables by reducing to composite form,
so we stay with the standard % Ewout's suggestion
smoothed dual function $g_\mu(\lambda)=\inf_x L_\mu(x;\lambda)$, whose gradient is
\[
\nabla g_\mu(\lambda) = \begin{bmatrix}
  -\delta\vec{1}-A^\T(y-Ax(\lambda)) \\
  -\delta\vec{1}+A^\T(y-Ax(\lambda)) \end{bmatrix}.
\]
% EJC: I changed signs below, please check.
The dual update is a true projection
\begin{equation}
	\begin{aligned}
          \lambda_{k+1} &= \argmin_{\lambda\in\R_+^{2n}} -g_\mu(\lambda_k) - \<\nabla g_\mu(\lambda_k),\lambda-\lambda_k\> + \thalf t_k^{-1} \|\lambda-\lambda_k\|_2^2 \\
          &= \argmin_{\lambda\in\R_+^{2n}} \| \lambda - \lambda_k
          - t_k \nabla g_\mu( \lambda_k ) \|_2
	\end{aligned}
\end{equation}
whose solution is simply the non-negative portion of a standard gradient step:
\[
\lambda_{k+1} = \textrm{Pos}( \lambda_k + t_k \nabla
g_\mu(\lambda_k) ), \quad [\textrm{Pos}(z)]_i
\triangleq \begin{cases} z_i, & z_i > 0, \\ 0, & z_i \leq
  0. \end{cases}
\]
In order to better reveal the similarities between the two models, let
us define $\bar{z}\triangleq \lambda_1-\lambda_2$ and
$\bar\tau\triangleq \vec{1}^\T(\lambda_1+\lambda_2)$.  Then we have
$\|\bar{z}\|_1 \leq \bar\tau$, and the Lagrangian and its minimizer
become
\begin{gather*}
	\mathcal{L}_\mu(x;\lambda) = \|x\|_1 + \thalf\mu\|x-x_0\|_2^2 - \<\bar{z},A^\T(y-Ax)\> - \delta\bar{\tau}, \\
	x(\lambda) = \ST(x_0-\mu^{-1}A^\T A\bar{z},\mu^{-1}), 
\end{gather*}
which are actually identical to the original norm-based model. The difference
lies in the dual update. It is possible to show that the dual update
for the original model is equivalent to
\begin{equation}
  z_{k+1} = \Sigma \cdot \argmin_{\lambda\in\R_+^{2n}} -g_\mu(\lambda_k) - \<\nabla g_\mu(\lambda_k),\lambda-\lambda_k\> + \thalf t_k^{-1} \|\Sigma(\lambda-\lambda_k)\|_2^2
\end{equation}
for $\Sigma = [ +I , -I ]$. In short, the dual function is linear in
the directions of $\lambda_1+\lambda_2$, so eliminating them from the
proximity term would yield true numerical equivalence to the original
model.

\section{Further Instantiations}
\label{sec:examples}

Now that we have seen the mechanism for instantiating a particular
instance of a compressed sensing problem, let us show how the
same approach can be applied to several other types of models.
Instead of performing the full, separate derivation for
each case, we first provide a template for our standard
form. Then, for each specific model, we show the necessary
modifications to the template to implement that particular case.

\subsection{A generic algorithm}

A careful examination of our derivations for the Dantzig selector,
as well as the developments in \S\ref{sec:projections}, provide a
clear path to generalizing Listing~\ref{alg:ds-alg} above. We require
implementations of the linear operators $\bar\cA$ and its adjoint
$\bar\cA^\T$, and values of the constants $\bar b$, $b_\tau$; recall
that these are the partitioned versions of $\cA$ and $b$ as described
in \S\ref{sec:composite}.  We also need to be able to perform
the two-sequence recursion (\ref{eq:gptwo}), modified to include
the step size multiplier $\theta_k$ and an adjustable centerpoint $x_0$
in the proximity function.

Armed with these computational primitives, we present in Listing
\ref{alg:generic-alg} a generic equivalent of the algorithm employed
for the Dantzig selector in Listing \ref{alg:ds-alg}. It is important
to note that this particular variant of the optimal first-order
methods may not be the best choice for every model; nevertheless each variant % Ewout's change
uses the same computational primitives.
\begin{algorithm}[H]
\caption{Generic algorithm for the conic standard form}
\label{alg:generic-alg}
\begin{algorithmic}[1]
\REQUIRE $z_0,x_0\in\R^n$, $\mu>0$, step sizes $\{t_k\}$
  \STATE $\theta_0\leftarrow 1$, $v_0\leftarrow z_0$
  \FOR{$k=0,1,2,\dots$}
  \STATE $y_k \leftarrow (1-\theta_k) v_k + \theta_k z_k$
  \STATE $x_k \leftarrow \argmin_x f(x)+ \mu d(x-x_0) - \<\bar\cA^\T(y_k),x\>$
  \STATE $z_{k+1} \leftarrow \argmin_z h(z)+ \tfrac{\theta_k}{2t_k}\|z-z_k\|^2 + \<\bar\cA(x_k)+ \bar b,z\>$
  \STATE $v_{k+1} \leftarrow (1-\theta_k) v_k + \theta_k z_{k+1}$
  \STATE $\theta_{k+1} \leftarrow 2/(1+(1+4/\theta_k^2)^{1/2})$
  \ENDFOR
\end{algorithmic}  
\end{algorithm}
In the next several sections, we show how to construct first-order % Ewout: remove will
methods for a variety of models. We will do so by replacing lines % Ewout: change recreating to replacing
4-5 of Listing~\ref{alg:generic-alg} with appropriate substitutions
and simplified expressions for each.

\subsection{The LASSO}\label{sec:LASSO}

The conic form for the smoothed LASSO is
\begin{equation}
	\begin{array}{ll}
		\text{minimize} & \|x\|_1 + \thalf \mu \| x - x_0 \|_2^2 \\
		\text{subject to} & (y-Ax,\epsilon)\in\mathcal{L}_2^n,
	\end{array}
\end{equation}
where $\mathcal{L}_2^n$ is the epigraph of the Euclidean norm. The Lagrangian is
\[
\lag_\mu(x;z,\tau) = \|x\|_1 + \thalf\mu\|x-x_0\|_2^2 - \langle z, y-Ax \rangle-\epsilon\tau.
\]
The dual variable $\lambda=(z,\tau)$ is constrained to lie in the
dual cone, which is also $\mathcal{L}_2^n$. Eliminating $\tau$ (assuming $\epsilon > 0$) yields
the composite dual form
\[
	\begin{array}{ll}
		\text{maximize} & \inf_x \|x\|_1 + \thalf\mu\|x-x_0\|_2^2 - \langle z, y-Ax \rangle-\epsilon\|z\|_2.
	\end{array}
\]
The primal projection with $f(x)=\|x\|_1$ 
is the same soft-thresholding operation used for
the Dantzig selector. The dual projection involving $h(z)=\epsilon\|z\|_2$,
on the other hand, is
\[
	z_{k+1} = \argmin_z \epsilon\|z\|_2 + \tfrac{\theta_k}{2t_k} \|z-z_k\|_2^2 + \< \tilde{x}, z \> =
	\Shrink(z_k-\theta_k^{-1}t_k \tilde{x},\theta_k^{-1}t_k\epsilon), 
\]
where $\tilde{x}\triangleq y - A x_k$ and $\Shrink$ is an $\ell_2$-shrinkage operation
\[
	\Shrink(z,t) \triangleq \max\{1-t/\|z\|_2,0\} \cdot z = \begin{cases} 0, & \|z\|_2 \leq t, \\ (1-t/\|z\|_2) \cdot z, & \|z\|_2 > t. \end{cases}
\]
The resulting algorithm excerpt is given in Listing~\ref{alg:lasso-alg}.
\begin{algorithm}[H]
\caption{Algorithm excerpt for LASSO}
\label{alg:lasso-alg}
\begin{algorithmic}[1]
	\setcounter{ALC@line}{3}
  \STATE $x_k \leftarrow \ST( x_0 - \mu^{-1} A^\T y_k, \mu^{-1} )$
  \STATE $z_{k+1} \leftarrow \Shrink(z_k-\theta_k^{-1}t_k(y-Ax_k),\theta_k^{-1}t_k\epsilon)$
\end{algorithmic}  
\end{algorithm}

\subsection{Nuclear-norm minimization}
\label{sec:nuclear}

Extending this approach to the nuclear-norm minimization problem
\begin{equation}
\label{eq:nuc_lasso} 
  \begin{array}{ll}
    \text{minimize}   & \quad \|X\|_* \\
    \text{subject to} & \quad \|y-\cA(X)\|_2\leq\epsilon
  \end{array}
\end{equation}
is straightforward. The composite smoothed dual form is
\[
	\begin{array}{ll}
		\text{maximize} & \inf_X \|X\|_* + \mu d(X) - \langle z, y-\cA(X) \rangle-\epsilon\|z\|_2,
	\end{array}
\]
and the dual projection corresponds very directly to the LASSO.
Choosing $d(X) = \thalf \|X - X_0\|_F^2$ leads to a primal projection
given by the soft-thresholding of singular values:  % changing per Ewout's suggestion
\begin{equation} \label{eq:matrix_shrinkage}
X_k = \SVT(X_0 - \mu^{-1} \cA^\T(y_k), \mu^{-1}).
\end{equation}
The $\SVT$ operation obeys
\[
\SVT(X,t) = U \cdot \ST(\Sigma,t) \cdot V^\T, 
\]
where $X = U \Sigma V^\T$ is any singular value decomposition of $Z$,
and $\ST(\Sigma)$
applies soft-thresholding to the singular values (the diagonal
elements of $\Sigma$). This results in the algorithm excerpt presented
in Listing \ref{alg:nn-lasso}.
\begin{algorithm}[H]
\caption{Algorithm for nuclear-norm minimization (LASSO constraint)}
\label{alg:nn-lasso}
\begin{algorithmic}[1]
  \setcounter{ALC@line}{3}
  \STATE $X_k \leftarrow \SVT(X_0 - \mu^{-1} \cA^\T(y_k), \mu^{-1})$
  \STATE $z_{k+1} \leftarrow \Shrink(z_k-\theta_k^{-1}t_k(y - \cA(X_k)),\theta_k^{-1}t_k\epsilon)$
\end{algorithmic}  
\end{algorithm}

Another constraint of interest is of Dantzig-selector type \cite{MatrixOracle} so that one is interested in 
\begin{equation}
\label{eq:nuc_ds} 
  \begin{array}{ll}
    \text{minimize}   & \quad \|X\|_* \\
    \text{subject to} & \quad \|\cA^\T(\cA(X)- y)\| \le \delta.
  \end{array}
\end{equation}
The cone of interest is, therefore, $\cK = \{ (X,t) : \|X\| \le t\}$ and the dual cone is $\cK^* = \{(X,t) : \|X\|_* \le t\}$. The derivation
proceeds as before, and the composite dual problem becomes
\[
	\begin{array}{ll}
		\text{maximize} & \inf_X \|X\|_* + \thalf\mu \|X-X_0\|_F^2 - \langle Z,\cA^\T(y-\cA(X))\rangle-\delta\|Z\|_*.
	\end{array}		
\]
The gradient of the smooth portion has a Lipschitz continuous gradient
with constant at most $\mu^{-1} \|\cA^\T\cA\|^2$, and
singular value thresholding is now used to perform the dual
projection. The resulting excerpt is given in Listing \ref{alg:nn-ds}.
\begin{algorithm}[H]
\caption{Algorithm for nuclear-norm minimization (Dantzig-selector constraint)}
\label{alg:nn-ds}
\begin{algorithmic}[1]
  \setcounter{ALC@line}{3}
  \STATE $X_k \leftarrow \SVT(X_0 - \mu^{-1} \cA^\T(\cA(Y_k)), \mu^{-1})$
  \STATE $Z_{k+1} \leftarrow \SVT(Z_k-\theta_k^{-1}t_k\cA^\T(y - \cA(X_k)),\theta_k^{-1}t_k\delta)$
\end{algorithmic}  
\end{algorithm}

\subsection{\texorpdfstring{$\ell_1$-analysis}{l1-analysis}}
\label{sec:W}

We are now interested in 
\begin{equation}
\label{eq:l1-analysis} 
  \begin{array}{ll}
    \text{minimize}   & \quad \|W x\|_1 \\
    \text{subject to} & \quad \|y-Ax\|_2 \le \epsilon,
    \end{array}
\end{equation}
where the matrix $W$ is arbitrary.  This problem is frequently
discussed in signal processing and is sometimes referred to as the
method of {\em $\ell_1$-analysis}. As explained in the introduction,
this is a challenging problem as stated, because a generalized
projection for $f(x)=\|Wx\|_1$ does not have an analytical form.

Let us apply our techniques to an alternative conic formulation % fixing per Ewout's suggestion
\[
  \begin{array}{ll}
    \text{minimize}   & \quad t \\
    \text{subject to} & \quad \|Wx\|_1 \le t,\\
    & \quad \|y-Ax\|_2 \le \epsilon, 
    \end{array}
\]
where $t$ is a new scalar variable.
The dual variables are $\lambda=(z^{(1)},\tau^{(1)},z^{(2)},\tau^{(2)})$, 
where
\[
\|z^{(1)}\|_\infty \le \tau^{(1)}, \quad \|z^{(2)}\|_2 \le \tau^{(2)},
\]
and the Lagrangian is given by 
\[
\mathcal{L}(x,t;z^{(1)},\tau^{(1)},z^{(2)},\tau^{(2)}) = 
t - \<z^{(1)}, Wx\> - \tau^{(1)} t - \<z^{(2)}, y - Ax\> - \epsilon \tau^{(2)}.
\]
The Lagrangian is unbounded unless $\tau^{(1)} = 1$; and
we can eliminate $\tau^{(2)}$ in our standard fashion as well.
These simplifications yield a dual problem
\[
	\begin{array}{ll}
		\text{maximize}   & \< y, z^{(2)} \> - \epsilon \|z^{(2)}\|_2 \\
		\text{subject to} & A^\T z^{(2)} - W^\T z^{(1)} = 0, \\
		                  & \| z^{(1)} \|_\infty \leq 1. 
	\end{array}		
\]

To apply smoothing to this problem, we use a standard proximity function
$d(x)=\thalf\|x-x_0\|^2$. (Setting $\tau^{(1)}=1$ causes $t$
to be eliminated from the Lagrangian, so it need not appear in our
proximity term.) The dual function becomes
\[
g_\mu(z^{(1)},z^{(2)}) = \inf_x \thalf \mu \|x-x_0\|_2^2 - \<z^{(1)}, Wx\> - \<z^{(2)}, y - Ax\> - \epsilon \| z^{(2)} \|
\]
and the minimizer $x(z)$ is simply
\[
x(z) = x_0 + \mu^{-1} (W^\T z^{(1)} - A^\T z^{(2)}).
\]
Now onto the dual projection
\[
z_{k+1} = \argmin_{z:~\|z^{(2)}\|_\infty\leq 1} \epsilon\|z^{(2)}\|_2 + \tfrac{\theta_k}{2t_k}\|z-z_k\|^2 + \< \tilde{x}, z \>, 
\]
where $\tilde{x}=(Wx(z),y-Ax(z))$. This will certainly converge if the step
sizes $t_k$ are chosen properly.  However, if $W$ and $A$ have
significantly different scaling, the performance of the algorithm may
suffer. Our idea is to apply different step sizes $t_k^{(i)}$ to each
dual variable
\[
z_{k+1} = \argmin_{z:~\|z^{(2)}\|_\infty\leq 1} \epsilon\|z^{(2)}\|_2  + \< \tilde{x}, z \> + \thalf\theta_k \sum_{i=1}^2 (t_k^{(i)})^{-1} \| z^{(i)}-z_k^{(i)} \|_2^2
\]
in a fashion that preserves the convergence properties of the method.
The minimization problem over $z$ is separable, and the solution is given by
\begin{subequations}
\begin{align} 
  z^{(1)}_k & = \Trunc( y^{(1)}_k - \theta_k^{-1} t^{(1)}_k \tilde{x}^{(1)} , \theta_k^{-1} t^{(1)}_k ) \\
  z^{(2)}_k & = \Shrink(y^{(2)}_k - \theta_k^{-1} t^{(2)}_k \tilde{x}^{(2)}, \theta_k^{-1} t^{(2)}_k \epsilon),
\label{update:z2}
\end{align}
\end{subequations}
where the truncation operator is given element-wise by
% adding Ewout's comments above
\[
\Trunc(z,\tau) = \sgn(z) \cdot \min\{|z|,\tau\} = \begin{cases} z,
  \quad & |z| \le \tau, \\ \tau\, \sgn(z), &|z| \ge \tau. \end{cases}
\]
In our current tests, we fix $t_k^{(2)}=\alpha t_k^{(1)}$, where we
choose $\alpha=\|W\|^2/\|A\|^2$, or some estimate thereof. This is 
numerically equivalent to applying a single step size to a scaled
version of the original problem, so convergence guarantees remain.
In future work, however, we intend
to develop a practical method for adapting each
step size separately.

The algorithm excerpt is given in Listing~\ref{alg:l1-analysis}.
\begin{algorithm}[H]
\caption{Algorithm excerpt for $\ell_1$-analysis}
\label{alg:l1-analysis}
\begin{algorithmic}[1]
  \setcounter{ALC@line}{3}
  \STATE $x_k \leftarrow x_0 + \mu^{-1} ( W y_k^{(1)} - A^\T y_k^{(2)} )$
  \STATE \begin{tabular}{@{}l}$z_{k+1}^{(1)}\leftarrow \Trunc(y^{(1)}_{k} - \theta_k^{-1} t^{(1)}_k W x_k,\theta_k^{-1} t^{(1)}_k)$\\$z_{k+1}^{(2)}\leftarrow\Shrink(y^{(2)}_k - \theta_k^{-1} t^{(2)}_k (y-Ax_k), \theta_k^{-1}t^{(2)}_k \epsilon)$\end{tabular}
\end{algorithmic}  
\end{algorithm}

\subsection{Total-variation minimization} \label{sec:TV}

We now wish to solve
\begin{equation}
\label{eq:TV} 
  \begin{array}{ll}
    \text{minimize}   & \|x\|_{\text{TV}} \\
    \text{subject to} & \|y-Ax\|_2 \le \epsilon
    \end{array}
\end{equation}
for some image array $x \in \R^{n^2}$ where $\|x \|_{\text{TV}}$
is the total-variation introduced in \S\ref{sec:motivation}.
We can actually cast this as a \emph{complex} $\ell_1$-analysis problem
\begin{equation*} 
  \begin{array}{ll}
    \text{minimize}   & \|Dx\|_1 \\
    \text{subject to} & \|y-Ax\|_2 \le \epsilon
  \end{array}
\end{equation*}
where $D:\R^{n^2}\rightarrow\mathbb{C}^{(n-1)^2}$ is a matrix 
representing the linear operation that places 
horizontal and vertical differences into the 
real and imaginary elements of the output, respectively:
\begin{gather*}
  [Dx]_{ij} \triangleq (x_{i+1,j} - x_{i,j}) + \sqrt{-1} \cdot
  (x_{i,j+1} - x_{i,j}), \quad 1\leq i<n,\,1\leq j<n.
\end{gather*}
Writing it in this way allows us to adapt our $\ell_1$-analysis derivations directly.
The smoothed dual function becomes
\[
g_\mu(z^{(1)},z^{(2)}) = \inf_x \thalf \mu \|x-x_0\|_2^2 -
\<z^{(1)}, Dx\> - \<z^{(2)}, y - Ax\> -
\epsilon \| z^{(2)} \|_2,
\]
where $z^{(2)}\in\R^m$ is identical to the previous problem, and 
$z^{(1)}\in\mathbb{C}^{(n-1)^2}$ satisfies $\|z^{(1)}\|_\infty\leq 1$.
Supporting a complex $z^{(1)}$ requires two modifications. First, we must
be careful to use the real-valued inner product
\[
	\<z^{(1)}, Dx\>\triangleq\Re((z^{(1)})^HDx) = ( \Re(D^Hz^{(1)}) )^T x
\]	
Second, the projection requires a complex version of the truncation operation:
\[
[\CTrunc(z, \tau)]_k = \min\{1,\tau/|z_k|\} \cdot z_k = \begin{cases} z_k, & |z_k| \leq \tau, \\ \tau z_k / |z_k|, & |z_k| \geq \tau. \end{cases}
\] 

The algorithm excerpt is given in Listing~\ref{alg:tv}.
\begin{algorithm}[H]
\caption{Algorithm excerpt for TV minimization}
\label{alg:tv}
\begin{algorithmic}[1]
  \setcounter{ALC@line}{3}
  \STATE $x_k \leftarrow x_0 + \mu^{-1} ( \Re( D^*y_k^{(1)} ) - A^*y_k^{(2)} )$
  \STATE \begin{tabular}{@{}ll}
  $z_{k+1}^{(1)}\leftarrow \CTrunc(y^{(1)}_{k} - \theta_k^{-1} t_k^{(1)} D x_k,\theta_k^{-1} t_k^{(1)})$ \\
  $z_{k+1}^{(2)}\leftarrow \Shrink(y^{(2)}_k - \theta_k^{-1} t^{(2)}_k (y-Ax_k), \theta_k^{-1} t^{(2)}_k \epsilon)$
  \end{tabular}
\end{algorithmic}  
\end{algorithm}

\subsection{Combining \texorpdfstring{$\ell_1$}{l1} analysis and total-variation
minimization} \label{sec:analysis_tv}

We could multiply our examples indefinitely, and we close this section
by explaining how one could solve the problem \eqref{eq:W+TV}, namely
that of finding the minimum of the weighted combination $\|Wx\|_1 +
\lambda \|x\|_{\text{TV}}$ subject to quadratic constraints. This
problem can be recast as
\begin{equation}  \label{eq:analysis_tv}
  \begin{array}{ll}
    \text{minimize}   & \quad t + \lambda s \\
    \text{subject to}
    & \quad \|Wx\|_1 \le t \\
    & \quad \|Dx\|_1 \le s \\
    & \quad \|Ax - y\|_2 \le \epsilon
  \end{array}
\end{equation}
and the strategy is exactly the same as before. The only difference
with \S\ref{sec:W} and \S\ref{sec:TV} is that the dual variable
now belongs to a direct product of three cones instead of
two. Otherwise, the strategy is the same,
\iffalse
\footnote{If one uses a
  separated quadratic functional to upper bound the smooth part $\gs$,
  one would only need to find a tight majorization of $\nabla^2 \gs$
  by a diagonal matrix.}
\fi
and the path is so clear that we prefer
leaving the details to the reader, who may also want to consult the
user guide which accompanies the software release \cite{Templates}.

\section{Implementing first-order methods}
\label{sec:implement}

So far we have demonstrated how to express compressed sensing problems
in a specific conic form that can be solved using optimal first-order
methods. In this section, we discuss a number of practical matters that
arise in the implementation of optimal first-order methods.  This work
applies to a wider class of models than those presented in this paper;
therefore, we will set aside our conic formulations and present the
first-order algorithms in their more native form.

\subsection{Introduction}

The problems of interest in this paper can be expressed in an unconstrained composite form
\begin{equation}
	\label{eq:composite}
	\begin{array}{ll}
		\text{minimize} & \phi(z) \triangleq g(z) +  h(z),
	\end{array}
\end{equation}
where $g,h:\R^n\rightarrow (\R\cup+\infty)$ are convex functions
with $g$ smooth and $h$ nonsmooth. (To be precise, the dual functions
in our models are concave, so we consider their convex negatives here.) Convex 
constraints are readily supported by including their corresponding indicator functions 
into $h$. 

First-order
methods solve (\ref{eq:composite}) with repeated calls to a
generalized projection, such as
\begin{equation}
	\label{eq:projection}
	z_{k+1} \leftarrow \argmin_z g(z_k) + \langle \nabla g(z_k), z - z_k \rangle + \tfrac{1}{2t_k} \|z-z_k\|^2 + h(z),
\end{equation}
where $\|\cdot\|$ is a chosen norm and $t_k$ is the step size
control. Proofs of global convergence depend upon the right-hand
approximation satisfying an upper bound property
\begin{equation}
	\label{eq:qbound}
	\phi(z_{k+1}) \leq g(z_k) + \langle \nabla g(z_k), z_{k+1} - z_k \rangle + \tfrac{1}{2t_k} \|z_{k+1}-z_k\|^2 + h(z_{k+1}). 
\end{equation}
This bound is certain to hold for sufficiently small $t_k$;
but to ensure global convergence, $t_k$ must be bounded
away from zero. This is typically accomplished by assuming that the
gradient of $g$ satisfies a generalized Lipschitz criterion,
\begin{equation}
	\label{eq:lipschitz2}
	\| \nabla g(x) - \nabla g(y) \|_* \leq L \| x - y \| \quad \forall x,y\in\dom \phi, 
\end{equation}
where $\|\cdot\|_*$ is the dual norm; that is,
$\| w \|_* = \sup \{ \< z, w \> \,|\, \|z\|\leq 1 \}$. Then
the bound (\ref{eq:qbound}) is guaranteed to hold
for any $t_k\leq L^{-1}$. Under these conditions,
convergence to $\epsilon$ accuracy---that
is, $\phi(z_k) - \inf_z \phi(z) \leq \epsilon$---is obtained
in $\order(L/\epsilon)$ iterations
for a simple algorithm based on (\ref{eq:projection}), and $\order(\sqrt{L/\epsilon})$
for the so-called optimal methods \cite{Nesterov83,Nes88,Nesterov07,tseng_2008}.
These optimal methods vary the calculation (\ref{eq:projection}) slightly,
but the structure and complexity remain the same. 

\subsection{The variants}
\label{subsec:variants}

Optimal first-order methods have been a subject of much study in the
last decade by many different authors.  In 2008, Tseng presented a
nearly unified treatment of the most commonly cited methods, and
provided simplified proofs of global convergence and complexity
\cite{tseng_2008}.  The elegance of Tseng's effort seems
underappreciated, possibly due to the fact that his premature passing
delayed formal publication. His work has greatly eased our efforts to
compare the performance of various algorithms applied to our conic
models.

We constructed implementations of five of the optimal
first-order variants as well as a standard projected gradient
algorithm.  To simplify discussion, we have given each variant a 2-3
character label.  Listing~\ref{alg:nest1} depicts N07, a variation of
the method described by Nesterov in \cite{NesterovBook,Nesterov07}.
\begin{algorithm}[H]
\caption{Nesterov's 2007 algorithm (N07).}
\label{alg:nest1}
\begin{algorithmic}[1]
\REQUIRE $z_0\in\dom\phi$, Lipschitz estimate $L$
  \STATE $\bar{z}_0\leftarrow z_0$, $\theta_0\leftarrow 1$
  \FOR{$k=0,1,2,\dots$}
  \STATE $y_k \leftarrow (1-\theta_k) z_k + \theta_k \bar{z}_k$
  \STATE $\bar{z}_{k+1} \leftarrow \argmin_z  \< \theta_k^2 \sum_{i=0}^k \theta_i^{-1} \nabla g(y_i), z \> + \thalf \theta_k^2 L \| z -z_0 \|^2 + h(z)$ 
  \STATE $z_{k+1} \leftarrow \argmin_z \<\nabla g(y_k),z\> + \thalf L \| z - y_k \|^2 + h(z)$
  \STATE $\theta_{k+1} \leftarrow 2 / ( 1 + ( 1 + 4/\theta_k^2 )^{1/2} )$
  \ENDFOR
\end{algorithmic}  
\end{algorithm}
\noindent The other variants can
be described described simply by replacing lines 4-5 as follows.
\begin{itemize}
\item TS: Tseng's single-projection simplification of N07 \cite{tseng_2008}.
\begin{algorithmic}[1]
  \setcounter{ALC@line}{3}
  \STATE $\bar{z}_{k+1} \leftarrow \argmin_z  \< \theta_k^2 \sum_{i=0}^k \theta_i^{-1} \nabla g(y_i), z \> + \thalf \theta_k^2 L \| z -z_0 \|^2 + h(z)$ 
  \STATE $z_{k+1} \leftarrow (1-\theta_k) z_k + \theta_k \bar{z}_{k+1}$
\end{algorithmic}
\item LLM: Lan, Lu, and Monteiro's modification of N07 \cite{LLM09}.
\begin{algorithmic}[1]
  \setcounter{ALC@line}{3}
  \STATE $\bar{z}_{k+1} \leftarrow \argmin_z \< \nabla g(y_k), z \> + \thalf \theta_k L \| z -\bar{z}_k \|^2 + h(z)$ 
  \STATE $z_{k+1} \leftarrow \argmin_z \<\nabla g(y_k),z\> + \thalf L \| z - y_k \|^2 + h(z)$
\end{algorithmic}
\item AT: Auslender and Teboulle's method from \cite{auslender_teboulle_2006}. 
\begin{algorithmic}[1]
  \setcounter{ALC@line}{3}
  \STATE $\bar{z}_{k+1} \leftarrow \argmin_z \< \nabla g(y_k), z \> + \thalf \theta_k L \| z -\bar{z}_k \|^2  + h(z)$ 
  \STATE $z_{k+1} \leftarrow (1-\theta_k) z_k + \theta_k \bar{z}_{k+1}$
\end{algorithmic}  
\item N83: Nesterov's 1983 method \cite{Nesterov83,Nesterov05}; see
  also FISTA \cite{FISTA}.
\begin{algorithmic}[1]
  \setcounter{ALC@line}{3}
	  \STATE $z_{k+1} \leftarrow \argmin_z \<\nabla g(y_k),z\> + \thalf L \| z - y_k \|^2+ h(z)$
	  \STATE Compute $\bar{z}_{k+1}$ to satisfy $z_{k+1} = (1-\theta_k) z_k + \theta_k \bar{z}_{k+1}$.
\end{algorithmic}
\item GRA: The classical projected gradient generalization.
\begin{algorithmic}[1]
  \setcounter{ALC@line}{3}
	  \STATE $z_{k+1} \leftarrow \argmin_z \<\nabla g(y_k),z\> + \thalf L \| z - y_k \|^2 + h(z)$
	  \STATE $\bar{z}_{k+1} \leftarrow z_{k+1}$
\end{algorithmic}  
\end{itemize}
Following Tseng's lead, we have rearranged steps and renamed variables,
compared to their original sources, so that the similarities
are more apparent. 
This does mean that simpler expressions of
some of these algorithms are possible, specifically for
TS, AT, N83, and GRA. Note in particular that GRA does not use the
parameter $\theta_k$.

Given their similar structure, it should not be surprising that
these algorithms, except GRA,
achieve nearly identical theoretical iteration performance. 
Indeed, it can be shown that if $z^\star$ is an optimal point for (\ref{eq:composite}), then
for any of the optimal variants,
\begin{equation}
	\label{ebound}
	\phi(z_{k+1}) - \phi(z^\star) \leq \thalf L\theta_k^2\|z_0-z^\star\|^2  \leq 2Lk^{-2}\|z_0-z^\star\|^2.
\end{equation}
Thus the number of iterations required to reach $\epsilon$ optimality
is at most $\lceil \sqrt{2L/\epsilon}\|z_0-z^\star\|^2\rceil$ (again,
except GRA). Tighter bounds can be constructed in some cases but the
differences remain small.

Despite their obvious similarity, the algorithms do have some key differences
worth noting. First of all, the sequence of points $y_k$ generated by the N83
method may sometimes lie outside of $\dom \phi$. 
This is not an issue for our applications, but it might for those where
$g(z)$ may not be differentiable everywhere.
Secondly, N07 and LLM require two
projections per iteration, while the others require only one. 
Two-projection methods would be preferred only if the added cost
results in a comparable reduction in the number of iterations required.
Theory does not support this trade-off,
but the results in practice may differ; see \S\ref{sec:compare_algos} for
a single comparison.

\subsection{Step size adaptation}
\label{sec:stepsize}

All of the algorithms involve the global Lipschitz constant $L$. Not
only is this constant often difficult or impractical to obtain, 
the step sizes it produces are often too conservative, since the
global Lipschitz
bound (\ref{eq:lipschitz2}) may not be
tight in the neighborhood of the solution trajectory.
Reducing $L$ artificially
can improve performance, but reducing it too much can cause
the algorithms to diverge. Our
experiments suggest that the transition between convergence and divergence
is very sharp. 

A common solution to such issues is backtracking: replace the
global constant $L$ with a per-iteration estimate $L_k$ that is 
increased as local behavior demands it. Examining
the convergence proofs of Tseng reveals that the following 
condition is sufficient to preserve convergence
(see \cite{tseng_2008}, Propositions 1, 2, and 3):
\begin{equation}
	\label{verify2}
	g(z_{k+1}) \leq 
g(y_k) + \< \nabla g(y_k), z_{k+1} - y_k \> + \thalf L_k \| z_{k+1} - y_k \|^2. 
\end{equation}
If we double the value of $L_k$ every time a violation of
(\ref{verify2}) occurs, for instance, then $L_k$ will satisfy $L_k\geq
L$ after no more than $\lceil \log_2(L/L_0)\rceil$ backtracks, after
which the condition must hold for all subsequent iterations. Thus
strict backtracking preserves global convergence. A simple improvement
to this approach is to update $L_k$ with $\max\{2L_k,\hat{L}\}$, where
$\hat{L}$ is the smallest value of $L_k$ that would satisfy
(\ref{verify2}). To determine an initial estimate $L_0$, we can select
any two points $z_0,z_1$ and use the formula
\[
L_0 = \| \nabla g(z_0) - \nabla g(z_1) \|_* / \| z_0 - z_1\|.
\]

Unfortunately, our experiments reveal that (\ref{verify2}) suffers
from severe cancellation errors when $g(z_{k+1})\approx g(y_k)$, often
preventing the algorithms from achieving high levels of accuracy. More
traditional Armijo-style line search tests also suffer from this
issue. We propose an alternative test that maintains fidelity at much
higher levels of accuracy:
\begin{gather}
	\label{verify3}
	\left| \< y_k - z_{k+1}, \nabla g(z_{k+1}) - \nabla g(y_k) \>
        \right| \leq \thalf L_k \| z_{k+1} - y_k \|_2^2.
\end{gather}
It is not difficult to show that (\ref{verify3})
implies (\ref{verify2}), so provable convergence is
maintained. It is a more conservative test, however, producing
smaller step sizes. So for best performance we prefer a hybrid
approach: for instance, use (\ref{verify2}) when
$g(y_k)-g(z_{k+1})\geq \gamma g(z_{k+1})$ for some small $\gamma>0$, and
use (\ref{verify3}) otherwise to avoid the cancellation error issues.

A closer study suggests a further improvement. Because
the error bound (\ref{ebound}) is proportional to $L_k\theta_k^2$, simple
backtracking will cause it to rise unnecessarily. 
This anomaly can be rectified by
modifying $\theta_k$ as well as $L_k$ during a backtracking step.
Such an approach was adopted by Nesterov for N07 in \cite{Nesterov07}; and
with care it can be adapted to any of the variants.
Convergence is preserved if
\begin{equation}
	\label{alphaupdate}
	L_{k+1}\theta_{k+1}^2 /  (1-\theta_{k+1}) \geq L_k\theta_k^2
\end{equation}
(\emph{c.f.} \cite{tseng_2008}, Proposition 1), which implies that the $\theta_k$ update in Line 6 of
Listing~\ref{alg:nest1} should be % not ``satisfy''; Ewout's suggestions
\begin{equation}
  \theta_{k+1} \leftarrow 2 / (1 + (1 + 4 L_{k+1} / \theta_k^2 L_k)^{1/2}).
\end{equation}
With this update the monotonicity of the
error bound (\ref{ebound}) is restored.  For N07 and TS, the update
for $\bar{z}_{k+1}$ must also be modified as follows:
\begin{equation}
  \bar{z}_{k+1} \leftarrow \argmin_z  
  \< \theta_k^2 L_k \textstyle\sum_{i=0}^k (L_i\theta_i)^{-1} \nabla g(y_i), z \> 
  + \thalf \theta_k^2 L_k \| z -z_0 \|^2 + h(z).  
\end{equation}

Finally, to improve performance we may consider \emph{decreasing}
the local Lipschitz estimate $L_k$ when conditions permit.
We have chosen a simple approach: attempt a slight decrease of $L_k$ at each
iteration; that is, $L_k=\alpha L_{k-1}$ for some fixed $\alpha\in(0,1]$. 
Of course, doing so will guarantee that
occasional backtracks occur. With judicious choice of $\alpha$, we can
balance step size growth for limited amounts of backtracking, minimizing
the total number of function evaluations or projections. We have found 
that $\alpha=0.9$ provides good performance in many
applications. 

\subsection{Linear operator structure}
\label{sec:linop}

Let us briefly reconsider the special
structure of our compressed sensing models. In
these problems, it is possible to express the smooth portion of
our composite function in the form 
\[
	g(z) = \bar{g}(\cA^*(z))+\<b,z\>,
\]	
where $\bar{g}$ remains smooth, $\cA$ is a linear operator, and $b$
is a constant vector (see \S\ref{sec:composite}; we have dropped some
overbars here for convenience).
Computing a value of $g$ requires a single application
of $\cA^*$, and computing its gradient also
requires an application
of $\cA$. In many of our models, the functions $\bar{g}$ and $h$
are quite simple to compute,
so the linear operators account for the bulk
of the computational costs. It is to our benefit, then,
to utilize them as efficiently as possible.

For the prototypical algorithms of Section~\ref{subsec:variants}, each
iteration requires the computation of the value and gradient of $g$ at
the single point $y_k$, so all variants require a single application
each of $\cA$ and $\cA^*$. The situation changes when backtracking is
used, however. Specifically, the backtracking test (\ref{verify2})
requires the computation of $g$ at a cost of one application
of $\cA^*$; and the alternative test (\ref{verify3}) requires the
gradient as well, at a cost of one application of $\cA$.  Fortunately,
with a careful rearrangement of computations we can eliminate both of
these additional costs for single-projection methods TS, AT, N83, and
GRA, and one of them for the two-projection methods N07 and LLM. 

\begin{algorithm}
\caption{AT variant with improved backtracking.}
\label{alg:llm-bt}
\begin{algorithmic}[1]
\REQUIRE $z_0\in\dom\phi$, $\hat{L}>0$, $\alpha\in(0,1]$, $\beta\in(0,1)$
  \STATE $\bar{z}_0\leftarrow z_0$, $z_{\cA0},\bar{z}_{\cA0}\leftarrow\bar\cA^*(z_0)$, $\theta_{-1}=+\infty$, $L_{-1}=\hat{L}$
  \FOR{$k=0,1,2,\dots$}
  \STATE $L_k\leftarrow \alpha L_{k-1}$
  \LOOP
  \STATE $\theta_k \leftarrow 2 / ( 1 + ( 1 + 4 L_k / \theta_{k-1}^2 L_{k-1} )^{1/2} )$ ($\theta_0\triangleq 1$)
  \STATE $y_k \leftarrow (1-\theta_k) z_k + \theta_k \bar{z}_k$, $y_{\cA,k} \leftarrow (1-\theta_k) z_{\cA,k} + \theta_k \bar{z}_{\cA,k}$
  \STATE $\bar{g}_k \leftarrow \nabla\bar{g}(y_{\cA,k})$, $g_k \leftarrow \cA \bar{g}_k + b$
  \STATE $\bar{z}_{k+1} \leftarrow \argmin_z \<g_k,z\> + h(z) + \tfrac{\theta_k}{2 t_k} \| z - y_k \|^2$, $\bar{z}_{\cA,k+1}\leftarrow\cA^*(z_{k+1})$
  \STATE $z_{k+1} \leftarrow (1-\theta_k) z_k + \theta_k \bar{z}_{k+1}$, $z_{\cA,k+1} \leftarrow (1-\theta_k) z_{\cA,k} + \theta_k \bar{z}_{\cA,k}$
  \STATE $\hat{L} \leftarrow 2 \left| \< y_{\cA,k} - z_{\cA,k+1}, \nabla\bar{g}(z_{\cA,k+1}) - \bar{g}_k \> \right| / \| z_{k+1} - y_k \|_2^2$
  \STATE \textbf{if} $L_k \geq \hat{L}$ \textbf{then break endif}
  \STATE $L_k \leftarrow \max\{ L_k / \beta, \hat{L} \}$
  \ENDLOOP
  \ENDFOR
\end{algorithmic}  
\end{algorithm}

Listing~\ref{alg:llm-bt} depicts this more efficient use of linear
operators, along with the step size adaptation approach, using the AT
variant. The savings come from two different additions. First, we
maintain additional sequences $z_{\cA,k}$ and $\bar{z}_{\cA,k}$ to
allow us to compute $y_{\cA,k}=\cA^*(y_k)$ without a call to $\cA^*$.
Secondly, we take advantage of the fact that
\[
	\< y_k - z_{k+1}, \nabla g(z_{k+1}) - \nabla g(y_k) \> =
	\< y_{\cA,k} - z_{\cA,k+1}, \nabla \bar{g}(z_{\cA,k+1}) - \nabla \bar{g}(y_{\cA,k}) \>, 
\]
which allows us to avoid having to compute the full gradient of $g$;
instead we need to compute only the significantly less expensive gradient
of $\bar{g}$.

\subsection{Accelerated continuation}
\label{sec:continuation}

Several recent algorithms, such as FPC and NESTA \cite{FPC2,NESTA}, have
empirically found that continuation schemes greatly improve performance.
The idea behind continuation is that we solve the problem of interest by
solving a sequence of similar but easier problems, using the results of
each subproblem to initialize or \emph{warm start} the next one.
Listing~\ref{alg:simpleContinuation} below depicts a standard continuation
loop for solving the generic conic problem in \eqref{eq:stdform_mu} with
a proximity function $d(x) = \frac{1}{2}\|x-x_0\|_2^2$. 
We have used a capital $X$ and loop count $j$ to distinguish these
iterates from the \emph{inner loop} iterates $x_k$ generated by a first-order method.
\begin{algorithm}[H]
    \caption{Standard continuation}
\label{alg:simpleContinuation}
\begin{algorithmic}[1]
\REQUIRE $Y_0$, $\mu_0>0$, $\beta < 1$
  \FOR{$j=0,1,2,\dots$}
  \STATE $X_{j+1} \leftarrow \argmin_{ \cA(x) + b \in \cK } f(x) + 
         \frac{\mu_j}{2}\| x - Y_j \|_2^2 $
  \STATE $Y_{j+1} \leftarrow X_{j+1}$ or $Y_{j+1} \leftarrow Y_j$
  \STATE $\mu_{j+1} \leftarrow \beta \mu_j$
  \ENDFOR
\end{algorithmic}
\end{algorithm}

Note that Listing~\ref{alg:simpleContinuation} allows for the updating of both
the smoothing parameter $\mu_j$ and the proximity center $Y_j$ at each
iteration. In many implementations, $Y_j \equiv X_0$ and only $\mu_j$ is 
decreased, but updating $Y_j$ as well will almost always be
beneficial. When $Y_j$ is updated in this manner, the algorithm 
is known as the proximal point method, which has been studied
since at least \cite{Rockafellar1976}. Indeed, one of the 
accelerated variants \cite{auslender_teboulle_2006} that is 
used in our solvers uses the proximal point framework to
analyze gradient-mapping type updates.
%\CHANGE{It turns out that we can do much better, and the results may have implications for all proximal point methods and related splitting methods \cite{EcksteinBertsekas}.  }
It turns out we can do much better by applying the same
acceleration ideas already mentioned.

Let us suggestively write
\begin{equation}
    h(Y) = \min_{x \in C } f(x) + \frac{\mu}{2}\|x-Y\|_2^2,
    \label{eq:h}
\end{equation}
where $\mu >0$ is fixed and $C$ is a closed convex set.
%\EJC{It may be strange that we return to our
  %conic formulation when we stated at the beginning of the section
  %that we were moving away from it.} \MCG{Agreed. We don't need it. $f(x)$
  %can be extended-valued.}
  This is an infimal convolution,
and $h$ is known as the Moreau-Yosida regularization of
$f$~\cite{HULemarechal}.  Define
\begin{equation}
  X_Y = \argmin_{x \in C} f(x) + \frac{\mu}{2}\|x-Y\|_2^2.
    \label{eq:X_Y}
\end{equation}
The map $ Y \mapsto X_Y$ is a proximity operator~\cite{Moreau1965}.

We now state a very useful theorem.
\begin{theorem}\label{thm:h}
  The function $h$ \eqref{eq:h} is continuously differentiable with
  gradient
\begin{equation}
    \nabla h(Y) = \mu ( Y - X_Y ).
\end{equation}
The gradient is Lipschitz continuous with constant $L = \mu$.
Furthermore, minimizing $h$ is equivalent to minimizing $f(x)$ subject
to ${x \in C}$. 
\end{theorem}

The proof is not difficult and can be found in Proposition I.2.2.4 and
\S XV.4.1 in \cite{HULemarechal}; see also exercises 2.13 and 2.14 in
\cite{BertsekasBook}, where $h$ is referred to as the \emph{envelope}
function of $f$. The Lipschitz constant is $\mu$ since $X_Y$ is a
proximity operator $P$, and $\Id - P$ is non-expansive for any
proximity operator.

The proximal point method can be analyzed in this framework.
%Consider minimizing $h(Y)$.  
%Note that $h(Y) \ge f(x^\star)$ for all
%$Y$, and $h(x^\star) = f(x^\star)$, where $x^\star$ is any optimal
%solution to \eqref{eq:stdform}.  So by solving the unconstrained
%problem $\min_Y h(Y)$, we recover an optimal solution to
%\eqref{eq:stdform}.  
Minimizing $h$ using gradient descent, with
step size $t = 1/L = 1/\mu$, gives
\begin{align*}
  Y_{j+1} = Y_j - t \nabla h(Y_j) &= Y_j - \frac{1}{\mu} \mu(Y_j -
  X_{Y_j} ) = X_{Y_j},
\end{align*}
which is exactly the proximal point algorithm. But since $h$
has a Lipschitz gradient, we can use the accelerated first-order
methods to achieve a
convergence rate of $\order(j^{-2})$, versus
$\order(j^{-1})$ for standard continuation.
In Listing
\ref{alg:accelContinuation}, we offer such an approach using a fixed
step size $t=1/L$ and the accelerated algorithm from Algorithm 2 in \cite{tseng_2008}.
This accelerated version of the proximal point algorithm has been analyzed
in \cite{guler:649} where it was shown to be stable with respect to inexact solves.
\begin{algorithm}[H]
    \caption{Accelerated continuation}
\label{alg:accelContinuation}
\begin{algorithmic}[1]
\REQUIRE $Y_0$, $\mu_0>0$
  \STATE $X_0\leftarrow Y_0$
  \FOR{$j=0,1,2,\dots$}
  \STATE  $ X_{j+1} \leftarrow \argmin_{ \cA(x) + b \in \cK } f(x) + 
         \frac{\mu_j}{2}\| x - Y_j \|_2^2 $
  \STATE $Y_{j+1} \leftarrow X_{j+1} + \frac{j}{j+3}( X_{j+1} - X_j ) $
  \STATE (optional) increase or decrease $\mu_{j}$
  \ENDFOR
\end{algorithmic}
\end{algorithm}

\begin{figure} 
\centering
\hspace{-7mm}
    \includegraphics[trim = 11mm 4mm 22mm 8mm, clip, width=4in]{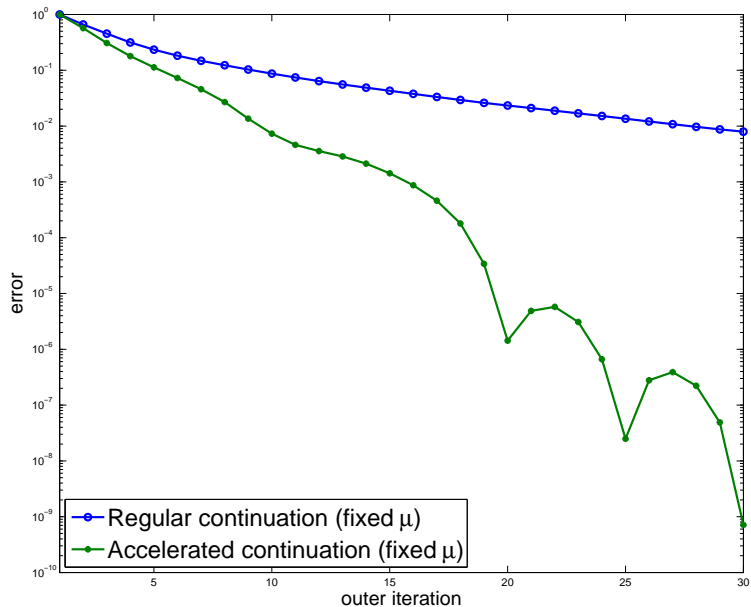}
\\[-10pt]
\caption{A comparison of simple continuation and accelerated
continuation applied to a small-scale LASSO model. The horizontal
axis gives the number of continuation steps taken; the vertical axis gives
the error between the continuation solution and the original, unsmoothed model.
The error is normalized to 1 at iteration 1.}
    \label{fig:continuation_a}
\end{figure} 
Figure~\ref{fig:continuation_a} provides an example of the potential
improvement offered by accelerated continuation. In this case, we have
constructed a LASSO model with a $80 \times 200$
i.i.d.\ Gaussian measurement matrix.  The horizontal axis gives the
number of continuation steps taken, and the vertical axis gives the
error $\|X_j - x^\star\|_2$ between each continuation iterate and the
optimal solution to the unsmoothed model. Both simple and accelerated
continuation have been employed, using a fixed smoothing parameter
$\mu_j\equiv\mu_0$.\footnote{If $\mu$ is fixed,
  this is not \emph{continuation} per se, but we still use the term
  since it refers to an \emph{outer} iteration.}
A clear advantage is demonstrated for accelerated
continuation in this case.

For models that exhibit the exact penalty property, updating the
proximity center $Y_j$ at each iteration yields an interesting result.
Examining the proof in Appendix~\ref{app:exact}, we see that the
property depends not on the size of $\mu$ but rather on the size of
$\mu\|x^\star_j-Y_j\|$, where $x^\star_j$ is projection of $Y_j$ on
the optimal set (\eg, if $x^\star$ is unique, then $x^\star_j\equiv
x^\star$).  Therefore, for \emph{any} positive $\mu$, the exact
penalty property will be obtained if $Y_j$ is sufficiently close to
the optimal set. This has some obvious and very useful consequences.
% raises some interesting questions for future
% research: for instance, can we bound the number of continuation steps
% required to reach the exact penalty region?  If not, are there ways to
% improve the accelerated continuation method so that it becomes
% possible?}

\iffalse
This has a very useful consequence, described in the following theorem:
%This has useful consequences: for instance, if $f$ is a contraction, the continuation
%process will yield the exact solution after a finite
%number of continuation steps.}
\begin{theorem}\label{thm:finite}
    Suppose the problem $\min_{x\in C} f(x)$ is an LP and enjoys the exact penalty property.
    Then if the smoothed subproblem is solved exactly, the outer iteration of a continuation method
    will terminate
    in a finite number of steps.
\end{theorem}
The proof is quite simple and presented in Appendix~\ref{sec:finite}.
\fi

The use of accelerated continuation does require some care, however,
particularly in the dual conic framework.  The issue is that we
solve the dual problem, but continuation is with the primal variable.
When $Y$ is updated in a complicated fashion, as in Listing
\ref{alg:accelContinuation}, or if $\mu$ changes, we no longer have a
good estimate for a starting dual variable $\lambda_j$, so the
subproblem will perhaps take many inner iterations.  In contrast, if
$Y$ is updated in a simple manner and $\mu$ is fixed, as in Listing
\ref{alg:simpleContinuation}, then the old value of the dual variable
is a good initial guess for the new iteration, and this warm start
results in fewer inner iterations.  For this reason, it is sometimes
advantageous to use Listing \ref{alg:simpleContinuation}.

%
%When solving the primal problem directly,
%continuation is a warm-start method because not only is
%$Y_j$ a new prox-center, but it is also a good initial guess for the primal
%solver.  With a dual solver, it is not always possible to find a good initial
%value for the dual variable, since the primal and dual variable are not 
%related in a simple fashion.  A notable exception is for regular
%continuation (Listing \ref{alg:simpleContinuation}) with constant $\mu$. 
%The update is simply $Y_{j+1} = X_j$, and $X_j$ was generated by a dual 
%solution $\lambda_j$.  Hence $\lambda_j$ is a good warm-start value for
%the next outer iteration.  In constrast, the case when $\mu$ changes, or 
%when $Y_{j+1}$ is updated in a more complicated fashion, as in accelerated
%continuation, there is not a good warm-start value for the dual.  In practice,
%$\lambda_j$ is used, but this no longer generates $Y_{j+1}$.
%Thus, the total number of inner iterations in the accelerated continuation
%may be higher than for regular continuation. 

\begin{figure} 
\centering
\hspace{-7mm}
    \includegraphics[trim = 14mm 5mm 20mm 8mm, clip, width=4in]{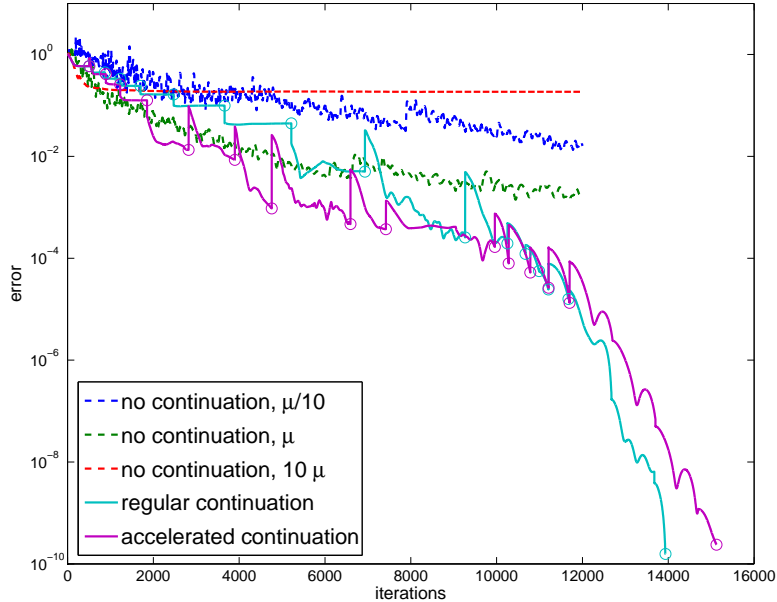}
\\[-10pt]
\caption{Comparing fixed smoothing and continuation strategies for
a Dantzig selector model. The horizontal axis gives the number
of inner iterations of each first-order method. For the two
continuation strategies, the circles depict the completion of
each continuation step.}
\label{fig:continuation_b}
\end{figure} 
Figure \ref{fig:continuation_b} shows an example where accelerated
continuation does not offer a clear benefit. The model solved is the
Dantzig selector, with a $64 \times 256$ partial DCT measurement
matrix, and the exact solution (see Appendix \ref{sec:testProblem})
has $50$ nonzeros and a dynamic range of $68$ dB which makes the test
quite challenging for the minimum $\ell_1$ solution is not the sparsest.
Both fixed values of $\mu$ and both type of continuation are employed.
The horizontal axis gives the total number of \emph{inner} iterations,
and the vertical axis gives the error $\|x_k - x^\star\|_2$ between
the current iterate and the solution to the (unsmoothed) model. When
no continuation is employed, the tradeoff between solution accuracy
(small $\mu$) and the the number of iterations to converge (large
$\mu$) is evident.  Continuation helps drastically, and high % removing ``clearly'' as Ewout suggests
levels of accuracy can be obtained quickly.  In this example, both
forms of continuation perform similarly.  

%When continuation is
%applied, high levels of accuracy are obtained; accelerated continuation
%is not clearly superior to standard continuation. 
%Again, this is
%partly because when $Y_j$ is updated in a complicated manner, the
%conic dual algorithm no longer has a good starting value for the
%dual variable, so the warm start is less effective than for standard
%continuation.

In our experience, %we have observed that 
accelerated 
continuation usually matches or outperforms regular continuation,
but the matter requires further study.  In particular, performance 
is affected by the stopping criteria of the inner iteration, and we plan to 
investigate this in future work.
The experiments in Figure \ref{fig:continuation_b} decreased the tolerance
by a factor of two every iteration, with the exception
of the final iteration which used a stricter tolerance.

\subsection{Strong convexity}
\label{sec:strong_convexity}

Suppose that our composite function $\phi$ is strongly convex; that is,
\[
	\phi(z) \geq \phi(z_0) + \thalf m_\phi \| z - z_0 \|_2^2 \quadŒ\forall z\in\dom\phi
\]
for some fixed $m_\phi>0$.  It is well known that a standard gradient
method will achieve linear convergence in such a case. Unfortunately,
without modification, none of the so-called optimal methods will do
so.  Thus in the presence of strong convexity, standard gradient
methods can actually perform better than their so-called optimal
counterparts.

If the strong convexity parameter $m_\phi\leq L_\phi$ is known or can
be bounded below, the N83 algorithm can be modified to recover linear
convergence \cite{Nesterov05,Nesterov07}, and will recover its
superior performance compared to standard gradient methods.  The
challenge, of course, is that this parameter rarely \emph{is}
known. In \cite{Nesterov07}, Nesterov provides two approaches for
dealing with the case of unknown $m_\phi$, one of which is readily
adapted to all of the first-order variants here. That approach is a
so-called \emph{restart} method: the algorithm is restarted after a
certain number of iterations, using the current iterate as the 
starting point for the restart; or equivalently,
the acceleration parameter $\theta_k$ is reset to $\theta_0=1$. In theory, the
optimal number of iterations between restarts depends on
$m_\phi/L_\phi$, but linear convergence
can be recovered even for sub-optimal choices of this iteration count
\cite{PARNES}.

To illustrate the impact of strong convexity on performance,
we constructed a strongly quadratic unconstrained function with
$m_\phi=0.07$ and $L_\phi=59.1$, and minimized it using the 
gradient method (GRA), the AT first-order method with and without restart,
and the N83 method tuned to the exact values of $(m_\phi,L_\phi)$. Backtracking
was employed in all cases; we have verified experimentally
that doing so consistently improves performance and does not
compromise the exploitation of strong convexity.
As can be seen,
while AT without restart initially performs much better than GRA, the
linear convergence obtained by GRA will eventually overtake it. The
N83 method is clearly superior to either of these, achieving linear
convergence with a much steeper slope. When restart is employed, AT
recovers linear convergence, and achieves near parity with N83 when
restarted every 100 iterations (the optimal value predicted by \cite{PARNES}
is every 112 iterations).
\begin{figure}
\centering
\includegraphics[trim = 15mm 5mm 22mm 11mm, clip, width=4.3in]{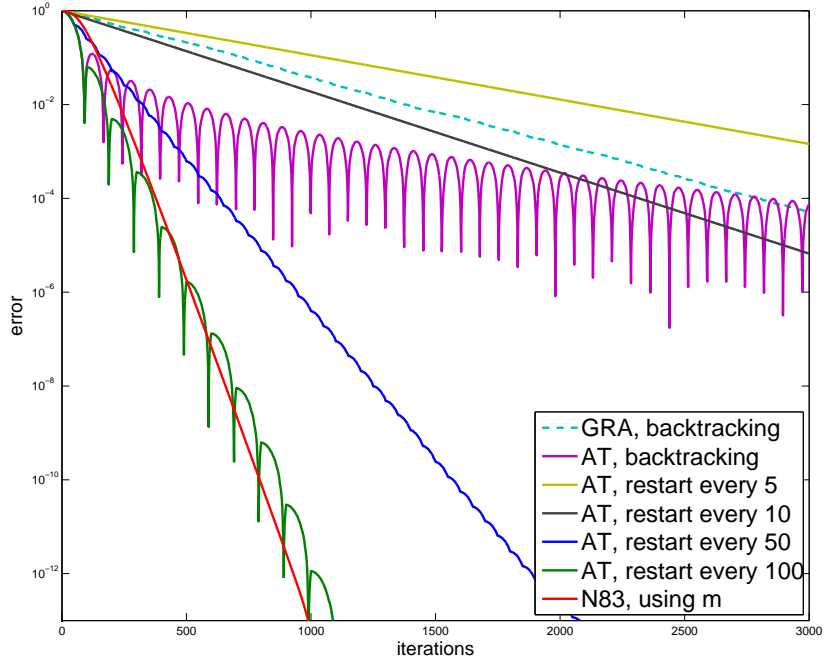}\\[-10pt]
\caption{First-order methods applied to a strongly convex problem: GRA, AT without
restart, AT with restart, and N83 tuned with knowledge of the strong convexity parameter.}
\label{fig:strong_convexity}
\end{figure}

The benefits of exploiting strong convexity seem clearly established.
In fact, while proofs of linear convergence assume global strong
convexity, our experiments in \S\ref{sec:numerics}
show that applying these methods to models with local strong convexity
can improve performance as well, sometimes significantly.
Unfortunately, the restart method requires manual tuning on a
case-by-case basis to achieve the best performance. This is certainly
acceptable in many applications, but further research is needed to
develop approaches that are more automatic.

\section{Numerical experiments}
\label{sec:numerics}

The templates we have introduced offer a flexible framework for
solving many interesting but previously intractable problems; for
example, to our knowledge, there are no first-order algorithms that
can deal with complicated objectives like $f(x) = \|\cWW x\|_1 +
\|x\|_{TV}$ for a non-diagonal and non-orthogonal $W$. This section
shows the templates in use to solve such real-world problems.  It also
describes some of the details behind the numerical experiments in
previous sections.

\subsection{Dantzig selector: comparing first-order variants}
\label{sec:compare_algos}

Other than Tseng's paper \cite{tseng_2008}, there has been little focus
on comparing the various accelerated methods.  Tseng's paper itself
presents few simulations that differentiate the algorithms. Since our
software uses interchangeable solvers with otherwise identical setups,
it is easy to compare the algorithms head-to-head applied to the same model.

% Old stats: 181 nonzeros, .0084 to 6.8 range of coefficients, mu = .27
% New stats: 603 nonzeros, 1.4e-4 to 6.49 range of coefficients (93 dB dynamic range),
%   mu is .284.  D now has 0.7468 dB dynamic range on D (13.1997 smallest, 14.38 largest)

For this comparison, we constructed a smoothed 
Dantzig selector model similar to the one employed in \S\ref{sec:exactRecovery} above.
The model used a partial DCT measurement matrix of size $512\times 2048$, 
a signal with 128 nonzero values, and an additive noise level of 30 dB SNR.
The smoothing parameter was chosen to be $\mu=0.25$, and we then employed the
techniques of Appendix~\ref{sec:testProblem} to perturb the model and obtain
a known exact solution.  This reference solution
had 341 nonzeros, a minimum magnitude of $0.002$ and a
maximum amplitude $8.9$.
The smoothed model was then solved using the 6 first-order
variants discussed here, using both a fixed step size of
$t=1/L=\mu/\|\cA\|^2$ and our proposed backtracking strategy,
as well as a variety of restart intervals.
%   Old text:
%we ensured that the
%exact solution had 181 nonzero entries with magnitudes ranging from
%0.0084 to 6.8. The model was then smoothed with a fixed parameter of
%$\mu=0.27$. 
%We solved this model using each of the 6 first-order
%variants discussed here, using both a fixed step size of
%$t=1/L=\mu/\|\cA\|^2$ and using our proposed backtracking strategy.

\begin{figure}
\centering
    \includegraphics[width=0.49\textwidth]{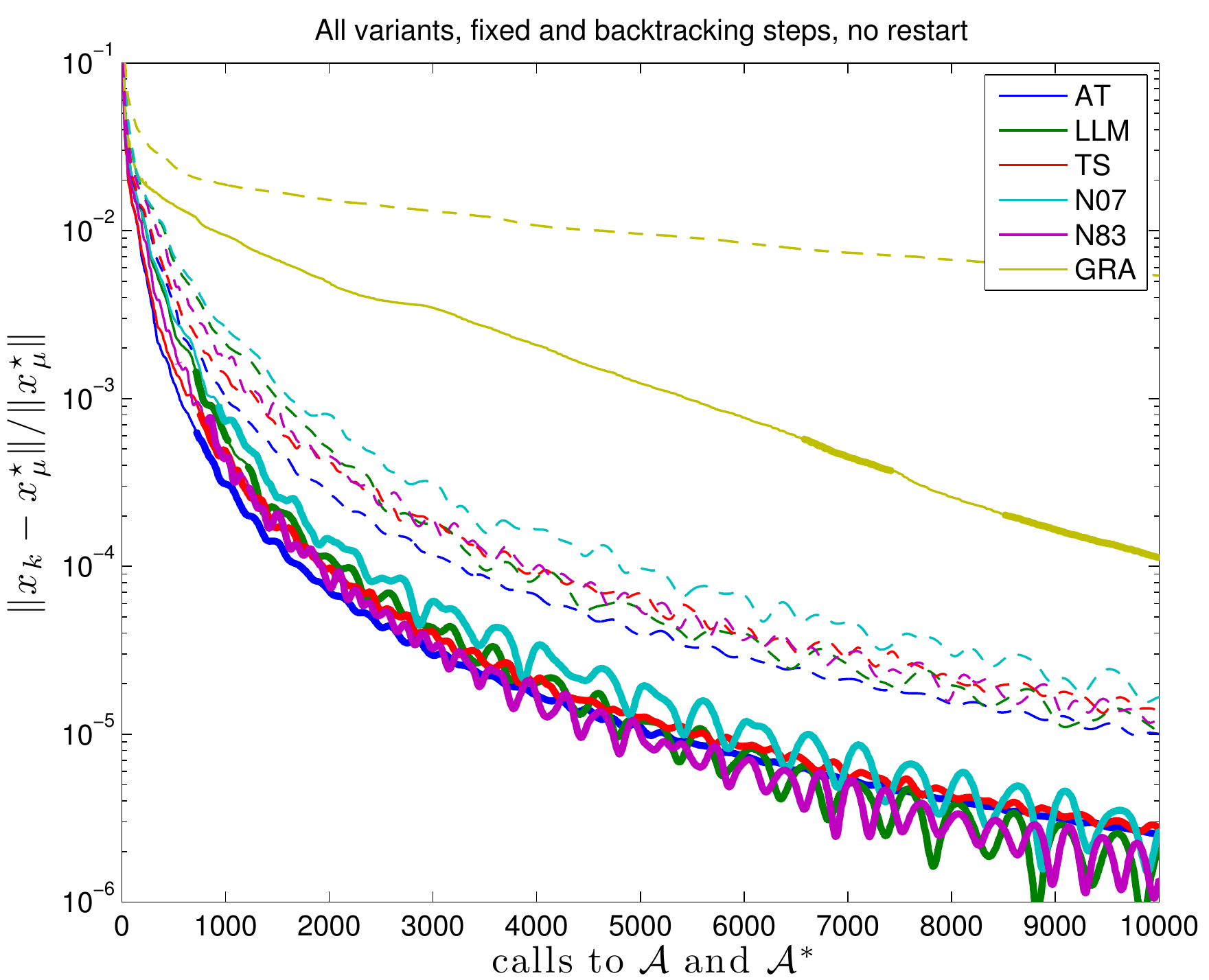}
    \includegraphics[width=0.49\textwidth]{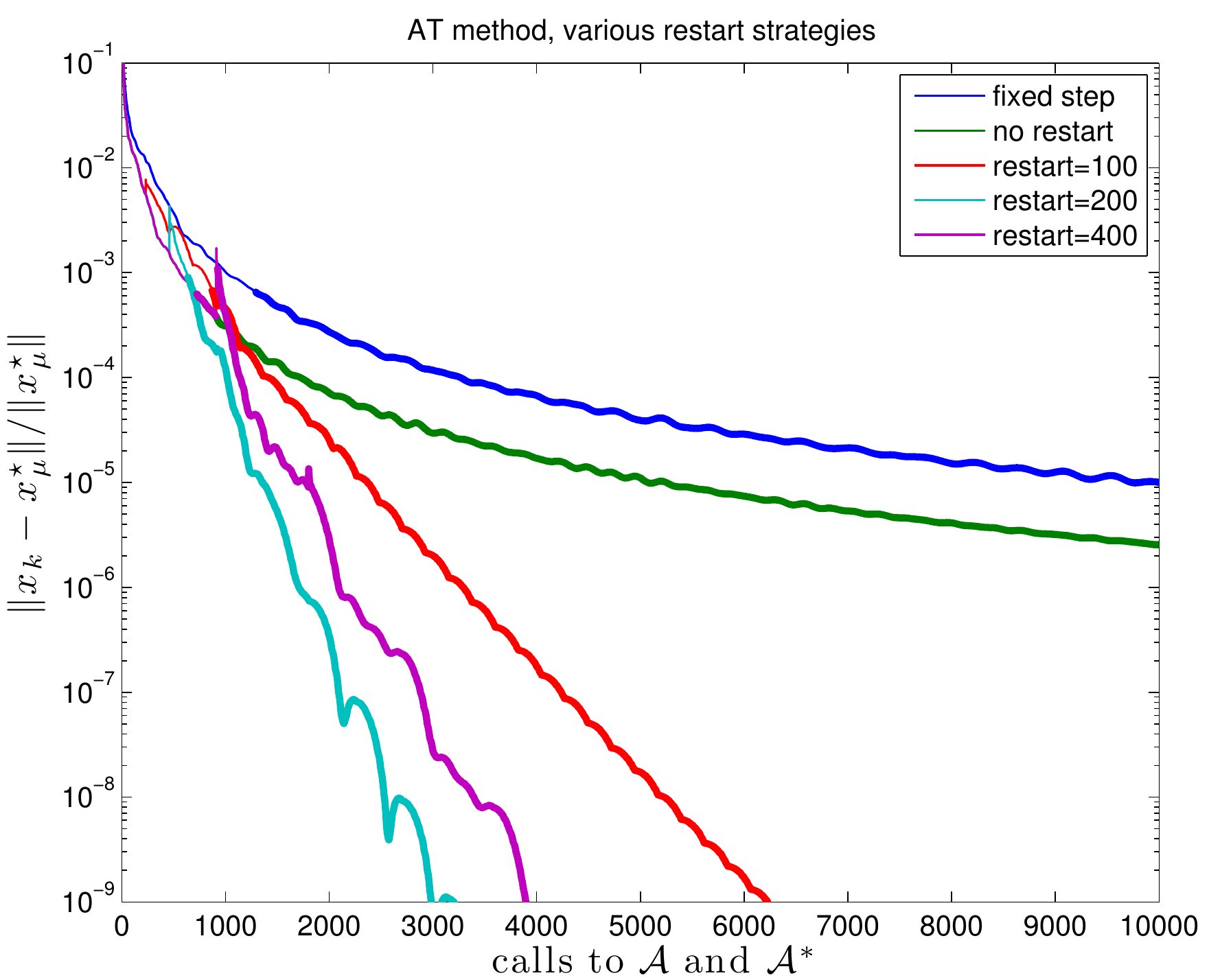}
\caption{Comparing first order methods applied to a smoothed Dantzig selector
model. Left: comparing all variants using a fixed step size
(dashed lines) and backtracking line search (solid lines). Right:
comparing various restart strategies using the AT method.}
\label{fig:compare_algos}
\end{figure}

The results of our tests are summarized by two plots in Figure~\ref{fig:compare_algos}.
The cost of the linear operator dominates, so the horizontal axes
give  the number of calls to either $\cA$ or $\cA^*$
taken by the algorithm. The vertical axes give the
relative error $\|x_k-x^\star_\mu\|/\|x^\star_\mu\|$.  Because this is a sparse recovery
problem, we are also interested in determining when the algorithms find
the correct support; that is, when they correctly identify the locations
of the 341 nonzero entries. Therefore, the lines in each plot are thicker
where the computed support is correct, and thinner when it is not.

The left-hand plot compares
all variants using both fixed step sizes and
backtracking line search, but with no restart.
Not surprisingly, the standard
gradient method performs significantly worse than all of the optimal
first-order methods. 
In the fixed step case, AT performs the
best by a small margin;
but the result is moot, as backtracking shows a significant
performance advantage.  
For example, using the AT variant with a fixed step size requires
more than 3000 calls to $\cA$ or $\cA^*$ to reach an error of $10^{-4}$;
with backtracking, it takes fewer than 2000. With backtracking, 
the algorithms exhibit very similar performance, with AT and TS exhibiting far
less oscillation than the others. All of the methods
except for GRA correctly identify the support (a difficult task
due to the high dynamic range) within 1000 linear operations.

The right-hand plot shows the performance of 
AT if we employ the restart method described in
\S\ref{sec:strong_convexity} for several choices of the restart interval.
We observe significant improvements in performance, revealing
evidence of local strong convexity. A restart
interval of 200 iterations yields the best results; in that case,
a relative error
of $10^{-4}$ is obtained after approximately 1000 linear operations, and
the correct support after only a few hundred operations. The other
variants (except GRA, which is unaffected by restart) 
show similar performance improvements when restart is
applied, although 
the two-projection methods (N07 and LLM) take about 50\% longer than 
the one-projection methods.

Of course, care should be taken when applying these results to other contexts.
For instance, the cost of the projections here is negligible; when they are
are more costly (see, for instance, \S\ref{sec:matrix_completion}),
two-projection methods (N07 and LLM) will be expected to fare worse.
But even among Dantzig selector models, we found significant variations
in performance, depending  upon sparsity, noise level, and smoothing.
% SRB adding
For some models, the two-projection methods perform
well; and in others, such as when we have local strong convexity, gradient
descent performs well (when compared to the other algorithms without
restart).  Overall, it seems there is no \emph{best} algorithm, but we
choose the AT algorithm as our default since in our experience it is
consistently one of the best and only requires one projection per
iteration.

\subsection{LASSO: Comparison with SPGL1}
\label{sec:lasso-spgl1}

As mentioned in the introduction, there are numerous algorithms for
solving the LASSO. Yet the algorithm produced by our dual conic
approach is novel; and despite
its apparent simplicity, it is competitive with the state of
the art. To show this, we compared the AT first-order variant with
SPGL1 \cite{vandenberg_friedlander_2009}, chosen
because recent and extensive tests in \cite{NESTA} suggest that it is
one of the best available methods.

The nature of the SPGL1 algorithm, which solves a sequence
of related problems in a root-finding-scheme, is such that it is fastest
when the noise parameter $\epsilon$ is large, and slowest when $\epsilon = 0$. 
To compare performance in both regimes, we constructed two tests. The
first is an ``academic''  test with an $s$-sparse signal and no
noise; however, we choose $s$ large enough so that the LASSO solution
does not coincide with the sparse solution, since empirically
this is more challenging for solvers. Specifically, the measurement
matrix $A$ is a $2^{13} \times 2^{14}$ partial DCT, while
the optimal value $x^\star$ was constructed to have $s=2^{12}$ nonzeros.
The second test uses Haar wavelet coefficients from the 
``cameraman'' test image (Figure \ref{fig:tv_analysis_test} (a))
which decay roughly according to a power law, and adds noise 
with a signal-to-noise ratio of 30 dB. 
The measurement matrix is also a partial DCT, this time of size 
$ 0.3\cdot 2^{16} \times 2^{16}$.

\begin{figure}
\centering
\includegraphics[width=0.49\textwidth,trim= 10mm 1mm 26mm 14mm,clip]{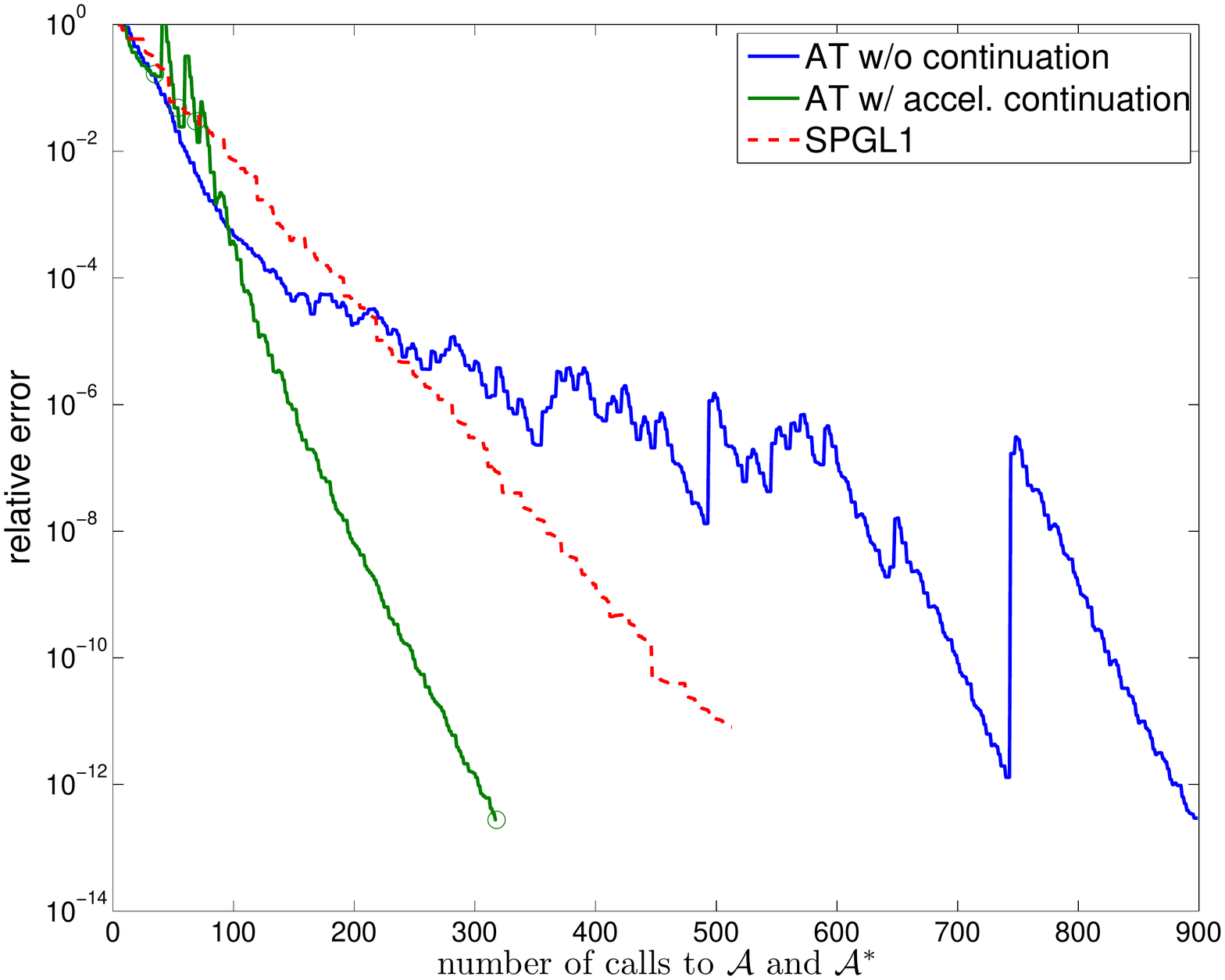}
\includegraphics[width=0.49\textwidth,trim= 10mm 1mm 26mm 14mm,clip]{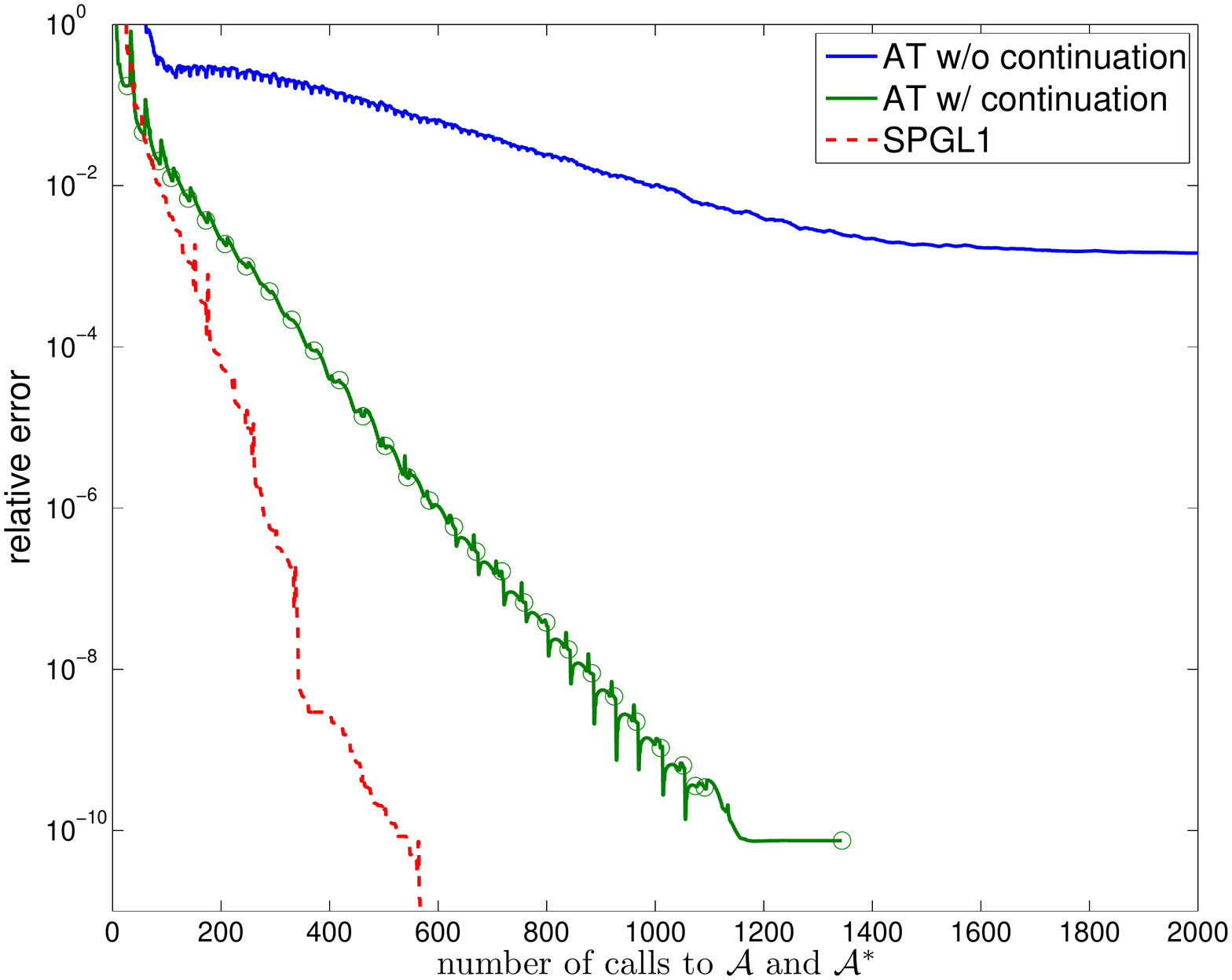} \\
\caption{Comparisons of the dual solver with SPGL1. The plot on the left involves
a noiseless basis pursuit model, while the plot on the right represents a noisy image model.}
\label{fig:spgl1_test}
% Notes: figure on right is from ``test 5'' and shows 2^14 randomly permuted wavelet coefficients (Haar wavelets, 3 levels), with DCT measurements, 30 dB SNR,
% image is cameraman.tif (builtin to Matlab), M = .3*N
% Figure on left is ``test 6'' and is noiseless, original signal has 10 dB dynamic range, N = 2^10 measurements, M = N/2, K = M/4, Gaussian measurement matrix
\end{figure}
Figure \ref{fig:spgl1_test} shows the results from both tests, each
plot depicting relative error $\|x_k-x^\star\|/\|x^\star\|$ versus the
number of linear operations. We see that both methods achieve several
digits of accuracy in just a few hundred applications of $A$ and its
adjoint.  SPGL1 outperforms a regular AT solver in the left-hand
``academic'' test; however, AT with accelerated continuation solves
the problem significantly faster than SPGL1.  The noiseless case
exploits our method's strength since the exact penalty property holds.

For the wavelet test, SPGL1 outperforms our method, even when we use
continuation. Although AT with continuation achieves high levels of
accuracy in fewer than 1000 operations, other tests confirmed that
SPGL1 is often a little better than our method, especially
for large $\epsilon$. But the dual conic approach is competitive in
many cases, and can be applied to a wider class of
problems\footnote{SPGL1 is also flexible in the choice of norm, but
  notably, it cannot solve the analysis problem due to difficulties in
  the primal projection.}.

\subsection{Wavelet analysis with total-variation} 
\label{sec:numerics_tv_wavelet}

\begin{figure}
\centering
\begin{tabular}{c c}
    \includegraphics[width=2.2in]{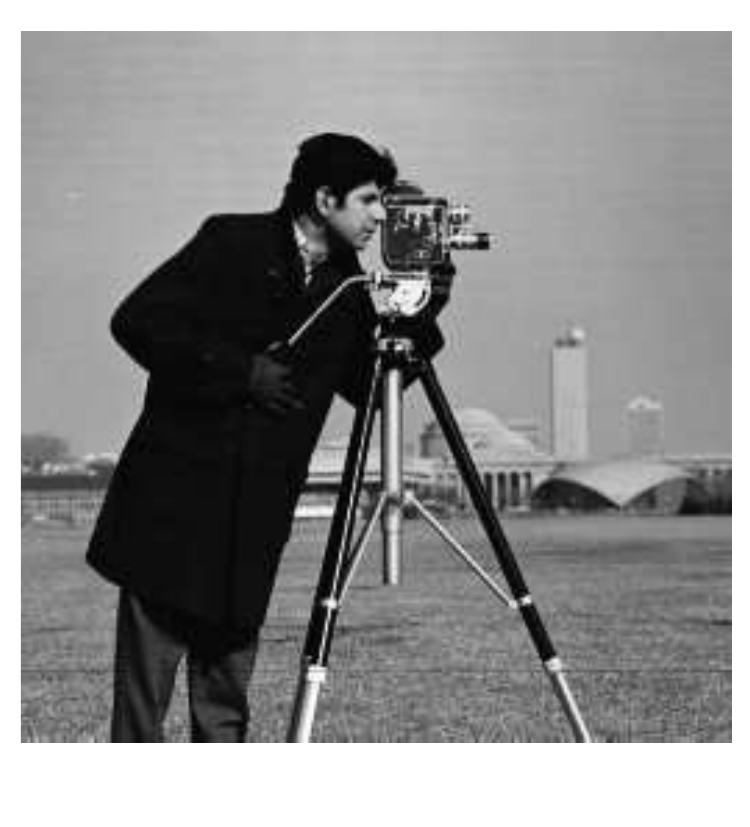}
    &
    \includegraphics[width=2.2in]{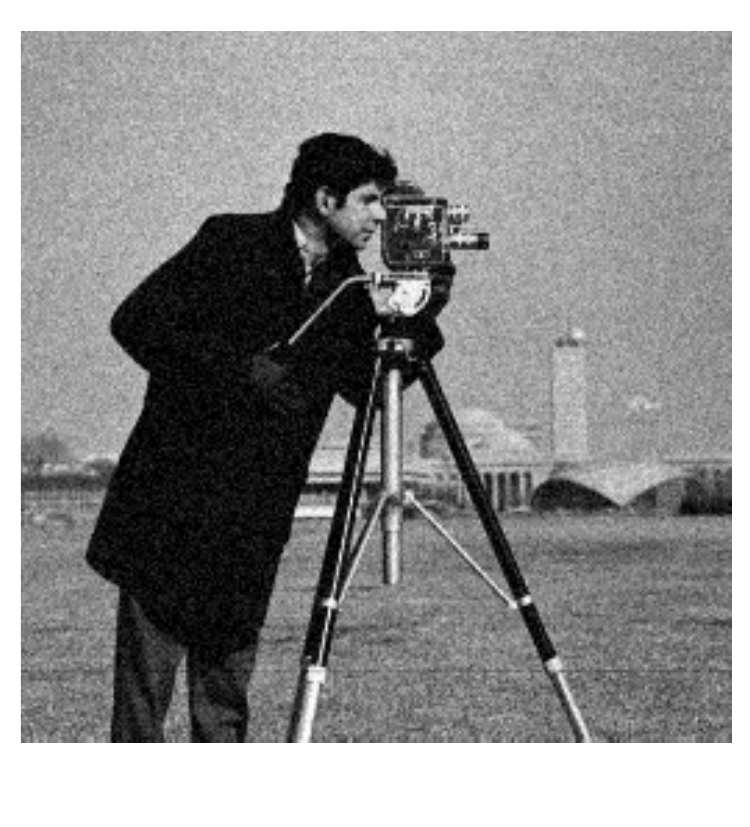} \\
[-16pt]
\small  (a) Original & \small (b) Noisy version (25.6 dB PSNR) \\
    \includegraphics[width=2.2in]{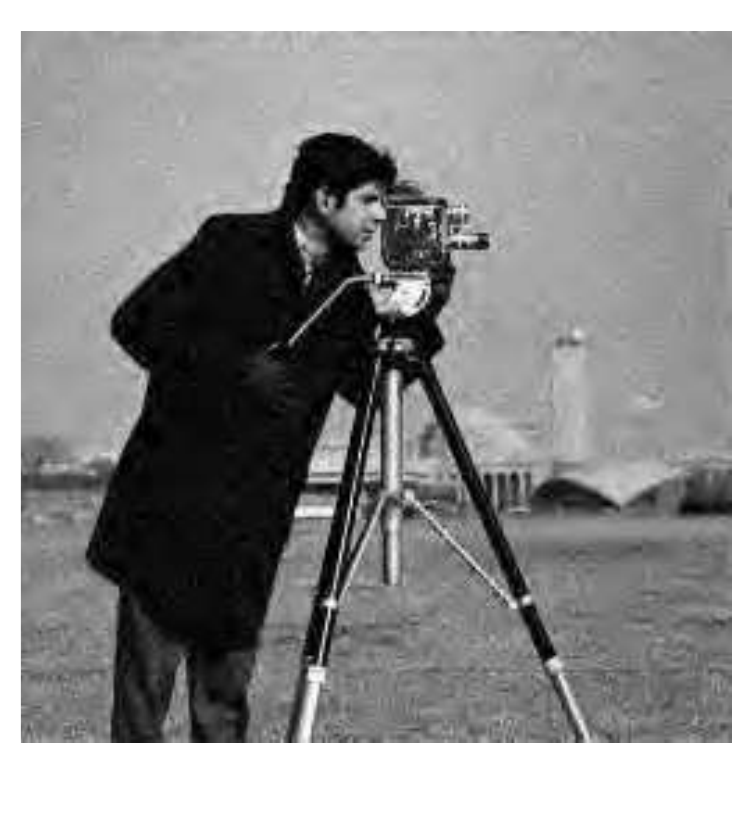}   &   \includegraphics[width=2.2in]{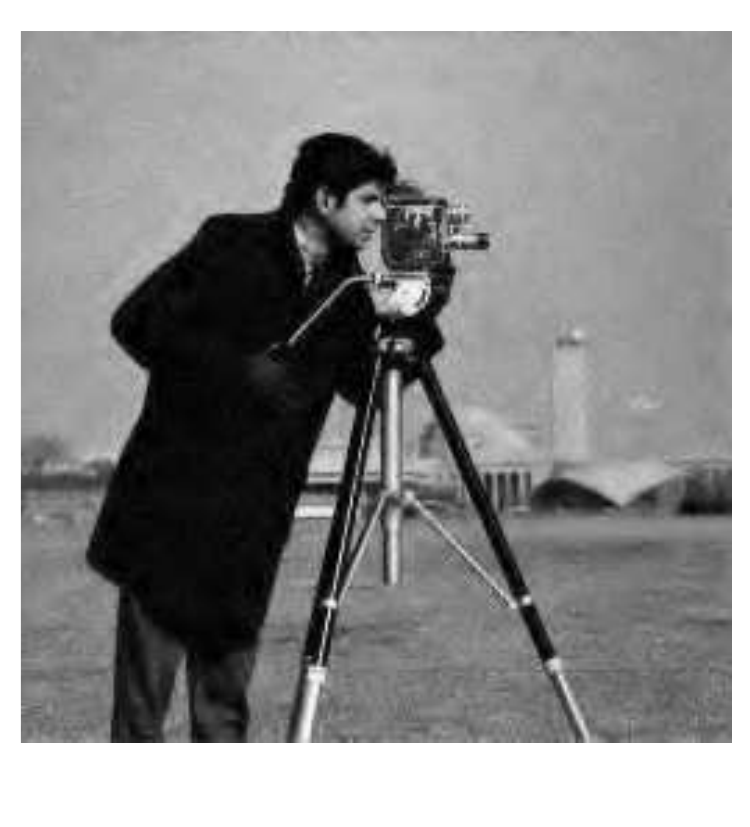} \\
[-16pt]
\small (c) Wavelet thresholded (28.3 dB PSNR)& \small (d) Wavelet $\ell_1$-analysis (29.0 dB PSNR) \\
    \includegraphics[width=2.2in]{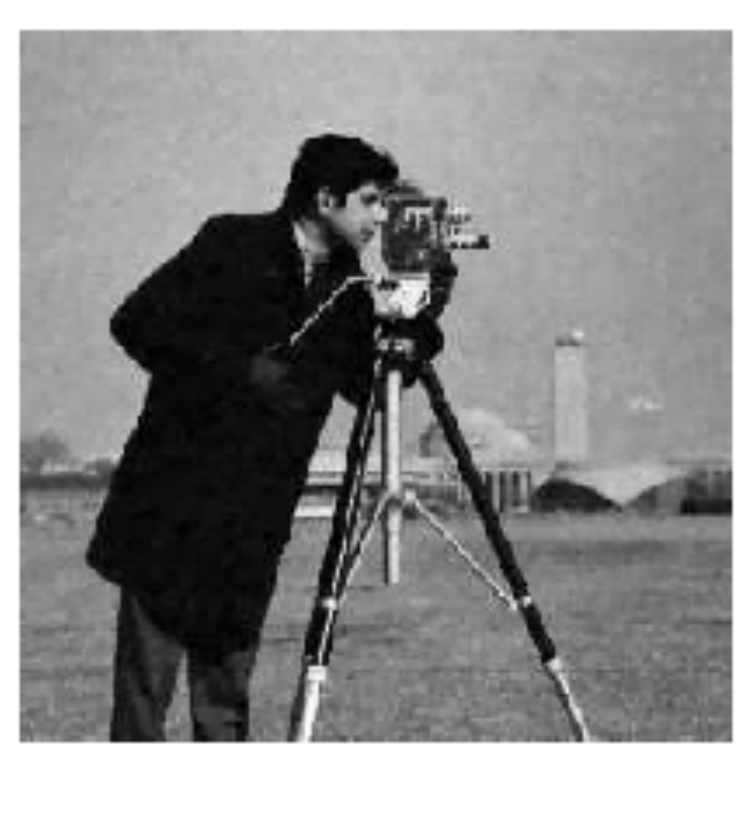}   &   \includegraphics[width=2.2in]{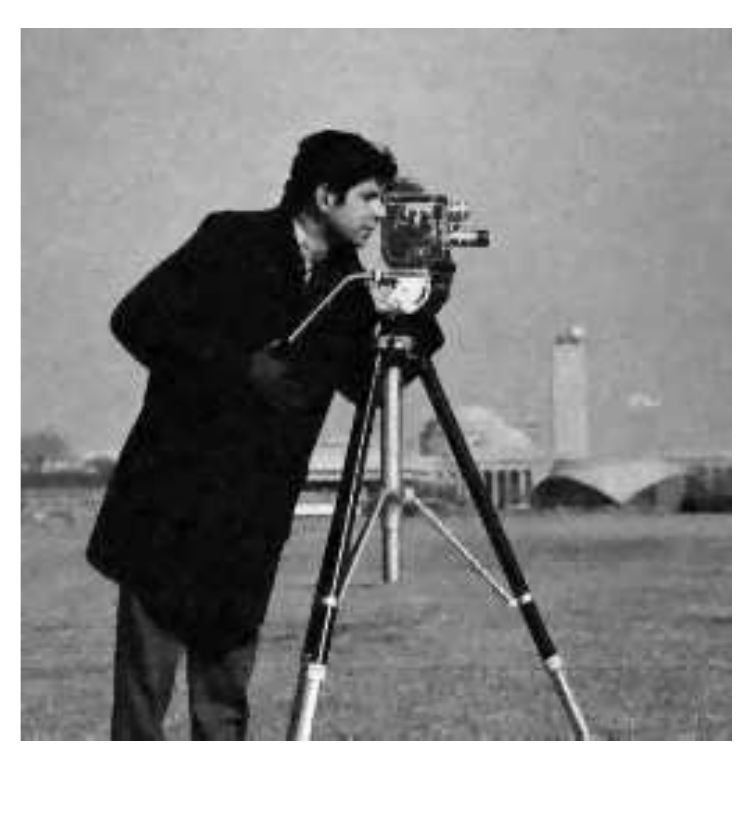} \\
[-16pt]
\small (e) TV minimization (30.9 dB PSNR) &\small (f) Wavelet analysis + TV minimization (31.0 dB PSNR) \\
\end{tabular}
\\[-10pt]
\caption{Denoising an $n=256^2$ image.}
\label{fig:tv_analysis_test}
\end{figure}

The benefit of our approach is highlighted by the fact
that we can solve complicated composite objective functions.  Using
the solver templates, it is easy to solve the $\ell_1$-analysis and TV
problem from \S\ref{sec:analysis_tv}.  We consider here a denoising
problem with full observations; \ie, $\cAA = \Id$.  Figure
\ref{fig:tv_analysis_test} (a) shows the original image $x_0$, to
which noise is added to give an image $y = x_0 + z$ with a
signal-to-noise ratio of 20 dB (see subplot (b)).  In the figure,
error is measured in peak-signal-to-noise ratio (PSNR), which for an
$n_1 \times n_2$ image $x$ with pixel values between in $[0,1]$ is
defined as
$$ \text{PSNR}(x) = 20\log_{10}\left( \frac{ \sqrt{n_1 n_2} }{\|x-x_0\|_F} \right) $$
where $x_0$ is the noiseless image.

To denoise, we work with a $9/7$ bi-orthogonal wavelet transform
$\cWW$ (similar to that in JPEG-2000) with periodic boundary
conditions (the periodicity is not ideal to achieve the lowest
distortion).  A simple denoising approach is to hard-threshold the
wavelet coefficients $\cWW x$ and then invert with $\cWW^{-1}$. Figure
\ref{fig:tv_analysis_test} (c) shows the result, where the
hard-threshold parameter was determined experimentally to give the
best PSNR.  We refer to this as ``oracle thresholding'' since we used
the knowledge of $x_0$ to determine the threshold parameter. As is
common with wavelet methods \cite{Starck03}, edges in the figure
induce artifacts.

Figures \ref{fig:tv_analysis_test} (d), (e) and (f) are produced
using the solvers in this paper, solving 
\begin{equation} \label{eq:weightedTV}
  \begin{array}{ll}
      \text{minimize}   & \alpha \|\cWW x\|_1 + \beta \|x\|_{\text{TV}} + \frac{\mu}{2}\|x-y\|_F^2 \\
    \text{subject to} & \|\cAA x  - y\|_2 \le \epsilon
    \end{array}
\end{equation}
Figure \ref{fig:tv_analysis_test} (d) employs wavelet analysis
only ($\alpha=1,\beta = 0$), (e) just TV ($\alpha = 0,\beta=1$), and (f)
both ($\alpha = 1,\beta=5$).  
For the best performance, the matrix $A$ was re-scaled so that all the dual variables are
of the same order; see \cite{Templates} for further discussion of scaling.

The Frobenius term in (\ref{eq:weightedTV}) is of course for smoothing purposes, and it is
possible to minimize its effect by choosing $\mu$ small or 
using the continuation techniques discussed. But
for denoising, its presence makes little
difference; in fact, it may give more visually pleasing results to
use a relatively large $\mu$. So to determine $\mu$, we started with an estimate like 
\[
\mu = \max( \alpha \|\cWW y\|_1, \beta \|y\|_{TV} ) / c 
\]
with $c \simeq 500$ and then adjusted to give reasonable results.
We ultimately employed 
$\mu = 1$ for (d), $\mu=50$ for (e), and $\mu=160$ for (f).

%For just wavelets, or just TV, as
%in \S\ref{sec:W} and \S\ref{sec:TV}, the Lipschitz constant is
%bounded by a scalar constant times the spectral norm of a block
%diagonal matrix, with blocks $\cAA \cAA^\T $ and either $\cWW \cWW^\T$ (for
%analysis) or $\cDD \cDD^\T$ (for TV).  It is desirable that all dual
%variables have the same scale, since there is a single scalar step
%size, and hence we can scale $\cAA$ (and $\bb$ and $\epsilon$) so that
%$\|\cAA\cAA^\T\| = \|\cWW \cWW^\T\| $ (or with $\cDD$ for TV).  For the
%composite problem, with $\alpha,\beta > 0$, the bound on the Hessian
%is block diagonal with three blocks: $\cAA \cAA^\T$, $\cWW \cWW^\T$ and
%$\cDD\cDD^\T$.  It is in general impossible to rescale to make all three
%norms equal (since $\frac{\alpha}{\beta}$ should remain fixed), but we
%can still scale $\cAA$ so that $\| \cAA\cAA^\T \| = \max( \alpha^2 \|\cWW
%\cWW^\T\|, \beta^2 \|\cDD \cDD^\T \| )$. 

The wavelet analysis run took 26 iterations, and was complete in about
5 seconds. As shown in image (d), it produced boundary effects that are
noticeable to the eye, and similar to those produced by thresholded image (c).
The TV model (e) and TV with wavelets model (f) took 37 iterations (3 seconds)
and 30 iterations (8 seconds), respectively. Both produced
better reconstructions, both by PSNR and by visual inspection.
The additional wavelet analysis term in plot (f) offers only minimal
improvement over TV alone, but this may be due to our simple choice of
wavelet transform.  For example, undecimated wavelets are common in
denoising and may give better results, but our point here is
simplicity and to point out the flexibility of the framework.

\subsection{Matrix completion: expensive projections}
\label{sec:matrix_completion}

We consider the nuclear-norm minimization problem \eqref{eq:nuc_lasso}
of a matrix $X \in \R^{n_1 \times n_2}$ in the conic dual smoothing
approach.  For matrix completion, the linear operator $\cA$ is the
subsampling operator revealing entries in some subset $E \subset [n_1]
\times [n_2]$. With equality constraints ($\epsilon = 0$) and $X_0=
0$, gradient ascent on the dual is equivalent to the SVT algorithm of
\cite{SVT}, a reference which also considered non-equality
constraints, \eg, of the form \eqref{eq:nuc_lasso}.

In addition to our own interest in this problem, one of the reasons 
we chose it for this article is that it differs from the others
in one key respect: its computational cost is dominated by
one of the projections, not by the linear operators. After all,
the linear operator in this case is no more than 
a set of memory accesses, while the primal projection requires the
computation of at least the largest singular values of a large matrix.
As a result, the considerations we bring to the design of an
efficient solver are unique in comparison to the other examples presented.

There are a number of strategies we can employ to reduce the cost of this
computation. The key is to exploit the fact that a nuclear-norm matrix completion
model is primarily of interest when its optimal value $X^\star$ is expected to have
low rank \cite{CR08}. Recall from \S\ref{sec:nuclear} that the update
of $X_k$ takes the form
\[
X_k = \SVT\bigl(X_0 - \mu^{-1} \cA^\T(\dual), \mu^{-1} \bigr).
\]
Our numerical experiments show that if
$\mu$ is sufficiently small, the ranks $r_k=\rank(X_k)$ will
remain within the neighborhood of $r^\star=\rank(X^\star)$. In fact,
with $\mu$ sufficiently small and $X_0=0$, the rank grows monotonically. 

There are a variety of ways we can exploit the low-rank structure
of the iterates $X_k$. By storing $X_k=U_k\Sigma_k V_k$
in factored form, we can reduce the storage costs from $\order(n_1n_2)$
to $\order(r_k(n_1+n_2))$. The quantity $\cA(X_k)$, used in the computation
of the gradient of the dual function, can be computed efficiently
from this factored form as well. Finally, using a Lanczos method
such as PROPACK \cite{propack}, the cost of computing the necessary
singular values will be roughly proportional to $r_k$. By combining
these techniques, the overall result is that the cost of singular
value thresholding is roughly linear in the rank of its result. These
techniques are discussed in more detail in \cite{Templates}.

Smaller values of $\mu$, therefore, reduce the ranks of the iterates
$X_k$, thereby reducing the computational cost of each iteration.
This leads to a unique tradeoff, however: as we already know, smaller values of $\mu$
increase the number of iterations required for convergence. In practice,
we have found that it is indeed best to choose a small value of $\mu$,
and not to employ continuation.

For our numerical example, we constructed a rank-10 matrix of size
$n_1=n_2=1000$, and randomly selected $10\%$ of the entries for
measurement---about 5 times the number of degrees of freedom in the
original matrix.  We solved this problem using three variations of GRA
and five variations of AT, varying the step size choices and the
restart parameter.  The smoothing parameter was chosen to be $\mu =
10^{-4}$, for which all of the methods yield a monotonic increase in
the rank of the primal variable.  The pure gradient method took about
$1.8$ minutes to reach $10^{-4}$ relative error, while the AT method
without restart took twice the time;  % Ewout suggests eliminating this: to reach $10^{-4}$; 
the error is
in the Frobenius norm, comparing against the true low-rank matrix,
which is the solution to the unperturbed problem.
%The pure gradient methods took
%between 4 and 6 minutes to complete 500 iterations

\begin{figure}
\centering
\includegraphics[trim = 15mm 5mm 22mm 11mm, clip, width=0.6\textwidth]{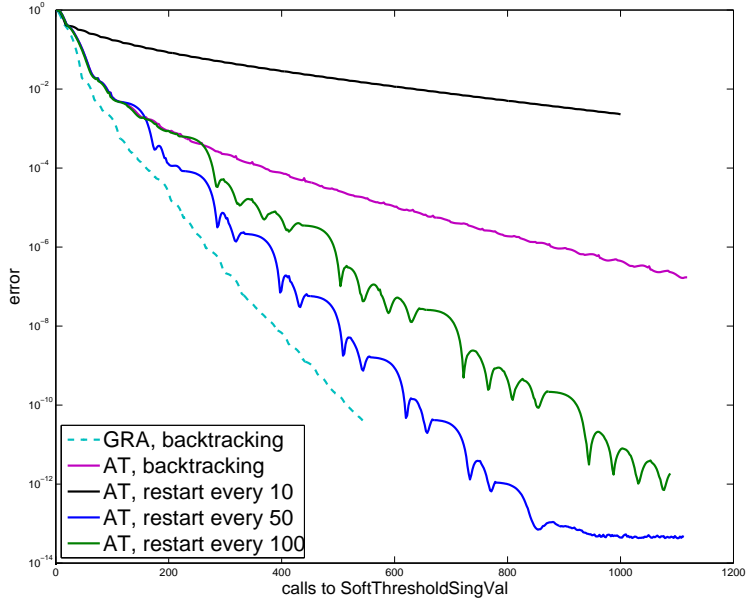}\\[-10pt]
\caption{Noiseless matrix completion using various first-order methods.}
\label{fig:nuclear_strong_convexity}
\end{figure}
The results of our experiments are summarized in Figure~\ref{fig:nuclear_strong_convexity}.
The horizontal axis now gives the number of $\SVT$ operations, a more accurate
measure of the cost in this case; the cost of a $\SVT$ call is not fixed,
but we observe that the rank of the iterates quickly reaches 10 and most $\SVT$ calls
have roughly the same cost.
The salient feature of this figure is that while AT initially outperforms GRA, the linear
convergence exhibited by GRA allows it to overtake AT at moderate levels
of precision. Employing a restart method with AT improves its performance
significantly; and for a restart intervals of 50 iterations, its performance
approaches that of GRA, but it does not overtake it.

What accounts for this performance reversal? First, we observe that
AT requires two projections per iteration while GRA requires only one.
This difference is inconsequential
in our other examples, but not for this one; and it is due to our use
of the backtracking line search.
Switching to a fixed step size eliminated the extra projection, but the overall 
performance suffered significantly. We hope to identify a new line
search approach that avoids the added cost revealed here.

Second, the similarities of Figures \ref{fig:strong_convexity} and
\ref{fig:nuclear_strong_convexity} suggest the presence of
strong convexity. Using an estimate of
the decay rate for gradient descent with step size $t=1/L_f$, and
comparing with known decay rate estimates \cite{NesterovBook} gives an
estimate of $m_f = 0.0024$.  For this value of $m_f$ (and with
$L_f=1$), the optimal restart number $K_{opt}$ from \cite{PARNES} is
about 80, which is consistent with the plot. 
It is easy to verify, however, that the smoothed dual function is \emph{not}
strongly convex.

The results suggest, then, that \emph{local} strong convexity is present.
The authors in \cite{PARNES} argue that for compressed sensing
problems whose measurement matrices satisfy the restricted isometry property, the primal
objective is locally strongly convex when the primal variable is
sufficiently sparse. The same reasoning may apply to the matrix
completion problem; the operator $\cA^\T \cA$ is nearly isometric
(up to a scaling) when restricted to the set of low-rank matrices.
The effect may be amplified by the fact that we have
taken pains to ensure that the iterates remain low-rank.
Overall, this effect is not yet well understood,
but will also be explored in later work.

We note that if the underlying matrix is not low-rank,
then the local strong
convexity effect is not present, and gradient descent does not
outperform the AT method.  
In this case, our experiments also
suggest that the restart method has little effect. Likewise,
when inequality constraints are added we observe that the unusually
good performance of gradient descent vanishes.
For both the non-low-rank and noisy cases,
then, it may be beneficial to employ an optimal first-order
method and to use continuation.

\begin{figure}
    \centering
    \includegraphics[width=0.6\textwidth]{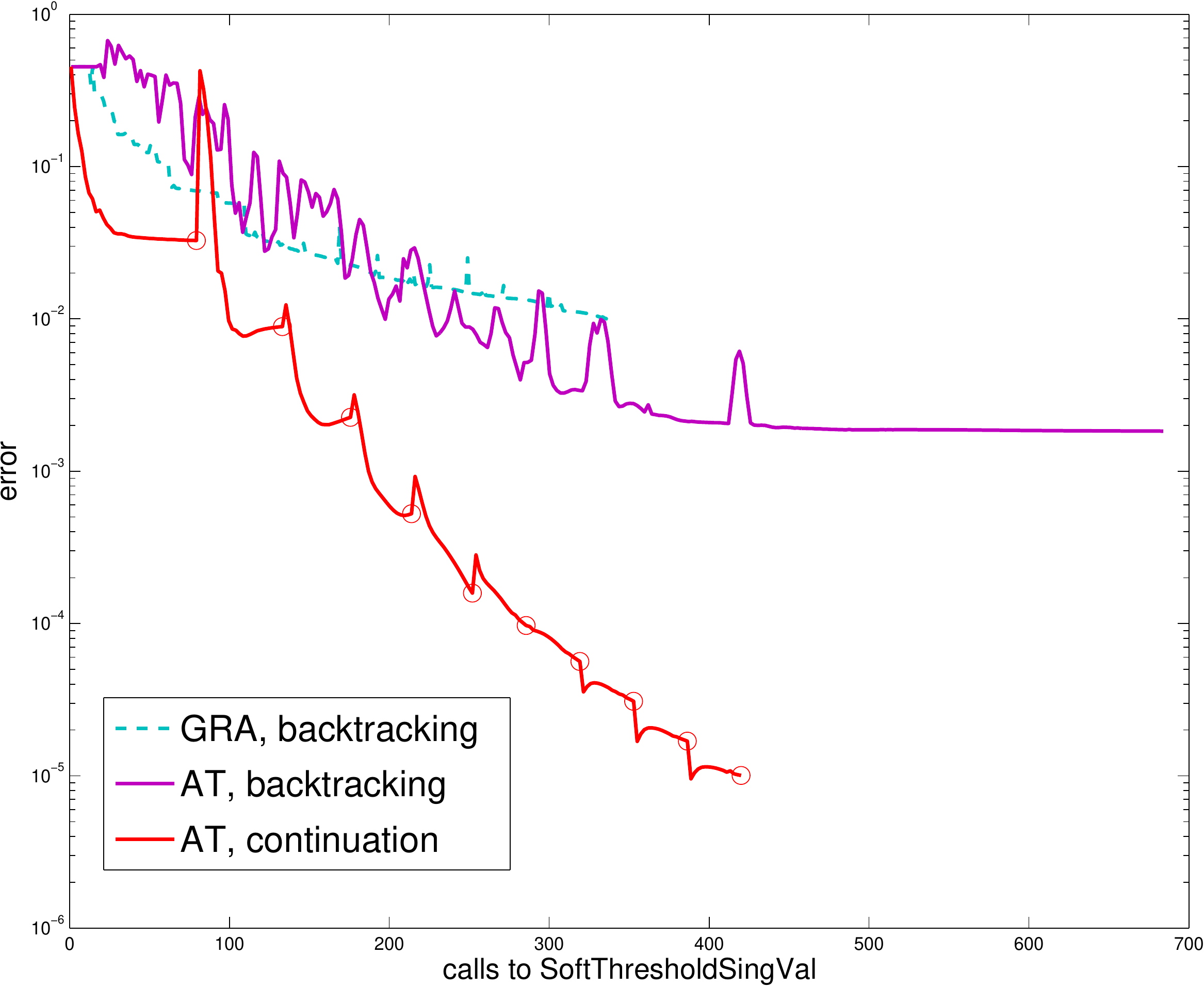}
    \caption{Noisy matrix completion on a non-low-rank matrix, using various first-order methods.}
    \label{fig:nuclear_hardcase}
\end{figure}
Figure \ref{fig:nuclear_hardcase} demonstrates these claims on a noisy matrix completion problem.
We constructed a $50\times 45$ matrix of rank 20, sampled $67\%$ of the entries,
and added white noise to yield a 30 dB SNR. By using this small problem size, we are able to use
CVX~\cite{cvx} to compute a reference solution; this took 5.4 minutes. The figure
depicts three different approaches to solving the problem: GRA and AT each with
$\mu=5\cdot 10^{-4}$ and no continuation; and AT with $\mu=10^{-2}$ and
accelerated continuation. The restart approach no longer has a beneficial effect, so it
is not shown. Each solver was run for 500 iterations, taking between $2.5$ and
$4$ seconds, and the plot depicts relative error versus the number of singular
value decompositions. As we can see, the advantage of GRA has been lost, and
continuation provides significant improvement: to achieve a relative
error of $10^{-2}$, the continuation method required only 
100 SVDs.

\subsection{\texorpdfstring{$\ell_1$-analysis}{l1-analysis}} 
\label{sec:analysis}

We present a brief example of solving the $\ell_1$-analysis problem
\eqref{eq:l1-analysis}, which is used in sparse recovery when the
primal variable $x$ is not sparse itself but rather sparse or
compressible in some other domain. The domain may be well-known, such
as frequency space ($W$ is a DFT) or a wavelet basis ($W$ is a
wavelet transform), and for these cases $W$ is invertible and even
orthogonal and the problem may be solved using a method like SPGL1.
However, these bases are often chosen out of convenience and not
because they best represent the signal.  More appropriate choices may
be {\em overcomplete} dictionaries $W \in \R^{p \times n}$ with $p \gg
n$, such as the undecimated wavelet transform, or the multilevel Gabor
dictionary.  The experiment below uses the Gabor dictionary with
$p=28n$.

To solve the LASSO problem using a dictionary $W$, the two
common approaches are {\em analysis}:
\begin{equation} \label{eq:l1-analysis2} 
\begin{array}{ll}
  \text{minimize} & \|Wx\|_1 \\
  \text{subject to} & \|y - A x\|_2 \le \epsilon
\end{array}
\end{equation}
and {\em synthesis}:
\begin{equation} \label{eq:l1-synthesis}
\begin{array}{ll}
  \text{minimize} & \|\alpha\|_1 \\
  \text{subject to} & \|y - A W^\T\alpha\|_2 \le \epsilon, 
\end{array}
\end{equation}
with decision variable $\alpha$. Similarly, the Dantzig selector
approach would have the same objectives but constraints of the form
$\|A^\T(y-Ax)\|_\infty \le \delta$ (analysis) and
$\|A^\T(y-AW^\T\alpha)\|_\infty \le \delta$ (synthesis).  When $W$ is
not orthogonal, the two approaches are generally different in
non-trivial ways, and furthermore, the solutions to {\em synthesis}
may be overly sensitive to the data
\cite{EladPrior}.  % As a simple example (using MATLAB
% notation for matrices), let $W = [I;I]$, and $x = e_1$ where $e_1$ is
% the first unit vector.  By analyzing $x$, we have $Wx = [e_1;e_1]$
% which is $2$-sparse.  Via synthesis, we have $x = W^\T\alpha_1 =
% W^\T\alpha_2$ for $\alpha_1 = [e_1;0]$ and $\alpha_2 = [0;e_1]$, and
% even $x = W^\T( t\alpha_1 + (1-t)\alpha_2)$.  The synthesis
% coefficients \emph{can} be sparser, but not unique and less stable.

The differences between analysis and synthesis are not very well
understood at the moment for two reasons. The first is that the
analysis problem has not been studied theoretically.  An exception is
the very recent paper \cite{analysis_stanford} which provides the
first results for $\ell_1$-analysis. The second is that there are no
existing efficient first-order algorithms to solve the analysis
problem, with the exception of the recent NESTA \cite{NESTA}
and C-SALSA \cite{CSALSA} algorithms, which both
work on the LASSO version, and only when $AA^\T=I$.

\begin{figure}
\centering
%trim option's parameter order: left bottom right top
%\setlength\fboxsep{0pt}
%\fbox{
 %\includegraphics[trim = 14mm 2mm 20mm 8mm, clip, width=4.2in]{figs/RMPI_fig_July27_Dantzig}
 %\includegraphics[trim = 14mm 2mm 20mm 8mm, clip, width=4.2in]{figs/RMPI_fig_July23}
 \includegraphics[width=4.2in]{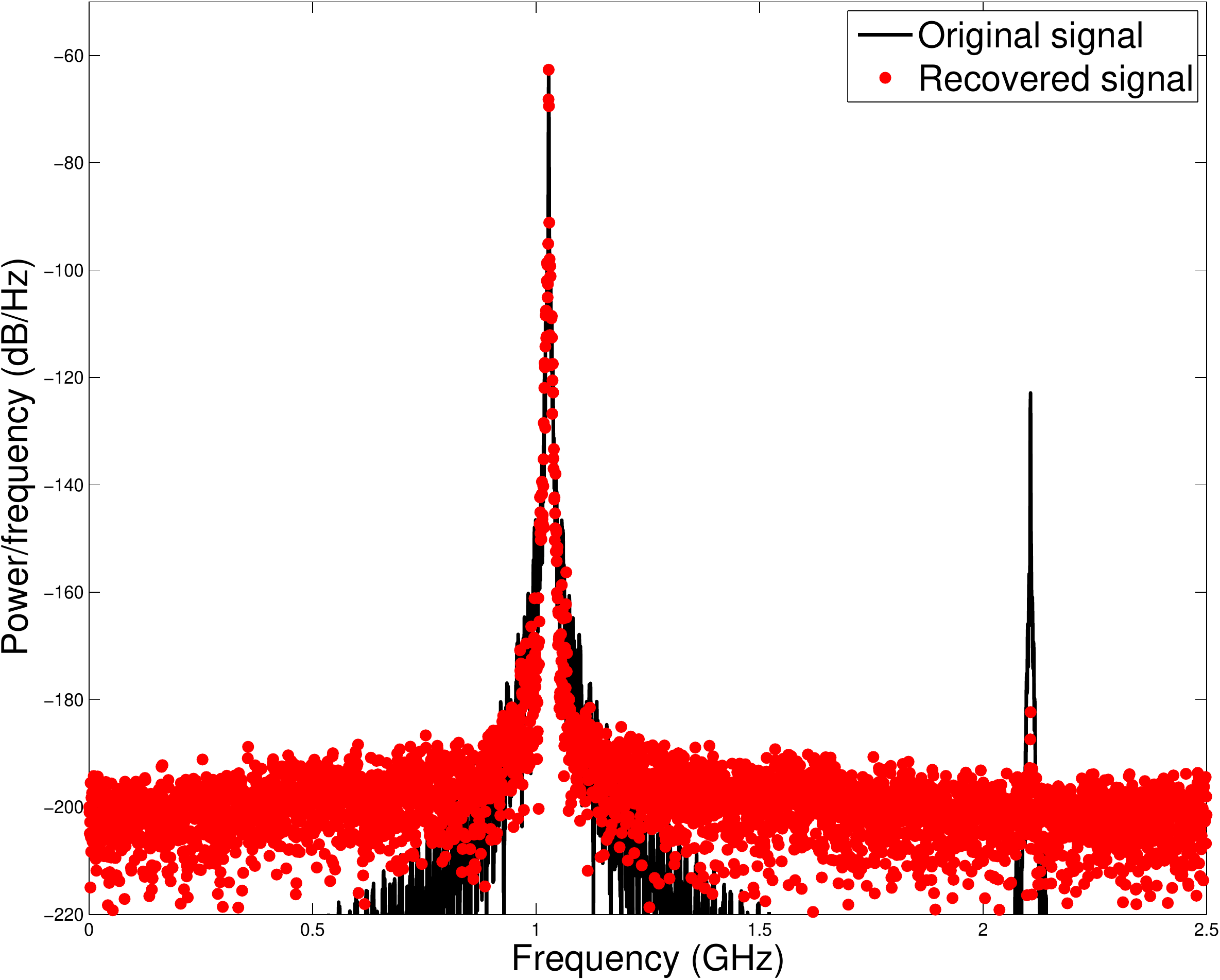}
 \includegraphics[width=4.2in]{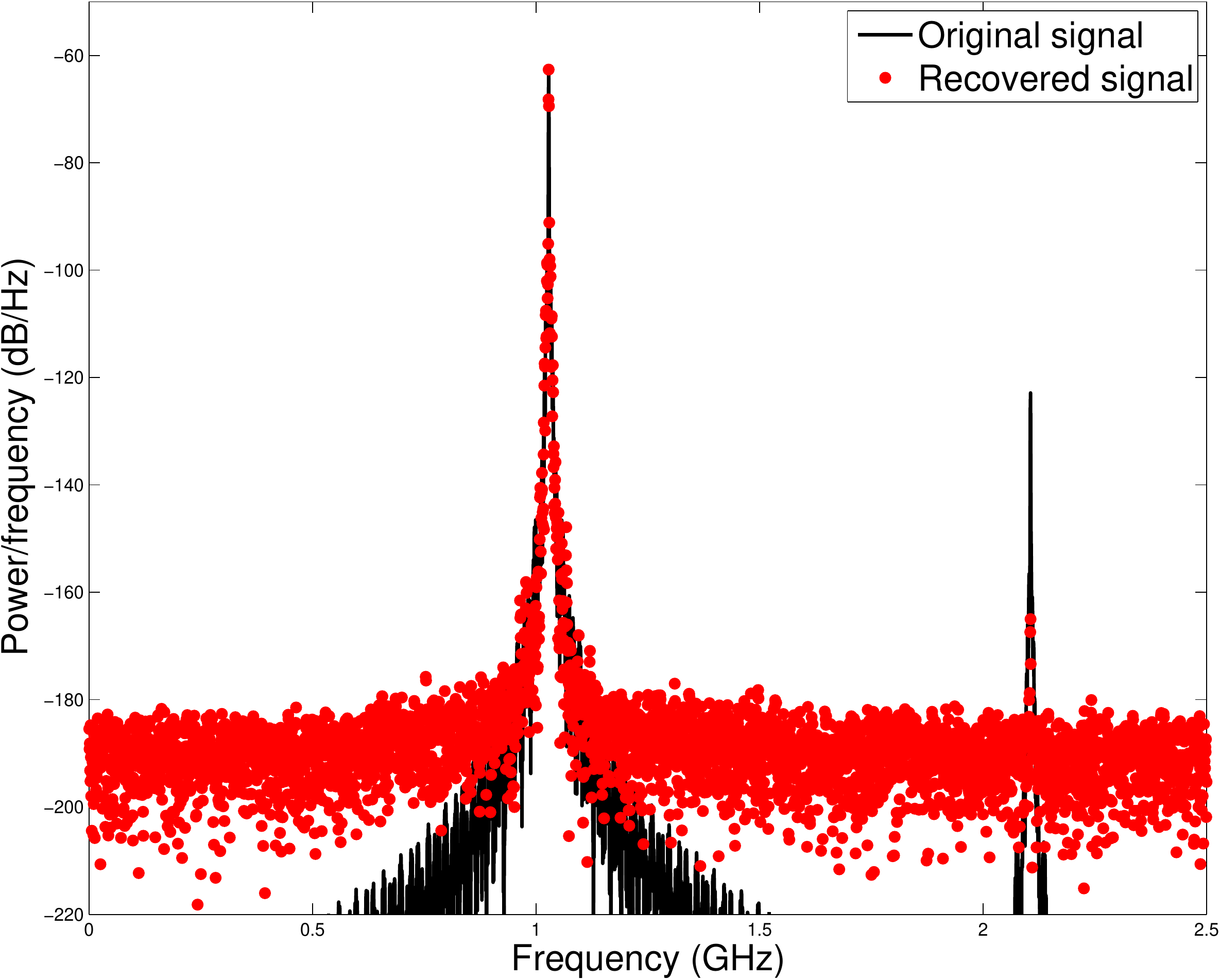}
    %}
%[-16pt]
%\\[-10pt]
 \caption{Recovery of a small pulse $60$ dB below a large pulse. The first plot
 employs a Dantzig selector, the bottom employs LASSO.}
\label{fig:RMPI}
\end{figure}

% Details: (see example_RMPI.m)
% Big Pulse: 
% Small Pulse:
% Both: frequency and phase uniformly random
% N = 8192 = 
% SNR of signal is 60.1 dB; SNR of big tone is 60.1 dB; SNR of small tone is -4.4 dB
% N = 8192, d = 228864, d/N = 27.9 for dictionary

With the dual conic approach, it is now possible to solve the smoothed
analysis problem, and by using continuation techniques, the effect of
the smoothing is negligible. To illustrate, we constructed a realistic
test of the recovery of two radio-frequency
radar pulses that overlap in time, as depicted in Figure~\ref{fig:RMPI}.
The first pulse is large, while
the second pulse has 60 dB smaller amplitude.  Both have carrier
frequencies and phases that are chosen uniformly at random, and noise
is added so that the small pulse has a signal-to-noise ratio of $0.1$ dB.

The signal is recovered at Nyquist rate resolution for $2.5$ GHz
bandwidth, and the time period is a little over $1600$ ns, so that $n =
8,\!192$.  The sensing matrix $A$ is modeled as a block-diagonal
matrix with $\pm 1$ entries on the blocks, representing a system that
randomly mixes and integrates the input signal, taking 8 measurements
every $100$ ns, which is $12.5 \times$ below the Nyquist rate.  Thus
$A$ is a $648 \times 8,\!192$ matrix and $W$ is a $228,\!864 \times
8,\!192$ matrix.  We applied both the Dantzig selector and LASSO
models to this problem.

To solve the problems, we employed the AT variant with accelerated
continuation. At
each iteration, the stopping tolerance is decreased by a factor of
$1.5$, and the continuation loop is ended when it no longer produces
significantly different answers, which is usually between 2 and 8
iterations.  The value of $\mu$ is set to 
\begin{equation}
  \mu = 0.1 \frac{ \|W x_{LS}\|_2 }{ \frac{1}{2}\|x_{LS}\|^2 }, 
\end{equation}
where $x_{LS}$ is the least-squares solution to $Ax=y$, which is easy
to calculate since $ m \ll n$ (and independent of $p$).

To enhance the results, we employed an outer \emph{reweighting} loop
\cite{Candes_reweighting}, in which we replace the $\|Wx\|_1$ term
with $\|RWx\|_1$ for some diagonal weight matrix $R$. Each reweighting
involves a full solve of either the Dantzig selector or the LASSO.  This
loop is run until convergence, which is typically about 2 to 5
iterations.  The results are plotted in the frequency domain in Figure
\ref{fig:RMPI}.  The large pulse is the dominant spike at about $1.1$
GHz, and it is easily recovered to very high precision (the relative
$\ell_2$ error is less than $5\cdot 10^{-3}$).  The small pulse is at
about $2.2$ GHz, and because it has an SNR of $0.1$ dB and the
measurements are undersampled by a factor of $12$, it is not possible
to recover it exactly, but we can still detect its presence and
accurately estimate its carrier frequency.

Table \ref{table:RMPI} reports the computational results of the test,
first solving either the Dantzig selector or the LASSO, and then taking
a single reweighting step and re-solving.
Each call of the algorithm takes about 100 iterations, and the
Dantzig selector or the LASSO is solved in about 1 minute, which
%algorithm is called 13 to 15 times due to the continuation and
%reweighting steps.  The answer is computed in about 2 minutes, which
is impressive since, due to the extremely large size of $W$, this
problem is intractable using an interior-point method.

\begin{table}
\begin{tabular}{|*{9}{l|}}
    \hline
    &
  \multicolumn{4}{|c|}{Dantzig selector} &
  \multicolumn{4}{|c|}{LASSO} \\
  \hline
  Reweighting step & \# cont. & iter. & time & error &
     \# cont. & iter. & time & error \\
  \hline
0&7 & 821 & 89.2 s & $2.7\cdot 10^{-3}$ &
  6 & 569 & 54.8 s & $2.2\cdot 10^{-3}$ \\
  \hline
  1&8 & 1021 & 114.6 s & $2.1\cdot 10^{-3}$ &
    7 &  683 & 67.8  s & $1.8\cdot 10^{-3}$ \\
    %\hline
    %\hline
%total &15 &  1842 & 203.8 s &  $2.1 10^{-3}$  &
 %13 & 1252 & 122.6 s &  $1.8 10^{-3}$ \\
  \hline
\end{tabular}
\caption{Details on the simulation used in Figure \ref{fig:RMPI}.
For both the Dantzig selector and LASSO versions, the algorithm
was run with continuation and reweighting.  The ``cont.''~column
is the number of continuation steps, the ``iter.''~column is the number of iterations (over all
the continuation steps), the ``time''~column is in seconds, and the ``error''~column
is the relative $\ell_2$ error.  Each row is one solve of the Dantzig selector or the LASSO.}
\label{table:RMPI}
\end{table}

\section{Software: TFOCS}
\label{sec:software}

The work described in this paper has been incorporated into a software
package, Templates for First-Order Conic Solvers (TFOCS, pronounced
\emph{tee-fox}), which will be made publicly available at
\url{http://tfocs.stanford.edu}.  As its name implies, this package is
a set of templates, or building blocks, that can be used
to construct efficient, customized solvers for a variety of models.

To illustrate the usage of the software, let us show
how to construct a simple solver for the smoothed Dantzig selector
model described in \S\ref{sec:dantzig}. We will begin by
assuming  that we are given the problem data
\verb@A@, \verb@b@, \verb@delta@, and a fixed smoothing
parameter \verb@mu@.
The basic solver templates require two functions to 
complete their work. The first function computes
the value and gradient of $\gs$,
the smooth component of the composite dual:
\begin{code}
function [ val, grad, x ] = g_dantzig( A, y, x0, mu, z )
x    = SoftThreshold( x0 - (1/mu) * A' * ( A * z ), 1/mu );
grad = A' * ( y - A * x );
val  = z' * grad - norm( x, 1 ) - 0.5 * mu * norm( x - x0 ) .^ 2;
\end{code}
The second function computes the generalized projection associated with
the nonsmooth component of the dual, which is in this case $h(z)=\delta\|z\|_1$.
\begin{code}
function [ z, h ] = h_dantzig( delta, z_old, grad, L )
z = SoftThreshold( z_old - (1/L)*grad, delta/L );
h = delta * norm( z, 1 );
\end{code}
Both of these functions depend on the soft-thresholding operator:
\begin{code}
function y = SoftThreshold( x, t )
y = sign( x ) .* max( abs( x ) - t, 0 );
\end{code}
Armed with these functions, the following code solves the smoothed
Dantzig selector, using the Auslender/Teboulle first-order variant
and the default choices for line search and stopping criteria.
\begin{code}
function x = Dantzig_smoothed( A, y, delta, mu )
[m,n] = size(A);
x0    = zeros(n,1); z0 = zeros(n,1);
g_sm  = @(z) g_dantzig( A, y, x0, mu, z );
h     = @(z,g,L) h_dantzig( delta, z, g, L );
[z,x] = solver_AT( g_sm, h, z0 );
\end{code}
The last line of code calls the AT solver, which upon completion
delivers its solution to both the smoothed
dual and the primal problem.

This simple solver is likely to be unnecessarily inefficient if 
\verb@A@ exhibits any sort of fast operator structure, but this
can be remedied rather simply. Note that the solver itself needs
to have no knowledge of \verb@A@ or \verb@y@ above; its interaction
with these quantities comes only through calls to \verb@g_dantzig@.
Therefore, we are free to \emph{rewrite} \verb@g_dantzig@ in a more numerically
efficient manner: for instance, if \verb@A@ is derived from a Fourier
transforms, we may substitute fast Fourier transform operations for
matrix-vector multiplications. This simple change will reduce
the cost of each iteration % will be reduced   % Ewout caught this
from $\order(mn)$ to $\order(n\log n)$, and the storage requirements from
$\order(n^2)$ to $\order(n)$. For large-scale problems these savings can
be quite significant.

A further improvement in performance is possible through careful
management of the linear operator calculations as described in \S\ref{sec:linop}.  
The TFOCS solver templates can perform this
management automatically.  To take advantage of it, we replace \verb@g_dantzig@
with two functions: one which implements the underlying linear 
operator, and one which implements the remainder of the smooth
function. The details of how to accomplish this are best left
to the user guide \cite{Templates}.

\iffalse
smooth function without
the linear operations,
\begin{code}
function [ f, x ] = g2_dantzig( x0, mu, ATAz )
x = SoftThreshold( x0 + (1/mu) * ATz, 1/mu );
f = ATz' * x - norm( x, 1 ) - 0.5 * mu * norm( x - x0 ) .^ 2 ;
\end{code}
and one which implements nothing \emph{but} the underlying affine operation.
It should obey the following calling sequence, but of course be implemented
efficiently:
\begin{code}
function z = affine_dantzig( x, mode )
global A y
if nargin == 0,
    z = size(A,2)*[1,1];     % Return the operator size
elseif nargin > 1 && isequal( mode, 'adjoint' ),
    z = - A' * ( A * x );    % Perform the adjoint operation
else
    z = A' * ( y - A * x );  % Perform the forward operation
end    
\end{code}
Given these tools, the modified Dantzig solver then becomes
\begin{code}
function x = Dantzig_smoothed( A_op, delta, mu )
sz    = A_op(); n = sz(1);
x0    = zeros(n,1); z0 = zeros(n,1);
g2_sm = @(z) g2_dantzig( x0, mu, z );
h     = @(z,g,L) h_dantzig( delta, z, g, L );
[z,x] = solver_AT( { g2_sm, A_op, 'adjoint' }, h, z0 );
\end{code}
where \verb@A_op@ is a handle to the \verb@affine_dantzig@ function. Passing
the quantities \verb@g2_sm@, \verb@A_op@, and the string \verb@'adjoint'@ in
a cell array communicates to the solver the information it needs to take
full advantage of the linear operator structure.
\fi

Of course, as part of the complete package, we have supplied a
solver for the Dantzig selector that exploits each of
these efficiencies, and others. Similar solvers have been
created for the LASSO, TV, and other models discussed here.
But the flexibility of the lower-level templates will allow
these same efficiencies to be achieved for many models that
we have not discussed here. In fact, evidently the templates
are not restricted to our specific conic form (\ref{eq:stdform})
or its dual, and we hope the software finds application
outside of the compressed sensing domain as well.

With the software release is a detailed user guide \cite{Templates}
that covers the usage, and documents the most
popular problem formulations and provide some examples.
We refer the reader to the user guide for further software details.

\section{Conclusion}

% \TODO{A listing of points to mention:}
% \begin{itemize}
%     \item Continuation: work out a theory that allows for inexact solves (and so gives recommendations for stopping tolerances)
%     \item Choose $\mu$ in a better fashion
%      \item investigate why the ``restart'' method works so well for
%      some problems.  It seems to work well only when we're close to a
%      solution and the answer is very sparse.  In this setting, the
%      assumptions of the PARNES paper may be valid.  We have a different
%      formulation than them, but it may be that we have ``local'' strong
%      convexity whenever we're close to an answer that is very sparse or
%      very low-rank. 
%     \item Scaled norms
%         As in \S\ref{sec:analysis_tv} and \S\ref{sec:numerics_tv_wavelet}.

% \end{itemize}

We have developed a convenient framework for constructing first-order
methods, which is flexible and handles a variety of convex cone
problems, including problems that did not have efficient  % taking Ewout's changes
algorithms before. On the implementation side, we have
introduced ideas which lead to novel, stable, and efficient
algorithms. When comparing our implementation on specific problems
such as the LASSO, which have been intensively studied, our techniques
appear surprisingly competitive with the state of the art.

The templates from this paper are flexible in a manner that has not
yet been seen in sparse recovery software.  Not only are the solvers
interchangeable and hence easy to benchmark one algorithm against
another, but they work with a wide variety of formulations which we
hope will greatly enable other researches.  It is also our goal that
the software will be easy to use for non-experts, and to this end our
future research will be to improve the usability of the software and
to provide sensible default parameters.  One major topic to pursue is
choosing the smoothing parameter $\mu$.  Further efforts will also be
made to improve the line search to take advantage of strong convexity,
to better manage expensive projections,
and to use scaled norms so that dual variables are all on the same
scale (this is an issue only for objectives with several additive
terms, as in the TV with analysis problem in \S\ref{sec:analysis_tv}
and \S\ref{sec:numerics_tv_wavelet}).

A subsequent paper will cover these issues, as well as further
investigating local strong convexity and how to take advantage of this
in an optimal manner, and improving the accelerated continuation
scheme by taking into account the inexact solves.  The software and
user guide \cite{Templates} will be kept up-to-date and supported.

\appendix
\section{\texorpdfstring{Exact Penalty}{Appendix: Exact Penalty}}
\label{app:exact}

The general approach proposed in this paper is to add a strongly
convex term to the primal objective function in order to smooth the
dual objective.  We now draw a parallel with augmented Lagrangian and
penalty function methods that eliminate constraints by incorporating a
penalty term $p$ into the objective, \eg, $\min_{x:Ax=b} f(x)$
becomes $\min_{x} f(x) + \lambda p(Ax-b ) $.  For reasonable choices
of $p$, the two problems become equivalent as $\lambda \rightarrow
\infty$.  Remarkably, for some special choices of $p$, the two
problems are equivalent for a large but \emph{finite} value of
$\lambda$.  Usually, $p$ is a non-differentiable function, such as $p
= \|\cdot\|_1$ or $p = \|\cdot\|_2$ (not $p=\|\cdot\|_2^2$); see
Bertsekas' book~\cite{BertsekasBook} for a discussion.

Our approach uses a different type of perturbation, since our goal is
not to eliminate constraints but rather to smooth the objective. But
we use the term ``exact penalty'' because, for some problems, we have
a similar result: the smoothed and unsmoothed problems are (nearly)
equivalent for some $\mu > 0$.  The result below also departs from
traditional exact penalty results because our perturbation
$\thalf\|x-x_0\|_2^2$ \emph{is} smooth.

Below we present the proof of Theorem~\ref{thm:exactRecovery}, which
exploits the polyhedral constraint set of linear programs.

\paragraph{Proof of Theorem \ref{thm:exactRecovery}.} 
We consider the linear program
\begin{equation}
\begin{array}{ll}
	\text{minimize} & \<c, x\> \\
	\text{subject to} & x \in \cP, 
\end{array}	
\tag{LP}
\label{eq:lp}
\end{equation}
where $\cP$ is a convex polyhedron, and its perturbed version
\begin{equation}
\begin{array}{ll}
	\text{minimize} & \<c, x\> +  \thalf\mu \<x-x_0, x-x_0\>_Q\\
	\text{subject to} & x \in \cP,   
\end{array}	
\tag{QP}
\label{eq:qp}
\end{equation}
where $\<x,y\> \triangleq \<x, Qy\>$ for some positive semidefinite $Q$.
Let $E^\star$ be the solution
set of \eqref{eq:lp} ($E^\star$ is a vertex or a face of the feasible
polyhedron) and let $x^\star$ be any point in $E^\star$ such that
$\<x-x_0, x-x_0\>_Q$ is minimum. (Note that when $Q$ is the identity,
$x^\star$ is the usual projection onto the convex set $E^\star$.)
With $f_\mu(x) \triangleq \<c, x\> + \thalf\mu \<x-x_0,x-x_0\>_Q$, it
is of course sufficient to show that $x^\star$ is a local minimizer of
$f_\mu$ to prove the theorem. Put differently, let $x \in \cP$ and
consider $x^\star + t(x-x^\star)$ where $0 \le t \le 1$. Then it suffices to
show that $\lim_{t \goto 0^+} f_\mu(x^\star + t(x-x^\star)) \ge
f_\mu(x^\star)$. We now compute
\begin{align*}
  f_\mu(x^\star + t(x-x^\star)) & = \<c, x^\star + t(x-x^\star)\> +
  \thalf\mu \<x^\star + t(x-x^\star) - x_0,x^\star + t(x-x^\star) - x_0\>_Q\\
  & = f_\mu(x^\star) + t\<c,x-x^\star\> + \mu t \<x - x^\star, x^\star
  - x_0\>_Q + \thalf\mu t^2 \<x -
  x^\star,x-x^\star\>_Q. 
  \end{align*}
Therefore, it suffices to establish that for all $x \in \cP$, 
\[
\<c,x-x^\star\> + \mu \<x - x^\star, x^\star - x_0\>_Q \ge 0,
\]
provided $\mu$ is sufficiently small. Now $x \in \cP$ can be expressed
as a convex combination of its extreme points (vertices) plus a
nonnegative combination of its extreme directions (in case $\cP$ is
unbounded). Thus, if $\{v_i\}$ and $\{d_j\}$ are the finite families
of extreme points and directions, then 
\[
x = \sum_{i} \lambda_i v_i + \sum_j \rho_j d_j, 
\]
with $\lambda_i \ge 0$, $\sum_i \lambda_i = 1$, $\rho_j \ge 0$. The
extreme directions---if they exist---obey $\<c,d_j\> \ge 0$ as
otherwise, the optimal value of our LP would be $-\infty$. Let $I$
denote those vertices which are {\em not} solutions to \eqref{eq:lp},
\ie, such that $\<c, v_i - \hat x\> > 0$ for any $\hat x \in
E^\star$.  Likewise, let $J$ be the set of extreme directions obeying
$\<c, d_j\> > 0$. With this, we decompose $x - x^\star$ as
\[
x-x^\star = \Bigl(\sum_{i \notin I} \lambda_i (v_i - x^\star) +
\sum_{j \notin J} \rho_j d_j\Bigr) + \Bigl(\sum_{i \in I} \lambda_i (v_i -
x^\star) + \sum_{j \in J} \rho_j d_j\Bigr). 
\]
It is not hard to see that this decomposition is of the form 
\[
x - x^\star = \alpha (\hat x - x^\star) +  \Bigl(\sum_{i \in I} \lambda_i (v_i -
x^\star) + \sum_{j \in J} \rho_j d_j\Bigr), 
\]
where $\hat x \in E^\star$ and $\alpha$ is a nonnegative scalar. This
gives
\[
\<c,x-x^\star\> = \sum_{i \in I} \lambda_i \<c, v_i - x^\star\> +
\sum_{j \in J} \rho_j \<c, d_j\>.
\]
Let $\alpha_i \ge 0$ (resp.~$\beta_j \ge 0$) be the cosine of the
angle between $c$ and $v_i - x^\star$ (resp.~$c$ and $d_j$).
Then for $i \in I$ and $j \in J$, we have $\alpha_i > 0$ and $\beta_j > 0$.
We can write
\[
\<c,x-x^\star\> = \sum_{i \in I} \lambda_i \alpha_i \|c\|_2 \|v_i -
x^\star\|_2 + \sum_{j \in J} \rho_j \beta_j \|c\|_2 \|d_j\|_2, 
\]
\newcommand{\iprodMed}[2]{\Bigl\langle #1 , #2 \Bigr\rangle} which is
strictly positive.  Furthermore,
\begin{align*}
  \iprod{x - x^\star}{ x^\star - x_0}_Q & = \alpha \iprod{\hat x -
    x^\star}{ x^\star - x_0}_Q + \iprodMed{\sum_{i \in I} \lambda_i (v_i
    - x^\star) + \sum_{j \in J}
    \rho_j d_j}{ x^\star - x_0}_Q\\
  & \ge \iprodMed{\sum_{i \in I} \lambda_i (v_i - x^\star) + \sum_{j
      \in J}
    \rho_j d_j}{ x^\star - x_0}_Q\\
  & \ge - \|Q(x^\star - x_0)\|_2 \Bigl( \sum_{i \in I} \lambda_i \|v_i-
  x^\star\|_2 + \sum_{j \in J} \rho_j \|d_j\|_2\Bigr).
\end{align*}
The second inequality holds because $x^\star$ minimizes $\<\hat x -
x_0, \hat x - x_0\>_Q$ over $E^\star$, which implies $\<\hat x -
x^\star, x^\star - x_0\>_Q \ge 0$. The third is a consequence of the
Cauchy-Schwartz inequality. In conclusion, with $\mu_Q \triangleq \mu
\|Q(x^\star-x_0)\|_2$, we have 
\[
  \<c,x-x^\star\> + \mu \<x - x^\star, x^\star - x_0\>_Q \ge \sum_{i
    \in I} (\alpha_i \|c\|_2 - \mu_Q) \lambda_i \|v_i -
  x^\star\|_2 \\
  + \sum_{j \in J} (\beta_j \|c\|_2 - \mu_Q) \rho_j \|d_j\|_2,
\]
and it is clear that since $\min_i \alpha_i > 0$ and $\min_j \beta_j >
0$, selecting $\mu$ small enough guarantees that the right-hand side
is nonnegative.

% \paragraph{Remark A.2 (Dantzig selector).}  If the solution to
% $$\min \|x\|_1 \quad\st\quad \|A^\T(b-Ax)\|_\infty \le \eps$$ 
% is unique, then it is the solution to 
% $$ \min_x \|x\|_1 + \frac{\mu}{2} \|x-x_0\|_2^2 \quad\st\quad \|A^\T(b-Ax)\|_\infty \le \eps $$
% provided $\mu$ is sufficiently small (simply transform the Dantzig
% selector to an equivalent LP in standard form, then choose $P$ to be a
% diagonal matrix that selects only the $x$ variable).  

% \paragraph{Remark A.3 (Noiseless basis pursuit)\cite{YinPenalty}.}  
% If the solution to
% $$\min \|x\|_1 \quad\st\quad Ax=b $$
% is unique, then it is the solution to
% $$ \min_x \|x\|_1 + \frac{\mu}{2} \|x-x_0\|_2^2 \quad\st\quad Ax = b. $$
% This is not true for noise-aware basis pursuit (\ie, constraints of
% the form $\|Ax-b\|_2 \le \eps$ constraints) since it is no longer an
% LP and the constraint set is not polyhedral.
% Futhermore, it is possible to construct a counter-example in the noise-aware case,
% so the result is sharp.
\
\section{\texorpdfstring{Creating a Synthetic Test Problem}{Appendix:
Creating a Synthetic Test Problem}} \label{sec:testProblem}

It is desirable to use test problems that have a precisely known
solution.  In some cases, such as compressed sensing problems in the
absence of noise, the solution may be known, but in general this is not
true. A common practice is to solve problems with an interior-point
method (IPM) solver, since IPM software is mature and accurate. 
However, IPMs do not scale well with the problem size, and cannot take
advantage of fast algorithms to compute matrix-vector products.  Another
disadvantage is that the output of an IPM is in the interior of the
feasible set, which in most cases means that it is not sparse.

Below, we outline procedures that show how to generate problems with
known exact solutions (to machine precision) for several common
problems.  The numerical experiments earlier in this paper used this
method to generate the test problems.  This is inspired by
\cite{Nesterov07}, but we use a variation that gives much more control
over the properties of the problem and the solution.

\paragraph{Basis pursuit.} Consider the basis pursuit problem and its
dual \begin{equation*} \begin{array}{ll} \text{minimize} & \|x\|_1 \\
\text{subject to} & A x = \bb, \end{array}	\quad \begin{array}{ll}
\text{maximize}   & \<\bb, \lambda\> \\ \text{subject to} & \|A^\T
\lambda\|_\infty \le 1. \end{array} \end{equation*} At the optimal
primal and dual solutions $x^\star$ and $\lambda^\star$, the KKT
conditions hold: \begin{align*} Ax^\star = \bb \quad \|A^\T
\lambda^\star\|_\infty \le 1 \quad (A^\T \lambda^\star )_T =
\text{sign}( x^\star_T ), \end{align*} where $T = \supp{x^\star}$.

To generate the exact solution, the first step is to choose $A$ and
$\bb$.  For example, after choosing $A$, $\bb$ may be chosen as $\bb =
A\tilde{x}$ for some $\tilde{x}$ that has interesting properties
(\eg, $\tilde{x}$ is $s$-sparse or is the wavelet coefficient sequence
of an image).  Then any primal dual solver is run to high accuracy to
generate $x^\star$ and $\lambda^\star$.  These solutions are usually
accurate, but not quite accurate to machine precision.

The idea is that $x^\star$ and $\lambda^\star$ are exact solutions to a
slightly perturbed problem.  Define $T = \supp{\hat{x}}$; in sparse
recovery problems, or for any linear programming problem, we have $|T|
\le m$ where $m$ is the length of the data vector $\bb$. The matrix $A$
is modified slightly by defining $\tilde{A} \leftarrow A D$ where $D$ is
a diagonal matrix.  $D$ is calculated to ensure that $\| (D A^\T
\lambda^\star)_{T^c} \|_\infty < 1 $ and $(DA^\T \lambda^\star)_T =
\text{sign}(x^\star) $.  If the original primal dual solver was run to
high accuracy, $D$ is very close to the identity. In practice, we
observe that the diagonal entries of $D$ are usually within $.01$ of
$1$.

The primal variable $x^\star$ is cleaned by solving $\tilde{A_T}
x^\star_T = \bb$;  this is unique, assuming $A_T$ has full column rank.
If the original problem was solved to high accuracy (in which case $D$
will have all positive entries), then the cleaning-up procedure does not
affect the sign of $x^\star$. The vectors $x^\star$ and $\lambda^\star$
are now optimal solutions to the basis pursuit problem using $\tilde{A}$
and $\bb$.

\paragraph{The LASSO.} A similar procedure is carried out
for the LASSO and its dual given by
\begin{equation*} \begin{array}{ll} \text{minimize} & \|x\|_1 \\
\text{subject to} & \|A x - \bb\|_2 \le \epsilon, \end{array}	\quad
\begin{array}{ll} \text{maximize}   & \<\bb, \lambda\> - \epsilon
\|\lambda\|_2 \\ \text{subject to} & \|A^\T \lambda\|_\infty \le 1.
\end{array} \end{equation*} 
Let $z$ be the variable such that $\bb =
Ax+z$, then strong duality holds if $$ \left( \<Ax^\star,
\lambda^\star\> - \|x^\star\|_1 \right) + \left( \<z, \lambda^\star\> -
\epsilon \|\lambda\|_2 \right) = 0. $$ The operator $\tilde{A}$ is
chosen as before, so $ \<Ax^\star, \lambda^\star\> - \|x^\star\|_1 =0$. 
Thus $z$ needs to satisfy $ \<z, \lambda^\star\> - \epsilon
\|\lambda\|_2 $, \ie, $z = \epsilon y/\|y\|_2 $.  This means that
$x^\star$ and $\lambda^\star$ are optimal solutions to the LASSO
with $\tilde{A}$ and $\tilde{\bb} = Ax^\star + z $.

\paragraph{Other problems.} For basis pursuit and the LASSO, it is possible to obtain exact solutions to the
\emph{smoothed} problem (with $d(x) = \frac{1}{2}\|x-x_0\|_2^2$). For
the Dantzig selector, an exact solution for the smoothed problem can
also be obtained in a similar fashion.  To find an exact solution to
the unsmoothed problem, we take advantage of the exact penalty
property from \S\ref{sec:exactRecovery} and simply find the smoothed
solution for a sequence of problems to get a good estimate of $x_0$
and then solve for a very small value of $\mu$.

\small

\subsection*{Acknowledgements}
This work has been partially supported by ONR grants N00014-09-1-0469
and N00014-08-1-0749, by a DARPA grant FA8650-08-C-7853, and by the
2006 Waterman Award from NSF. 
We would like to thank E. van den Berg for a careful reading of the manuscript.
SRB would like to thank Peter Stobbe
for his Hadamard Transform and Gabor dictionary code.

\bibliographystyle{plain}
%  for arxiv, using this style:
%\bibliographystyle{hplain}
\bibliography{cnest}

\end{document}